\documentclass[12pt,leqno]{amsart}

\usepackage{amsthm}
\usepackage[usenames]{color}

\usepackage{amsmath}
\usepackage{amscd}
\usepackage{amsmath}
\usepackage{amssymb}
\usepackage{latexsym}
\makeatletter

\@addtoreset{equation}{section}
\makeatother

\newtheorem{thm}{Theorem}[subsection]
\newtheorem{pr}[thm]{Proposition}
\newtheorem{df}[thm]{Definition}
\newtheorem{lm}[thm]{Lemma}
\newtheorem{cor}[thm]{Corollary}

\newtheorem{ex}[thm]{Example}

\newcommand{\sm}{\raisebox{2.33pt}{~\rule{6.4pt}{1.3pt}~}}

\input xy \xyoption{all} \CompileMatrices

\begin{document}

\title{Graded quotients of ramification groups\\
of local fields with imperfect residue fields}
\author{Takeshi Saito}

\maketitle
\begin{abstract}
We prove that the graded quotients
of the filtration by ramification groups
of any henselian discrete valuation field
of residue characteristic $p>0$
are ${\mathbf F}_p$-vector spaces.
We define an injection of
the character group of each
graded quotient to
a twisted cotangent space
defined using
a cotangent complex.
\end{abstract}

Let $K$ be a henselian discrete
valuation field of residue characteristic
$p>0$
and $L$ be a finite Galois extension
of $K$ of Galois group $G$.
Then, the decreasing filtration on $G$
by upper numbering
ramification groups $G^r
\subset G$
is defined in \cite{AS}.
The definition is rephrased purely
in the language of schemes in \cite{TJM}
and is briefly recalled
in Sections \ref{ssFun} and \ref{ssRam}.
The ramification groups $G^r$
are normal subgroups
indexed by positive rational numbers
$r>0$.
We set $G^{r+}=\bigcup_{s>r}G^s
\subset G^r$
and call ${\rm Gr}^rG=
G^r/G^{r+}$ the graded quotients.
The inertia subgroup $I\subset G$
and its $p$-Sylow subgroup
$P\subset I$ equal
$G^{1+}\subset G^1$.

\medskip
\noindent
{\bf Theorem \ref{thmgr}.}
{\em For $r>1$,
the graded quotients ${\rm Gr}^rG=
G^r/G^{r+}$ are ${\mathbf F}_p$-vector
spaces.}

\medskip
In the classical case where
the residue field is perfect,
Theorem \ref{thmgr} is proved in 
\cite[Corollaire 1 of
Proposition 7, Section 2, Chapitre IV]{CL}
by using the lower numbering filtration.
Many other cases have been proved
in \cite{X} by using $p$-adic differential
equations and in \cite{nonlog}
by geometric methods.

We give two proofs of
Theorem \ref{thmgr}.
The first proof is 
by the reduction
to the classical case.
A similar approach was proposed
in \cite{Bo}.
Using the first proof,
we give in \cite{TATA}
a characterization 
of the filtration by 
ramification groups
by a certain functoriality
and the classical case
where the residue field is
perfect.
We also deduce the integrality
of the conductor
of a representation
of $G$ from the classical case.

The second proof
is more geometric
and by the reduction to the case
where $r$ is an integer.
In the notation prepared below,
the second proof
gives rise to an injection
\begin{equation*}
{\rm ch}\colon
{\rm Hom}({\rm Gr}^rG,
{\mathbf F}_p)
\to 
{\rm Hom}_{\bar F}
({\mathfrak m}_{\bar K}^r
/{\mathfrak m}_{\bar K}^{r+},
H_1(L_{\bar F/S}))
\leqno{\rm (\ref{eqgrK})}
\end{equation*}
for every rational number $r>1$,
called the characteristic form.
The morphism (\ref{eqgrK}) is a generalization
of a non-logarithmic variant
of the refined Swan conductor
defined by Kato \cite[Corollary (5.2)]{Kato}
in the case where $G$ is abelian.
In a next article,
we plan to construct the characteristic cycle
of a constructible sheaf 
in mixed characteristic case
in the framework of \cite{LEX},
at least on a dense open subset
of the boundary
as an application of the injection (\ref{eqgrK}).

The two proofs are based
on the following idea.
For a discrete valuation field $K$,
suppose that we find an extension $K'$
such that 
the induced morphisms
on the graded quotients
of ramification groups be isomorphisms.
Then, 
the conclusion of Theorem \ref{thmgr}
for $K'$ would imply the same conclusion
for $K$.
In the first proof,
we find $K'$ with
perfect residue field
where the conclusion is a classical result.
In the second proof,
we find $K'$ such that
the index $r'$
corresponding to
a given rational number $r>1$ is an integer.
For integral indices,
we prove that
a geometric construction similar
to the equal characteristic case in \cite{nonlog} 
yields Theorem \ref{thmgr}
and an injection 
{\rm (\ref{eqgrK})}.
Finding such $K'$
is thus a crucial reduction step
in the proof of Theorem \ref{thmgr}.

To establish the reduction step,
we need to go back to a geometric
construction behind the definition
of ramification groups recalled in
Section \ref{ssFun}.
The construction
provides a certain finite Galois
covering of a vector space $\Theta^{(r)}$
over an algebraic closure
of the residue field
such that the Galois group
is the graded quotient ${\rm Gr}^rG$.
Since ${\rm Gr}^rG$ is a $p$-group,
if the induced morphism
on $\Theta^{(r)}$
for $K'$ is dominant,
the surjectivity of
the morphism on the fundamental groups
proved in Section \ref{ssExt}
induces a surjection
on the graded quotients as required.
The necessary functorial properties of
the vector space $\Theta^{(r)}$
are prepared in Section \ref{sIOm}
by studying certain normal sheaves
and modules of differentials.

The key point
in the second proof
in the case where $r>1$ is an integer
is the property
that the Galois covering of 
$\Theta^{(r)}$
of the Galois group
${\rm Gr}^rG$ 
carries a structure
of smooth group scheme
such that covering is
a morphism of group schemes.
This is proved by applying
a general criterion
prepared in Section \ref{ssGH}
to the geometric construction similar
to that in \cite{nonlog}.
The construction is a refinement of
that in \cite{III} and is
given in Section \ref{sskey}.
Further applying
the description
of extension groups of
a vector space over a field of
characteristic $p>0$ by ${\mathbf F}_p$
in Section \ref{ssExt},
we define the injection
{\rm (\ref{eqgrK})}.

We sketch the idea more precisely.
Both the first proof of Theorem \ref{thmgr}
and the construction of
(\ref{eqgrK}) are based on
the construction of 
the tangent space
of a local ring
at a geometric closed point
defined by an algebraic closure
of the residue field.
The construction is globalized
in \cite{LEX}.
Let $S={\rm Spec}\, {\mathcal O}_K$
and $\bar F$ be an algebraic closure
of the residue field $F
={\mathcal O}_K/{\mathfrak m}_K$.
We show in Proposition \ref{prH1} that
the cotangent complex
$L_{{\bar F}/S}$ is acyclic except at degree
$-1$ and
the cohomology group
$H_1(L_{{\bar F}/S})$ is an
$\bar F$-vector space fitting
in an exact sequence
\begin{equation}
0\to {\mathfrak m}_K/
{\mathfrak m}_K^2
\otimes_F\bar F
\to
H_1(L_{\bar F/S})
\to
\Omega^1_F\otimes_F\bar F
\to 0.
\label{eqH1bF}
\end{equation}
We define the tangent space
of $S={\rm Spec}\, {\mathcal O}_K$
at the geometric closed point
${\rm Spec}\, \bar F$
as the spectrum
\begin{equation}
\Theta_{K,\bar F}
={\rm Spec}\,
S(H_1(L_{\bar F/S}))
\label{eqTh}
\end{equation}
of the symmetric algebra
over $\bar F$.
If ${\mathcal O}_{K_0}\subset {\mathcal O}_K$
is a discrete valuation subring
with perfect residue field $k$,
we have a canonical surjection
\begin{equation}
H_1(L_{\bar F/S})\to 
\Omega^1_{{\mathcal O}_K/{\mathcal O}_{K_0}}
\otimes_{{\mathcal O}_K}\bar F
\label{eqH1LObF}
\end{equation}
and (\ref{eqH1LObF})
is an isomorphism
if the ramification index
$e_{K/K_0}$ is not 1,
by Proposition \ref{prLOm}.
The exact sequence 
(\ref{eqH1bF}) is canonically
identified with
the scalar extension of
that constructed in
\cite[Proposition 9.6.14]{GR}
by a different method
(\cite[Corollary 1.1.10]{LEX}).

We say that an extension $K'$ of
$K$ of discrete valuation fields
is dominant on the tangent spaces,
or tangentially dominant,
if the induced morphism
$S(H_1(L_{\bar F/S}))\to 
S(H_1(L_{\bar F'/S'}))$ of symmetric
algebras is an injection.
We derive Theorem \ref{thmgr}
from the following fact
by constructing a
tangentially dominant extension $K'$ of
$K$ with perfect residue field
in Proposition \ref{prperf}.

\medskip
\noindent
{\bf Corollary \ref{corRed1}} {\rm  (1)}
{\em 
Let $K$ be a henselian discrete valuation field and 
let $L$ be a finite Galois extension of $K$
of Galois group $G$.
Let $K'$ be a tangentially dominant extension of
henselian discrete valuation field,
let $L'=LK'$ be a composition field
and let $G'$ be the Galois group.
For $r>1$, the canonical injection
$G^{\prime r}
\to G^r$ is an isomorphism.
}

\medskip

We give an independent proof
of Theorem \ref{thmgr}
without the reduction to the classical case
and construct the injection (\ref{eqgrK})
for $r>1$.
Let ${\rm ord}_K$
also denote the extension
to a separable closure
$\bar K$ of the normalized
discrete valuation of $K$
and, for a rational number $r$,
set
${\mathfrak m}_{\bar K}^r
=\{x\in \bar K\mid
{\rm ord}_Kx\geqq r\}
\supset
{\mathfrak m}_{\bar K}^{r+}
=\{x\in \bar K\mid
{\rm ord}_Kx> r\}$
so that
the $\bar F$-vector space
${\mathfrak m}_{\bar K}^r
/{\mathfrak m}_{\bar K}^{r+}$
in  (\ref{eqgrK})
is of dimension 1.

We sketch the proof of Corollary \ref{corRed1} (1).
Using the property that the ramification groups
are compatible with quotient
\cite[Corollary 1.4.3]{TJM},
we reduce the proof to the 
case where $
G^{r+}=1$ and ${\rm Gr}^rG=G^r$.
The geometric construction of ramification groups recalled in Section \ref{ssFun}
provides a connected
$G^r$-torsor $X^{(r)\circ}$ over the
dual space $\Theta^{(r)}$
of the $\bar F$-vector space
on the target of (\ref{eqgrK}).
Let $k\to k'$ be
a morphism of
algebraically closed fields
of characteristic $p>0$
and let $E$ and $E'$ be
a $k$-vector space and
a $k'$-vector space.
For a $k'$-linear morphism
$E'\to E\otimes_kk'$,
the morphism 
$\pi_1(E',0)_{\text{pro-}p}
\to
\pi_1(E,0)_{\text{pro-}p}$
of the pro-$p$
quotients of the \'etale fundamental
groups is surjective if and only if the
morphism $E'\to E$ of schemes is dominant
by Corollary \ref{corpiE}.1.
Using this property and the fact that
$G^r$ is a $p$-group,
we prove
Corollary \ref{corRed1} (1).

We sketch the second proof of Theorem \ref{thmgr} and 
the construction of (\ref{eqgrK}).
In Section \ref{sGH}, we prove some properties Proposition \ref{prET} 
on finite \'etale isogenies
of smooth connected group schemes
and the classification
Proposition \ref{prExt}
of extensions by ${\mathbf F}_p$
of a vector space over a field
of characteristic $p>0$,
by the dual space.
Then,
it is reduced to showing that
the $G^r$-torsor $X^{(r)\circ}$
admits a group structure
such that the finite \'etale morphism
$X^{(r)\circ}\to \Theta^{(r)}$ is 
a morphism of group schemes over
$\bar F$.
The construction of a
group structure on $X^{(r)\circ}$
satisfying the required property
is reduced by Proposition
\ref{prisog} to a construction of
an automorphism
of the base change
$X^{(r)\circ}_{\bar F'}$
to an extension $\bar F'$
of the function field
of $\Theta^{(r)}$
compatible with
the translation of $\Theta^{(r)}_{\bar F'}$
by the generic point.

By using N\'eron's desingularization,
we reduce the construction to 
the case where
$K$ is essentially of finite
type and is separable
over a discrete valuation field
$K_0$ with
perfect residue field.
By constructing a suitable
ramified extension of $K$,
we further reduce it to the case
where the index $r\geqq2$ is an integer.

In the case
where $r\geqq2$ is an integer,
an automorphism
of the base change $X^{(r)\circ}_{\bar F'}$
satisfying the required property
is constructed in Section \ref{sskey} as follows.
First, we define a deformation $L'_1$ of 
the composition field $L_1=LK_1$ over
the local field $K_1$ at a point of
a smooth scheme over ${\mathcal O}_K$
appearing in the 
construction of the
$G^r$-torsor $X^{(r)\circ}$.
By the construction of $L'_1$,
the \'etale morphism $X^{(r)\circ}_1\to 
\Theta^{(r)}_1$ for the extension $L'_1$ of $K_1$
is given by the generic translation
of the scalar extension to ${\bar F'}$
of $X^{(r)\circ}\to 
\Theta^{(r)}$ for the original 
extension $L$ over $K$.
We prove that
over an unramified extension 
$K_2$ of $K_1$,
the composition field
$L'_2=L'_1K_2$ is isomorphic to
$LK_2$,
using the defining property
of ramification groups
and a descent property of
\'etaleness of finite covering
with a partial section
proved in Section \ref{ssNorm}.
This isomorphism
induces a required isomorphism from
$X^{(r)\circ}_1\to 
\Theta^{(r)}_1$ 
to
$X^{(r)\circ}_{\bar F'}\to 
\Theta^{(r)}_{\bar F'}$.

In the geometric case
where 
${\mathcal O}_K$ is the henselization
of a local ring ${\mathcal O}_{X,\xi}$
at the generic point $\xi$ of a smooth
divisor $D$
of a smooth scheme $X$
over a perfect field $k$,
such a group structure
is deduced in \cite{nonlog} from
the groupoid structure on
${\rm pr}_1,
{\rm pr}_2\colon X\times X\to X$
defined by
${\rm pr}_{13}\colon
(X\times X)\times_X(X\times X)
=
X\times X\times X\to X\times X$
and its modification along the
boundary obtained
by taking blow-up along $R=rD$
embedded in the diagonal
$X\subset X\times X$.
The construction in this article
is an approximation of the construction
in the geometric case.

The contents of each section
are briefly described at
the beginning of each section.

A crucial idea described above
of the proof
of the main result
in this article
was found during the preparation
of a talk at Chinese Academy of
Science in September 2019. 
The author thanks
W.\ Zheng for his hospitality 
during the stay.
The author thanks L.\ Illusie
for discussions 
and providing a historical remark
on cotangent complexes
and for careful reading 
and numerous comments
on the manuscript.
The author thanks O.\ Gabber
for indicating another construction
of tangent spaces and the reference
in \cite{GR}.
The author thanks 
N.\ Umezaki for discussion
on the attempt to construct
the automorphism of 
$X^{(r)\circ}_{\bar F'}$.
The research is partially supported
by Grant-in-Aid (B) 19H01780.

\tableofcontents

\section{Tangent spaces, immersions and differentials}
\label{sIOm}

First, we define a tangent space
at a geometric point of a scheme
as the dual of $H_1$ of the cotangent complex in Section \ref{ssLT}.
Using the tangent spaces,
we define a condition for
a morphism of schemes to
be tangentially dominant at a point.
This is weaker than formal smoothness
and will be used to construct
an extension of discrete valuation field
with perfect field in the first
proof of Theorem \ref{thmgr}.

In Section \ref{ssimm},
we study closed immersions
of the spectrum of a finite
extension of discrete valuation rings
to smooth schemes,
used in the geometric definition
of the filtration by ramification groups
in Section \ref{ssFun}.
We study the case where
the discrete valuation ring
is essentially of finite type
over a discrete valuation ring
with perfect residue field,
which will be used in a reduction
step to the case where the index
$r>1$ is an integer of
the second proof of Theorem \ref{thmgr}.

In Section \ref{ssTor},
we study the relation of
the conormal sheaf of 
closed immersion with the
tangent space introduced
in Section \ref{ssLT}
and with
the module of differentials.
The relation with the
tangent space will be used
in the reduction to the perfect
residue field case
and also to show the triviality
of a Galois action in Section \ref{ssFun}.
The relation with
the module of differentials
is used in the reduction
step to the case where the index
$r>1$ is an integer of
the second proof.

\subsection{Tangent space
at a point of a scheme}\label{ssLT}

We briefly recall basic facts on
cotangent complexes from
\cite[Chapitres II, III]{Ill}.
For a morphism of schemes
$X\to S$, the cotangent complex
$L_{X/S}
$ is defined 
\cite[Chapitre II, 1.2.3]{Ill}
as a chain complex of flat
${\mathcal O}_X$-modules, whose cohomology
sheaves are quasi-coherent.
If $X={\rm Spec}\, A$ and $S={\rm Spec}\, B$
are affine,
there exist a complex 
$L_{A/B}$ of flat $A$-modules
and a canonical quasi-isomorphism 
$L_{A/B}\otimes_A{\mathcal O}_X
\to L_{X/S}$.

The cohomology sheaf
${\mathcal H}_1(L_{X/S})$
is closely related to the module of imperfection studied
in \cite[Chapitre 0, Section 20.6]{EGA4}.
Later Grothendieck
introduced the truncation
$\tau_{[-1,0]}L_{X/S}$
denoted $L^{X/S}$
in his prenotes for Hartshorne's
seminar on Residues and Duality
and called 
${\mathcal H}_1(L^{X/S})$
the conormal module of
$X$ relatively to $S$.
In this section,
we study the case
where $X={\rm Spec}\, E$
for a field $E$ of characteristic $p>0$
and show that the $E$-vector space
$H_1(L_{E/S})$ plays a role
of the cotangent space of $S$
at the point defined by
${\rm Spec}\, E\to S$.
We study a globalization in
\cite[Section 1]{LEX}.

There is a canonical isomorphism
${\mathcal H}_0(L_{X/S})\to \Omega^1_{X/S}$
\cite[Chapitre II, Proposition 1.2.4.2]{Ill}.
This induces a canonical morphism
$L_{X/S}\to \Omega^1_{X/S}[0]$.
For a commutative diagram
\begin{equation}
\begin{CD}
X'@>>> S'\\
@VfVV
@VVV\\
X@>>> S,
\end{CD}
\label{eqXS}
\end{equation}
a canonical morphism
$Lf^*L_{X/S}\to L_{X'/S'}$
is defined \cite[Chapitre II, (1.2.3.2)$'$]{Ill}. 
For a morphism
$f\colon X\to Y$
of schemes over a scheme $S$,
a distinguished triangle
\begin{equation}
Lf^*L_{Y/S}\to L_{X/S}\to L_{X/Y}\to
\label{eqLXYS}
\end{equation}
is defined
\cite[Chapitre II, Proposition 2.1.2]{Ill}.

If $X\to S$ is a closed immersion 
defined by the ideal sheaf ${\mathcal I}_X
\subset {\mathcal O}_S$ and
if $N_{X/S}={\mathcal I}_X/{\mathcal I}_X^2$ denotes the conormal
sheaf,
there exists a canonical isomorphism
${\mathcal H}_1(L_{X/S})\to N_{X/S}$
\cite[Chapitre III, Corollaire 1.2.8.1]{Ill}.
This induces a canonical morphism
$L_{X/S}\to N_{X/S}[1]$.

\begin{lm}\label{lmNOm}
{\rm 1.
(\cite[Chapitre III, Proposition 1.2.9]{Ill})}
Let $f\colon X\to Y$ be
an immersion of schemes over
a scheme $S$.
Then, the boundary morphism
$\partial\colon N_{X/Y}\to 
f^*\Omega^1_{Y/S}$
of the distinguished triangle
$Lf^*L_{Y/S}\to L_{X/S}\to L_{X/Y}
\to $ sends
$g$ to $-dg$.

{\rm 2.} {\rm (\cite[Chapitre III, 
Proposition 3.1.2 (i)$\Rightarrow$(ii)]{Ill})}
Let $X\to S$ be a smooth morphism.
Then, 
the canonical morphism
$L_{X/S}\to \Omega^1_{X/S}[0]$
is a quasi-isomorphism.

{\rm 3.} {\rm (\cite[Chapitre III, 
Proposition 3.2.4 (iii)]{Ill})}
If $X\to S$ is a regular immersion,
the canonical morphism
$L_{X/S}\to N_{X/S}[1]$
is a quasi-isomorphism.

{\rm 4. (\cite[9.5.8]{17}
cf.\ \cite[Chapitre III, 
Remarque 3.1.4]{Ill})}
Let $f\colon X\to S$ be a formally
smooth morphism.
Then, $\Omega^1_{X/S}$ is a flat
${\mathcal O}_X$-module
and $H_1(L_{X/S})=0$.
\end{lm}

\begin{lm}\label{lmH0}
Let $F\to E$ be an extension
of fields.
Then,
the cotangent complex
$L_{E/F}$ is acyclic except at degree $[-1,0]$.
If $E$ is a separable extension
of $F$,
the canonical morphism
$L_{E/F}\to 
\Omega^1_{E/F}[0]$
is a quasi-isomorphism.
If $F$ is of characteristic $p>0$
and if $k\subset F^p$
is a subfield,
the minus of the boundary 
morphism of 
the distinguished triangle
$L_{F/k}\otimes_FE
\to L_{E/k}\to L_{E/F}\to $
defines an isomorphism
\begin{equation}
-\partial\colon H_1(L_{E/F})
\to
{\rm Ker}(\Omega^1_{F/k}
\otimes_FE
\to \Omega^1_{E/k}).
\label{eqLEF}
\end{equation}
\end{lm}

\proof{
First, we show the case where
$E$ is a separable extension of $F$.
By taking the limit
\cite[(1.2.3.4)]{Ill},
we may assume that $E$ 
is finitely generated over $F$.
We may further assume that $E$ is 
the function field of a smooth scheme
$X$ over $F$.
Then, we have a quasi-isomorphism
$L_{X/F}\to \Omega^1_{X/F}[0]$ by 
Lemma \ref{lmNOm}.2.
Hence
$L_{E/F}\to \Omega^1_{E/F}[0]$
is a quasi-isomorphism
and $\Omega^1_{F/k}
\otimes_FE
\to \Omega^1_{E/k}$
is an injection in this case.

To show the remaining assertion,
we may assume that
$F$ is of characteristic $p>0$.
Let $k_0\subset F$ be a perfect subfield,
say $k_0={\mathbf F}_p$.
Then, since $F$ and $E$ are
separable extensions of $k_0$,
the distinguished triangle
$L_{F/k_0}\otimes_FE
\to L_{E/k_0}\to L_{E/F}\to $
implies that
$L_{E/F}$ is acyclic except at degree $[-1,0]$
and defines an isomorphism 
$H_1(L_{E/F})\to 
{\rm Ker}(\Omega^1_{F/k_0}
\otimes_FE
\to \Omega^1_{E/k_0})$
(\ref{eqLEF}) for $k=k_0$.
If $k\subset F^p$
is a subfield,
$\Omega^1_{F/k}\to
\Omega^1_{F/F^p}$
and
$\Omega^1_{E/k}\to
\Omega^1_{E/F^p}$
are isomorphisms
and the right hand side
of (\ref{eqLEF})
is independent of such $k$,
including $k_0$.
\qed

}

\begin{pr}\label{prH1}
Let $S$ be a scheme,
$s\in S$ be a point
and $E$ be an extension of 
the residue field $F=k(s)$ at $s$.

{\rm 1.}
The distinguished triangle
\begin{equation}
\to L_{F/S}\otimes_FE
\to L_{E/S}\to L_{E/F}\to
\label{eqLESF}
\end{equation}
of cotangent complexes
and the isomorphism
{\rm (\ref{eqLEF})}
define an isomorphism
$H_0(L_{E/S})
\to \Omega^1_{E/F}$ and
an exact sequence
\begin{equation}
0\to {\mathfrak m}_s/
{\mathfrak m}_s^2
\otimes_FE
\to
H_1(L_{E/S})
\overset{-\partial}\longrightarrow
{\rm Ker}
(\Omega^1_F\otimes_FE
\to
\Omega^1_E)
\to 0
\label{eqH1L}
\end{equation}
of $E$-vector spaces.

{\rm 2.}
Assume
that $F$ is of characteristic $p>0$
and that $E$ is an extension of $F^{1/p}$.
Then,
the exact sequence {\rm (\ref{eqH1L})}
defines an exact sequence
\begin{equation}
0\to {\mathfrak m}_s/
{\mathfrak m}_s^2
\otimes_FE
\to
H_1(L_{E/S})
\overset{-\partial}\longrightarrow
\Omega^1_F\otimes_FE
\to 0.
\label{eqH1L1}
\end{equation}

{\rm 3.}
Let $E\subset E'$ be extensions of $F$.
Then, the canonical morphism
\begin{equation}
H_1(L_{E/S})\otimes_EE'
\to
H_1(L_{E'/S})
\label{eqEF1p}
\end{equation}
is an injection.
This is an isomorphism
either if $F$ is of characteristic $0$
or if
 $F$ is of characteristic $p>0$
 and $E$ is an extension of
 $F^{1/p}$.

{\rm 4.}
If $S$ is regular at $s$,
then
the cotangent complex
$L_{E/S}$
is acyclic except at degree $[-1,0]$.
Further if
$\Omega^1_{E/F}=0$,
the complex $L_{E/S}$
is acyclic except at degree $-1$.
\end{pr}

The exact sequence 
(\ref{eqH1L1}) 
in the case $E=F^{1/p}$
is canonically
identified with that
constructed in
\cite[Proposition 9.6.14]{GR}
by a different method
(\cite[Corollary 1.1.10]{LEX}).
The converse of
Proposition \ref{prH1}.4 holds
by  \cite[Corollary 9.6.45]{GR}.
The exact sequence (\ref{eqH1L})
implies that
$H_1(L_{E/S})$
is determined by
${\mathcal O}_{S,s}/
{\mathfrak m}_s^2
\to E$.

\proof{
1.
We have
$H_1(L_{F/S})=
{\mathfrak m}_s/
{\mathfrak m}_s^2$
and
$H_0(L_{F/S})=0$.
By Lemma \ref{lmH0},
we have
$H_0(L_{E/F})=\Omega^1_{E/F}$,
(\ref{eqLEF}) and
$H_2(L_{E/F})=0$.
Hence 
the distinguished triangle
(\ref{eqLESF})
defines an isomorphism
$H_0(L_{E/S})
\to \Omega^1_{E/F}$ and
the exact sequence (\ref{eqH1L}).

2.
Assume $F^{1/p}\subset E$.
Then, the morphism
$\Omega^1_F\otimes_FE
\to
\Omega^1_E$
is $0$
and we have
${\rm Ker}
(\Omega^1_F\otimes_FE
\to
\Omega^1_E)=
\Omega^1_F\otimes_FE$.
Hence  (\ref{eqH1L}) gives
(\ref{eqH1L1}).

3.
Since
$L_{E'/E}$ is acyclic except at degree
$[-1,0]$ 
by Lemma \ref{lmH0},
the distinguished triangle
$L_{E/S}\otimes_EE'\to
L_{E'/S}\to L_{E'/E}\to$
implies the injectivity
of (\ref{eqEF1p}).

If $F$ is of characteristic $0$,
we have
$H_1(L_{E'/E})=0$ and
(\ref{eqEF1p}) is a surjection.
If $F^{1/p}\subset E$,
the boundary mapping
$H_1(L_{E'/E})
\to \Omega^1_{E/F}\otimes_EE'$
is an injection by Lemma \ref{lmH0}.
Hence 
(\ref{eqEF1p}) is a surjection.

4.
The morphism
$L_{F/S}\to N_{F/S}[1]$
is a quasi-isomorphism by Lemma \ref{lmNOm}.3.
Hence by Lemma \ref{lmH0},
the distinguished triangle
(\ref{eqLESF})
shows that
$L_{E/S}$
is acyclic except at degree $[-1,0]$.
Further if $H_0(L_{E/S})=\Omega^1_{E/F}$
is $0$,
then $L_{E/S}$
is acyclic except at degree $-1$.
\qed

}

\medskip

For a morphism
$S'\to S$
of schemes
sending a point
$s'\in S'$ to $s\in S$
and a morphism $E\to E'$
of extensions 
extending the morphism
$F=k(s)\to F'=k(s')$ of
the residue fields,
by the functoriality of cotangent complexes,
we have a canonical morphism
\begin{equation}
H_1(L_{E/S})
\to
H_1(L_{E'/S'}).
\label{eqH1Lf}
\end{equation}

Let $S$ be a scheme,
$s\in S$ be a point
and $u\in {\mathcal O}_{S,s}$.
Assume that
the residue field $F=k(s)$ is of characteristic $p>0$
and let $E$ be an extension of $F$
containing a
$p$-th root $v\in E$ of $\bar u\in F$.
Then, 
\begin{equation}
\widetilde d u\in H_1(L_{E/S})
\label{eqdfdu}
\end{equation}
is defined in \cite[Definition 1.1.6]{LEX}.
We briefly recall the construction.
Let
$W\subset  
{\mathbf A}^1_S
=S\times_{{\rm Spec}\, {\mathbf Z}}
{\rm Spec}\, {\mathbf Z}[T]$
be the closed subscheme
defined on a neighborhood of
$s$ by $T^p-u$
and define a morphism
${\rm Spec}\, E\to W$
over $S$
by sending $T$ to $v$.
Then, the image of
$\tilde d u$ by
the injection
$H_1(L_{E/S})
\to H_1(L_{E/{\mathbf A}^1_S})$
is the image of
the section $u-T^p$
of the conormal sheaf
$N_{W/{\mathbf A}^1_S}$.
The construction of 
$\widetilde du$ is functorial.
The following property
is a special case of
\cite[Proposition 1.1.5]{LEX}.

\begin{lm}\label{lmdu}
Let $S$ be a scheme,
$s\in S$ be a point
and $u\in {\mathcal O}_{S,s}$.
Assume that the residue field
$F=k(s)$ is of characteristic $p>0$
and let $E$ be an extension
of $F$ containing a $p$-th root $v$
of $\bar u\in F$.

{\rm 1. (cf.\ \cite[9.6.12]{GR})}
Let $u,u'\in {\mathcal O}_{S,s}$
such that there exist 
$p$-th roots $v,v'\in E$ of $\bar u,\bar u'$.
Then, we have
\begin{align}
\widetilde d(u+u')
&\, =
\widetilde du+
\widetilde du'+
{\displaystyle\sum_{i=1}^{p-1}}
\dfrac{v^i}{i!}
\dfrac{v^{\prime (p-i)}}{(p-i)!}
\cdot \widetilde dp,
\label{eqadd}
\\
\widetilde d(uu')
&\, =
u'\cdot \widetilde du+
u\cdot \widetilde du'.
\label{eqLb}
\end{align}


{\rm 2.}
The image of
$\widetilde du \in H_1(L_{E/S})$
{\rm (\ref{eqdfdu})}
by ${-\partial}
\colon H_1(L_{E/S})
\to
\Omega^1_F\otimes_FE$
{\rm (\ref{eqH1L})}
is $d\bar u$.

{\rm 3.}
If $u\equiv 0\bmod {\mathfrak m}_s$,
we have
$\omega\equiv u\in 
{\mathfrak m}_s/
{\mathfrak m}_s^2
=H_1(L_{F/S})=
H_1(L_{E/S})$.
\end{lm}

By (\ref{eqadd}), if one of $u,u'$
is an element of the maximal ideal
${\mathfrak m}_s\subset
{\mathcal O}_{S,s}$,
we have
$\widetilde d(u+u')
=
\widetilde du+
\widetilde du'$.
By (\ref{eqLb}),
we have
$\widetilde du^i=iu^{i-1}\widetilde du$
by induction on integer $i\geqq 0$.
Hence for an integer $a$,
since $a-a^p$ is divisible by $p$,
we have
$$\widetilde d a=
\widetilde d a^p+
\widetilde d (a-a^p)
=
pa^{p-1}\cdot \widetilde d a+
p\cdot \widetilde d\, \dfrac{a-a^p}p
+
\dfrac{a-a^p}p\cdot 
\widetilde d p
=
\dfrac{a-a^p}p\cdot 
\widetilde d p.$$
In the case where
$v\notin F=k(s)$,
the following remark is due
to Luc Illusie.
Let $S'={\rm Spec}\,
{\mathcal O}_{S,s}[T]/(T^p-u)\to S$ 
and identify $s'={\rm Spec}\, F(v)$ 
with the closed point of $S'$.
Then, the morphism
$H_1(L_{F(v)/S})\to H_1(L_{F(v)/S'})$
defines a splitting of the exact sequence
$0\to H_1(L_{F/S})\otimes_FF(v)
\to H_1(L_{F(v)/S})\to H_1(L_{F(v)/F})\to 0$
defined by the distinguished
triangle
$L_{F/S}\otimes^L_FF(v)
\to L_{F(v)/S}\to L_{F(v)/F}\to$
and the image of
$\widetilde d u\in H_1(L_{F(v)/S})$
in 
$H_1(L_{F(v)/S'})$ is
$\widetilde d T^p=
p T^{p-1}\widetilde d T=0$
by the Leibniz rule.

\begin{pr}\label{prLOm}
Let $S\to S_0$ be a morphism
of schemes,
let $s\in S$ be a point
and $s_0\in S_0$ be its image.
Let $F=k(s)$ be the residue field and
assume that $k=k(s_0)$ is a
perfect field of characteristic $p>0$.
Let $E$ be an extension of $F$
and define a morphism
\begin{equation}
{-\partial}\colon H_1(L_{E/S})
\to
\Omega^1_{S/S_0}
\otimes_{{\mathcal O}_S}E.
\label{eqH1LOm}
\end{equation}
to be the minus of the boundary morphism 
of the distinguished triangle
$L_{S/S_0}\otimes^L_{{\mathcal O}_K}E
\to L_{E/S_0}\to L_{E/S}\to $
defined by ${\rm Spec}\, E\to S\to S_0$.

{\rm 1. (cf.\ \cite[Lemma 9.6.3]{GR})}
For $E=F$, 
the above distinguished triangle
defines an exact sequence
\begin{equation}
0\to
{\mathfrak m}_s/
({\mathfrak m}_s^2
+{\mathfrak m}_{s_0}
{\mathcal O}_{S,s})
\overset d\longrightarrow
\Omega^1_{S/S_0}
\otimes_{{\mathcal O}_S}F
\to 
\Omega^1_{F/k}
\to 0
\label{eqOmm}
\end{equation}
of $F$-vector spaces.

{\rm 2. (\cite[Proposition 1.1.5.4]{LEX})}
For $u\in {\mathcal O}_{S,s}$,
the mapping {\rm (\ref{eqH1LOm})}
sends $\widetilde d u$ to
$du$.

{\rm 3.}
The mapping {\rm (\ref{eqH1LOm})}
fits in the commutative diagram
\begin{equation}
\xymatrix{
H_1(L_{E/S})\ar[rd]^{(\ref{eqH1LOm})}
\ar[r]^{\rm(\ref{eqH1L})}
&
\Omega^1_{F/k}
\otimes_FE
\\
{\mathfrak m}_s/
{\mathfrak m}_s^2
\otimes_FE
\ar[r]^d\ar[u]^{\rm(\ref{eqH1L})}
&
\Omega^1_{S/S_0}
\otimes_{{\mathcal O}_S}E
\ar[u]
}
\label{eqH1Om}
\end{equation}
of $E$-vector spaces.
The kernel and the cokernel
of {\rm (\ref{eqH1LOm})}
are the images
of 
${\mathfrak m}_{s_0}/
{\mathfrak m}_{s_0}^2
\otimes_FE
\to 
{\mathfrak m}_s/
{\mathfrak m}_s^2
\otimes_FE$
and of
$\Omega^1_{F/k}
\otimes_FE
\to 
\Omega^1_{E/k}$.
Consequently,
the morphism 
{\rm (\ref{eqH1LOm})} is a surjection
(resp.~an injection)
if $F$ is a subfield of $E^p$
(resp. if ${\mathfrak m}_{s_0}
{\mathcal O}_{S,s}
\subset {\mathfrak m}_s^2$).
\end{pr}


\proof{
1.
By Lemma \ref{lmNOm}.1,
the minus of the boundary morphism 
of the distinguished triangle
$L_{S/S_0}\otimes^L_{{\mathcal O}_S}F
\to L_{F/S_0}\to L_{F/S}\to$
defines an exact sequence
\begin{equation}
H_1(L_{F/S_0})
\to
{\mathfrak m}_s/
{\mathfrak m}_s^2
\overset d\longrightarrow
\Omega^1_{S/S_0}
\otimes_{{\mathcal O}_S}F
\to 
\Omega^1_{F/k}
\to 0
\label{eq115}
\end{equation}
of $F$-vector spaces.
Since $\Omega^1_k=0$,
the exact sequence (\ref{eqH1L})
gives an isomorphism
${\mathfrak m}_{s_0}/
{\mathfrak m}_{s_0}^2
\otimes_kF\to
H_1(L_{F/S_0})$
and we obtain (\ref{eqOmm}).

3.
The commutative diagram
$$\begin{CD}
L_{S/S_0}\otimes^L
_{{\mathcal O}_S}F
\otimes_FE
@=
L_{S/S_0}\otimes^L
_{{\mathcal O}_S}E
@>>>
L_{F/S_0}
\otimes_FE
\\
@VVV@VVV@VVV\\
L_{F/S_0}
\otimes_FE
@>>>
L_{E/S_0}
@=
L_{E/S_0}
\\
@VVV@VVV@VVV\\
L_{F/S}
\otimes_FE
@>>>
L_{E/S}
@>>>
L_{E/F}
\\
@VVV@VVV@VVV
\end{CD}
$$
of distinguished triangles
defines the commutative 
diagram (\ref{eqH1Om})
by Lemma \ref{lmNOm}.1.
Similarly as in the proof of 1,
we obtain an exact sequence
\begin{equation}
{\mathfrak m}_{s_0}/
{\mathfrak m}_{s_0}^2
\otimes_kE
\to
H_1(L_{E/S})
\to
\Omega^1_{S/S_0}
\otimes_{{\mathcal O}_S}E
\to 
\Omega^1_{E/k}.
\label{eq116}
\end{equation}
Since $k$ is perfect,
we have
$\Omega^1_{F/k}
=\Omega^1_F$.
Hence the assertion on the 
kernel and the cokernel
follows from (\ref{eq116}),
(\ref{eq115})$\otimes_FE$,
(\ref{eqH1L}).

If $F^{1/p}\subset E$,
then 
$\Omega^1_F\otimes_FE
\to
\Omega^1_E$ is the $0$-mapping.
If ${\mathfrak m}_{s_0}
{\mathcal O}_{S,s}
\subset {\mathfrak m}_s^2$,
then 
${\mathfrak m}_{s_0}/
{\mathfrak m}_{s_0}^2
\to 
{\mathfrak m}_s/
{\mathfrak m}_s^2$
is the $0$-mapping.
\qed

}

\begin{df}\label{dfcot}
Let $S$ be a scheme,
$s\in S$ be
a point and $F=k(s)$
be the residue field.

{\rm 1.}
For an extension $E$ of $F$,
we call the spectrum
\begin{equation}
\Theta_{S,E}=
{\rm Spec}\, S(H_1(L_{E/S}))
\label{eqThA}
\end{equation}
of the symmetric algebra
over $E$ 
the tangent space of $S$ at $E$.
If $S={\rm Spec}\, {\mathcal O}_K$
for a discrete valuation ring ${\mathcal O}_K$
and if $s\in S$ is the closed point,
we also write
$\Theta_{S,E}=\Theta_{K,E}$.

{\rm 2.}
Let $S'\to S$ be a morphism
of schemes and
let $s\in S$ be the image of
$s'\in S'$.
Let $\bar F\to \bar F'$
be a morphism of algebraic closures
of the residue fields 
extending $F=k(s)\to F'=k(s')$.
We say that 
the morphism $S'\to S$ is 
dominant on the tangent spaces,
or {\em tangentially dominant},
at $s'$ if the morphism
$$
S(H_1(L_{\bar F/S}))
\to
S(H_1(L_{\bar F'/S'}))$$
is an injection.

{\rm 3.}
Let ${\mathcal O}_K$ and ${\mathcal O}_{K'}$
be discrete valuation rings.
We say that a morphism
${\mathcal O}_K\to {\mathcal O}_{K'}$  
is a morphism 
of discrete valuation rings
if it is faithfully flat.
We say that
a morphism
${\mathcal O}_K\to {\mathcal O}_{K'}$  
of discrete valuation rings
is tangentially dominant if
$S'={\rm Spec}\, {\mathcal O}_{K'}
\to S={\rm Spec}\, {\mathcal O}_K$
is dominant on the tangent space
at the closed point of $S'$.
\end{df}

For $E=F$,
we recover the definition of
the Zariski tangent space.
That a morphism $S'\to S$ is 
dominant on the tangent spaces at $s'$
means that
the morphism
$\Theta_{S',\bar F'}
\to
\Theta_{S,\bar F}$
of schemes is dominant.

\begin{lm}\label{lmH1}
Let $S$ be a scheme,
$s\in S$ a point
and $E$ be an extension
of the residue field $F=k(s)$.

%
%
Let $S'\to S$ be a morphism
of schemes
sending $s'\in S'$ to $s\in S$
and let $E\to E'$ be a morphism
of extensions of residue fields.
If $S'\to S$ is formally smooth at $s'$,
the morphism {\rm (\ref{eqH1Lf})}
induces an injection
\begin{equation}
H_1(L_{E/S})\otimes_EE'
\to 
H_1(L_{E'/S'})
\label{eqHLSS}
\end{equation}
of $E'$-vector spaces.
\end{lm}

\proof{
By Proposition \ref{prH1}.3, 
we may assume that $E=E'$.
By Lemma \ref{lmNOm}.4,
we have 
$H_1(L_{S'/S}\otimes^L_SE')=0$.
Hence the distinguished triangle
$L_{S'/S}\otimes^L_SE'\to
L_{E'/S}\to L_{E'/S'}$
implies the injectivity
of (\ref{eqHLSS}).
\qed

}

\begin{pr}\label{prcot}
{\rm 1.}
Let $S'\to S$
be a morphism of schemes
and let $s\in S$
be the image of $s'\in S'$.
We consider the following conditions:

{\rm (1)}
The morphism
$S'\to S$ is formally smooth at $s'$.

{\rm (2)}
The morphism
$S'\to S$
is dominant on the tangent spaces
at $s'$.

{\rm (3)}
The morphism
${\mathfrak m}_s/
{\mathfrak m}_s^2
\to
{\mathfrak m}_{s'}/
{\mathfrak m}_{s'}^2$
is injective.

\noindent
Then, we have
the implications
{\rm (1)}$\Rightarrow${\rm (2)}$\Rightarrow${\rm (3)}.

{\rm 2.}
Let ${\mathcal O}_K\to {\mathcal O}_{K'}$
be a morphism of discrete valuation
rings.
We consider the following conditions:

{\rm (1)}
The ramification index 
$e_{K'/K}$ is $1$
and $F'={\mathcal O}_{K'}/
{\mathfrak m}_{K'}$
is a separable extension of 
$F={\mathcal O}_K/
{\mathfrak m}_K$.

{\rm (2)}
The morphism
${\mathcal O}_K\to {\mathcal O}_{K'}$
is tangentially dominant.

{\rm (3)}
The ramification index 
$e_{K'/K}$ is $1$.

\noindent
Then, we have
the implications
{\rm (1)}$\Rightarrow${\rm (2)}$\Rightarrow${\rm (3)}.
\end{pr}

\proof{
1.
(1)$\Rightarrow$(2):
Let $\bar F\to \bar F'$
be a morphism of algebraic closures of
$F=k(s)$
and $F'=k(s')$.
Then, by Lemma \ref{lmH1}.2,
the morphism
$H_1(L_{\bar F/S})\otimes_{\bar F}\bar F'
\to
H_1(L_{\bar F'/S'})$ is an injection.
Hence
$S(H_1(L_{\bar F/S}))
\to
S(H_1(L_{\bar F/S}))
\otimes_{\bar F}\bar F'
\to
S(H_1(L_{\bar F'/S'}))$ is an injection.

(2)$\Rightarrow$(3):
By Proposition \ref{prH1},
the injection
$S(H_1(L_{\bar F/S}))
\to
S(H_1(L_{\bar F'/S'}))$ induces 
an injection
${\mathfrak m}_s
/{\mathfrak m}_s^2
\to 
{\mathfrak m}_{s'}
/{\mathfrak m}_{s'}^2$.

2.
Let
${\mathcal O}_K\to {\mathcal O}_{K'}$
be a morphism 
of discrete valuation rings.
Let $S'={\rm Spec}\,
{\mathcal O}_{K'}\to S={\rm Spec}\,
{\mathcal O}_K$ be the associated morphism
and $s'\in S'$ be the closed point.
Then
by \cite[Chapitre 0, Th\'eor\`eme (19.7.1)]{EGA4},
the condition (1) in 1.\
is equivalent to that
${\mathcal O}_{K'}/{\mathfrak m}_K
{\mathcal O}_{K'}$
is formally smooth over the residue field
$F={\mathcal O}_K/{\mathfrak m}_K$.
By \cite[Chapitre 0, Corollaire (19.6.5)]{EGA4},
this implies that
${\mathfrak m}_{K'}=
{\mathfrak m}_K
{\mathcal O}_{K'}$ and that
$F'={\mathcal O}_{K'}/{\mathfrak m}_{K'}$
is a separable extension of
$F$. Since the converse is true by
\cite[Chapitre 0, Th\'eor\`eme (19.6.1)]{EGA4},
the condition (1) in 1.\
is equivalent to
the condition (1) in 2.

Further,
the condition (3) in 1.\
is equivalent to
the condition (3) in 2.
\qed

}
\medskip

The implications
{\rm (1)}$\Rightarrow${\rm (2)}$\Rightarrow${\rm (3)}
are strict as follows.

\begin{ex}\label{ex11}
{\rm Let ${\mathcal O}_K\to {\mathcal O}_{K'}$ be a
morphism of discrete valuation rings.
Assume that the residue field
$F$ is of characteristic $p>0$
and let $u\in {\mathcal O}_K$
and $v\in {\mathcal O}_{K'},
w\in {\mathfrak m}_{K'}$
be elements satisfying
$u=v^p+w$.
Then, by Lemma {\rm \ref{lmdu}}.1,
we have
\begin{equation}
\widetilde du=pv^{p-1}\widetilde d v+
\widetilde dw
=\widetilde dw
\label{equvw}
\end{equation}
in $H_1(L_{\bar F'/S})$.
Hence, if $du\neq 0$ in $\Omega^1_F$
and $w\in {\mathfrak m}_{K'}^2$,
the morphism
${\mathcal O}_K\to {\mathcal O}_{K'}$ is
not tangentially dominant.

Concretely, 
for $K=k(u)((t))$,
the separable extension $K'=K(v)$ 
of degree $p$
defined by $v^p-t^2v=u$ 
is not tangentially dominant.
Hence the condition (3) is
satisfied but (2) is not.
}
\qed
\end{ex}

The following Proposition
gives an example of
extension satisfying
the condition (2) 
but not (1).

\begin{pr}\label{prperf}
Let $K$ be a discrete valuation field.
Then, there exists a
tangentially dominant morphism
${\mathcal O}_K\to {\mathcal O}_{K'}$
of discrete valuation rings
such that
the residue field $F'$ is perfect.
\end{pr}

\proof{
It suffices to consider the case where the
residue field $F$ is of
characteristic $p>0$.
Let $t\in {\mathcal O}_K$
be a uniformizer and let
$(u_i)_{i\in I}$
be a $p$-basis of the residue
field $F$.
For each $i\in I$,
take a lifting of $u_i\in F$
in ${\mathcal O}_K$ and let it also
denoted $u_i$ by abuse of notation.
The localization ${\mathcal O}_{K_0}$
of the polynomial ring
${\mathcal O}_K[w_{i,0}\ (i\in I)]$
at the prime ideal $(t)$
is a discrete valuation ring.
We set
$v_{i,0}=u_i-w_{i,0}t\in {\mathcal O}_{K_0}$ for each $i\in I$
so that $(v_{i,0},w_{i,0})_{i\in I}$
is a $p$-basis of
the residue field $F_0=F(w_{i,0}\ (i\in I))$
of $K_0$.

We define an inductive
system $({\mathcal O}_{K_n})_{n\in {\mathbf N}}$
of discrete valuation rings inductively by
${\mathcal O}_{K_{n}}
={\mathcal O}_{K_{n-1}}
[v_{i,n},w_{i,n}\ (i\in I)]
/(v_{i,n}^p-v_{i,n-1},
w_{i,n}^p-w_{i,n-1}\ (i\in I))$.
The inductive limit
${\mathcal O}_{K'}=\varinjlim_{n\in {\mathbf N}}
{\mathcal O}_{K_{n}}$
is a discrete valuation ring
and the residue field
$F'$ of $K'$
is a perfect closure of $F_0$
and is perfect.

The uniformizer
$t$ defines a basis of
${\mathfrak m}_K/
{\mathfrak m}_K^2$
and $(du_i)_{i\in I}$
defines a basis of
$\Omega^1_F$.
Hence by the exact sequence
$
0\to {\mathfrak m}_K/
{\mathfrak m}_K^2
\otimes_F\bar F
\to
H_1(L_{\bar F/S})
\to
\Omega^1_F\otimes_F\bar F
\to 0$
(\ref{eqH1L1})
and Lemma \ref{lmdu},
$\widetilde dt$ together with
$(\widetilde du_i)_{i\in I}$
defines a basis of $H_1(L_{\bar F/S})$.
By (\ref{equvw}),
the images of
$\widetilde du_i \in H_1(L_{\bar F/S})$
in $H_1(L_{\bar F'/S'})
={\mathfrak m}_{K'}/
{\mathfrak m}_{K'}^2
\otimes_{F'}\bar F'$
is $w_{i,0}t$.
Since $(w_{i,0})_{i\in I}$
are algebraically independent
over $\bar F$,
the morphism
$S(H_1(L_{\bar F/S}))\to
S(H_1(L_{\bar F'/S'}))$
is an injection.
\qed

}

\subsection{Immersions
to smooth schemes}\label{ssimm}

\begin{lm}\label{lmQT}
Let ${\mathcal O}_K\to {\mathcal O}_L$
be a finite morphism of
discrete valuation rings
and set
$S={\rm Spec}\, {\mathcal O}_K$ and
$T={\rm Spec}\, {\mathcal O}_L$.

{\rm 1.}
Let
$T\to Q$
be an immersion to a smooth scheme over
$S$.
Then,
after replacing $Q$ by
a neighborhood of $T$,
there exist a smooth scheme $P$ over $S$
and a cartesian diagram
\begin{equation}
\begin{CD}
Q@<<< T\\
@VVV @VVV\\
P@<<<S
\end{CD}
\label{eqPQTS}
\end{equation}
of schemes over $S$
such that
$Q\to P$
is quasi-finite and flat.

{\rm 2.}
Let $T\to Q$ and $T\to Q'$
be closed immersions
to smooth schemes over $S$.
Then, 
after replacing $Q$ by
a neighborhood of $T$,
there exists a commutative
diagram 
\begin{equation}
\xymatrix{
T\ar[r]\ar[rd]\ar[d]
&Q'\\
Q&Q''\ar[l]\ar[u]}
\label{eqQQQ}
\end{equation}
of schemes over $S$
such that 
$Q''\to Q'$ is \'etale.

{\rm 3.}
Let $T\to Q$ and $T\to Q'$
be closed immersions
to smooth schemes over $S$
and let 
\begin{equation}
\xymatrix{
T\ar[r]\ar[d]
&Q'\ar[ld]\\
Q&}
\label{eqTQQ'}
\end{equation}
be a commutative diagram
of schemes over $S$
such that $Q'\to Q$
is smooth.
Then, 
after replacing $Q'$ by
a neighborhood of $T$,
there exists
a commutative diagram
\begin{equation}
\xymatrix{
T\ar[r]\ar[d]
&Q'\ar[d]\ar[ld]\\
Q&{\mathbf A}^n_Q
\ar[l]
}
\label{eqTQQ}
\end{equation}
of schemes over $S$
such that
$Q'\to {\mathbf A}^n_Q$
is \'etale,
${\mathbf A}^n_Q\to Q$
is the projection
and the composition
$T\to {\mathbf A}^n_Q$
is the composition
of $T\to Q$
with the $0$-section
$Q\to {\mathbf A}^n_Q$.
\end{lm}

\proof{
1.
Since $T$ and $Q$ are regular,
the immersion
$T\to Q$ is a regular immersion
by \cite[Chap.\ 0 Corollaire (17.1.9) a)$\Rightarrow$ b)]{EGA4}.
After replacing $Q$ by
a neighborhood of $T$,
we may assume that
the closed subscheme
$T\subset Q$ is
defined by a regular
sequence $t_1,\ldots,t_n$.
We define
$Q\to P={\mathbf A}^n_S$
by $t_1,\ldots, t_n$
and $S\to P$
to be the $0$-section.
Then, we obtain a cartesian diagram
(\ref{eqPQTS}).
Since $T\to S$ is flat
and $t_1,\ldots, t_n$ is a regular
sequence,
after further replacing $Q$ by
a neighborhood of $T$,
the morphism $Q\to P$
is flat and quasi-finite
by \cite[Chapitre 0, Proposition (15.1.21) b)$\Rightarrow$ a)]{EGA4}.

2.
Since $Q$ is smooth over $S$,
there exists a commutative diagram
(\ref{eqQQQ}) such that
$Q''\to Q'$ is \'etale.

3.
Since the commutative diagram
(\ref{eqTQQ'}) defines a section
$T\to Q'\times_QT$
of a smooth morphism
$Q'\times_QT\to T$,
locally on a neighborhood
of 
$T\subset Q'\times_QT$,
there exist an integer $n\geqq0$
and an \'etale morphism
$Q'\times_QT\to {\mathbf A}^n_T$
such that
the composition
$T\to Q'\times_QT\to {\mathbf A}^n_T$
is the $0$-section.
Hence, it suffices to
take a lifting
$Q'\to {\mathbf A}^n_Q$
on a neighborhood of $T$.
\qed}
\medskip

\begin{df}\label{dfmin}
Let ${\mathcal O}_K\to {\mathcal O}_L$
be a finite morphism of
discrete valuation rings
and $E$ be the residue field of $L$.
We say that an immersion 
$T={\rm Spec}\, {\mathcal O}_L\to Q$
to a smooth scheme over
$S={\rm Spec}\, {\mathcal O}_K$ 
is {\em minimal} if the canonical surjection
$\Omega^1_{Q/S}
\otimes_{{\mathcal O}_Q}E\to 
\Omega^1_{T/S}
\otimes_{{\mathcal O}_T}E$
is an isomorphism.
\end{df}

\begin{lm}\label{lmmin}
Let ${\mathcal O}_K\to {\mathcal O}_L$
be a finite morphism of
discrete valuation rings
and set
$S={\rm Spec}\, {\mathcal O}_K$ and
$T={\rm Spec}\, {\mathcal O}_L$.

{\rm 1.}
There exists a minimal immersion
$T\to Q$ to a smooth scheme $Q$ over $S$.

{\rm 2.}
Let $T\to Q$ and $T\to Q'$
be closed immersions
to smooth schemes over $S$
and let 
\begin{equation}
\xymatrix{
T\ar[r]\ar[d]
&Q'\ar[ld]\\
Q&}
\label{eqTQQ'2}
\end{equation}
be a commutative diagram
of schemes over $S$.
If $T\to Q$
is minimal,
then $Q'\to Q$ is smooth
on a neighborhood of $T$.
\end{lm}

\proof{
1.
Let $a_1,\ldots,a_n\in {\mathcal O}_L$
be elements such that
$da_1,\ldots,da_n$
define a basis of
$\Omega^1_{{\mathcal O}_L/{\mathcal O}_K}
\otimes_{{\mathcal O}_L}E$.
Then the morphism
$T\to {\mathbf A}^n_S
={\rm Spec}\, {\mathcal O}_K
[X_1,\ldots,X_n]$
defined by $a_1,\ldots,a_n\in {\mathcal O}_L$
is unramified.
By \cite[Chapitre 4, Corollaire (18.4.7)]{EGA4},
there exist an \'etale scheme
$Q$ over ${\mathbf A}^n_S$
and an immersion $T\to Q$
lifting
$T\to {\mathbf A}^n_S$.
Then,
$\Omega^1_{Q/S}
\otimes_{{\mathcal O}_Q}E\to 
\Omega^1_{T/S}
\otimes_{{\mathcal O}_T}E$
sends a basis
$dX_1,\ldots,dX_n$
to 
$da_1,\ldots,da_n$.

2.
Since $T\to Q$ is minimal
and $\Omega^1_{Q'/S}
\otimes_{{\mathcal O}_Q}
E
\to
\Omega^1_{T/S}
\otimes_{{\mathcal O}_T}
E$ is a surjection,
the morphism
$\Omega^1_{Q/S}
\otimes_{{\mathcal O}_Q}
E
\to \Omega^1_{Q'/S}
\otimes_{{\mathcal O}_{Q'}}E$
is an injection of $E$-vector spaces.
Hence, $Q'\to Q$ is smooth
on a neighborhood of $T$.
\qed

}

\medskip

We will use in Section \ref{sskey}
the following
description of the fiber of the
conormal sheaf.

\begin{lm}\label{lmNQt}
Let ${\mathcal O}_K\to {\mathcal O}_L$ be a 
finite morphism of discrete valuation rings
and $T={\rm Spec}\, {\mathcal O}_L\to Q$
be a closed immersion to a smooth scheme
over $S={\rm Spec}\, {\mathcal O}_K$.
Let ${\mathcal I}_T\subset {\mathcal O}_Q$
be the ideal sheaf defining
the closed subscheme $T\subset Q$
and let
${\mathfrak m}_{Q,t}
\subset {\mathcal O}_{Q,t}$
and 
${\mathfrak m}_{T,t}
\subset {\mathcal O}_{T,t}$ be the maximal
ideals at the closed point $t\in T\subset Q$.

{\rm 1.}
The inclusion
${\mathcal I}_{T,t}\to {\mathfrak m}_{Q,t}$
of ideals of ${\mathcal O}_{Q,t}$
induces an isomorphism
\begin{equation}
N_{T/Q}\otimes_{{\mathcal O}_T}E
\to
{\rm Ker}({\mathfrak m}_{Q,t}/
{\mathfrak m}_{Q,t}^2
\to
{\mathfrak m}_{T,t}/
{\mathfrak m}_{T,t}^2).
\label{eqNTt1}
\end{equation}

{\rm 2.}
Let $D\subset Q$ be a smooth divisor meeting
$T$ transversally at $t$.
Then the isomorphism {\rm (\ref{eqNTt1})}
induces an isomorphism
\begin{equation}
N_{T/Q}\otimes_{{\mathcal O}_T}E
\to
N_{t/D}.
\label{eqNTt3}
\end{equation}
\end{lm}

\proof{
1.
Since the ideal ${\mathcal I}_{T,t}
\subset {\mathcal O}_{Q,t}$
of a regular local ring
is generated by a part of
a system of local parameters,
we have an exact sequence
$$
0\to
N_{T/Q}\otimes_{{\mathcal O}_T}E
\to
{\mathfrak m}_{Q,t}/
{\mathfrak m}_{Q,t}^2
\to
{\mathfrak m}_{T,t}/
{\mathfrak m}_{T,t}^2
\to 0$$
and hence the isomorphism
(\ref{eqNTt1}).

2.
By the assumption that
$D$ meets $T$ transversally,
the product 
${\mathfrak m}_{Q,t}/
{\mathfrak m}_{Q,t}^2
\to
{\mathfrak m}_{T,t}/
{\mathfrak m}_{T,t}^2
\times
{\mathfrak m}_{D,t}/
{\mathfrak m}_{D,t}^2
$
of the canonical surjections
is an isomorphism.
Since $N_{t/D}=
{\mathfrak m}_{D,t}/
{\mathfrak m}_{D,t}^2$,
the isomorphism
(\ref{eqNTt1})
induces
an isomorphism (\ref{eqNTt3}).
\qed
%

}

\begin{df}\label{dfdvr}
Let ${\mathcal O}_{K_0}\to {\mathcal O}_L$ be a morphism 
of discrete valuation rings.

{\rm 1.}
We say that $ {\mathcal O}_L$ 
is essentially of finite type
(resp.\ essentially smooth) over ${\mathcal O}_{K_0}$
if there exist a scheme $X$
of finite type
(resp.\ a smooth scheme $X$)
over ${\mathcal O}_{K_0}$,
a point $\xi\in X$
and a morphism
${\mathcal O}_{X,\xi}\to {\mathcal O}_L$
of discrete valuation rings
over ${\mathcal O}_{K_0}$
inducing an isomorphism
on the henselizations.

{\rm 2.}
Let ${\mathcal O}_L$ be a discrete
valuation ring essentially of finite
type over ${\mathcal O}_{K_0}$
and let
$T={\rm Spec}\, {\mathcal O}_L
\to Q_0$ be a morphism
to a scheme smooth over $S_0
={\rm Spec}\, {\mathcal O}_{K_0}.$
We say that
$T\to Q_0$ is 
essentially an immersion
if the following condition
is satisfied:
There exist a closed subscheme
$X\subset Q_0$,
a point $\xi\in X$
and a morphism
${\mathcal O}_{X,\xi}
\to {\mathcal O}_L$ of discrete valuation rings
inducing $T\to Q_0$ and
an isomorphism
on the henselizations.

We say that
an essential immersion
$T\to Q_0$ is divisorial
if $X\subset Q_0$ above
is a Cartier divisor.

{\rm 3.}
Let
\begin{equation}
\begin{CD}
{\rm Spec}\,
{\mathcal O}_{K'_0}=\ &
S'_0@<<< Q'_0@<<< T'&\ =
{\rm Spec}\, {\mathcal O}_{L'}\\
&@VVV@VVV@VVV\\
&S_0@<<<Q_0@<<< T
\end{CD}
\label{eqbscfT}
\end{equation}
be a commutative diagram of
schemes satisfying the following conditions:
The large rectangle
is defined by morphisms of
discrete valuation rings
and the compositions
of horizontal arrows are
essentially of finite type.
The left horizontal arrows
are smooth
and the right horizontal
arrows are essentially immersions.

Let $X\subset Q_0$ and
$X'\subset Q'_0$ be as in {\rm 2.}
Then, we say that the right square is
essentially cartesian,
if the diagram
\begin{equation}
\begin{CD}
Q'_0@<<< X'\\
@VVV@VVV\\
Q_0@<<< X
\end{CD}
\label{eqQX}
\end{equation}
is cartesian on a neighborhood
of the image of $T'$.
\end{df}

In Definition 
{\rm \ref{dfdvr}.1},
we may take $X$ to be regular.

\begin{lm}\label{lmLSS0}
Let ${\mathcal O}_{K_0}\to {\mathcal O}_K$
be a morphism of discrete
valuation rings essentially of
finite type.
Set $S={\rm Spec}\, {\mathcal O}_K$
and
$S_0={\rm Spec}\, {\mathcal O}_{K_0}$.

{\rm 1.}
Let $X$ be a scheme
of finite type,
$\xi\in X$ be
a point and
$S\to X$ be a morphism
over $S_0$
such that
${\mathcal O}_{X,\xi}\to {\mathcal O}_K$
induces an isomorphism
on henselizations as in Definition 
{\rm \ref{dfdvr}.1}.
Then, $X$ is regular on 
a neighborhood of $\xi$.

{\rm 2}.
Let $S\to Q_0$
be essentially an immersion
to a smooth scheme over $S_0$.
Let $X\subset Q_0$
be a closed subscheme
such that $S\to Q_0$
induces a morphism $S\to X$
inducing an isomorphism on
the henselizations.
Then,
the morphism
$L_{X/Q_0}\otimes_{{\mathcal O}_X}
{\mathcal O}_S\to L_{S/Q_0}$
induces a quasi-isomorphism
$L_{S/Q_0}\to 
N_{X/Q_0}\otimes_{{\mathcal O}_X}
{\mathcal O}_S[1]$.

{\rm 3.}
The cotangent complex
$L_{S/S_0}$ is of tor-amplitude
$[-1,0]$.

{\rm 4}.
If $K$ is a separable extension
of $K_0$,
the canonical morphism
$L_{S/S_0}\to 
\Omega^1_{S/S_0}[0]$
is a quasi-isomorphism.
\end{lm}

In the following,
we write
$N_{S/Q_0}
=H_1(L_{S/Q_0})$.
This is canonically isomorphic 
to the free ${\mathcal O}_K$-module
$N_{X/Q_0}\otimes_{{\mathcal O}_X}
{\mathcal O}_K$
of finite rank
by Lemma \ref{lmLSS0}.2.

\proof{
1.
Since $S\to S_0$
is flat and since
${\mathfrak m}_K\subset
{\mathcal O}_K$
is generated by a non-zero divisor,
we may assume that $X$ is flat over
$S_0$
and that $\xi$ is the generic point
of an integral effective Cartier divisor $D
\subset X$.
Further, we may assume
that the generic fiber
$X_{K_0}=X\times_S{\rm Spec}\, K_0$
is $X\sm D$ and is regular.
Since $D$ is of finite type
over the residue field
$k$ of ${\mathcal O}_{K_0}$,
we may also assume that $D$
is regular.
Hence $X$ is regular.

2.
By 1,
we may assume that $X\subset Q_0$
is regular.
Hence the immersion $X\to Q_0$
is a regular immersion and
the canonical morphism
$L_{X/Q_0}\to N_{X/Q_0}[1]$
a quasi-isomorphism
by Lemma \ref{lmNOm}.3.
Since the formation
of cotangent complexes
commutes with
\'etale morphisms and limits,
we obtain a quasi-isomorphism
$L_{S/Q_0}\to 
N_{X/Q_0}\otimes_{{\mathcal O}_X}
{\mathcal O}_S[1]$.

3.
By 1,
there exists an affine regular scheme
$X
\subset Q_0={\mathbf A}^n_{S_0}$ 
of finite type over $S_0$,
a morphism $S\to X$ over $S_0$
and $\xi\in X$ be as in
Definition \ref{dfdvr}.1.
Since the formation
of cotangent complexes
commutes with
\'etale morphisms and limits,
we may assume that
${\mathcal O}_K={\mathcal O}_{X,\xi}$.
The canonical morphisms
$L_{Q_0/S_0}\to
\Omega^1_{Q_0/S_0}[0]$
and
$L_{X/Q_0}\to
N_{X/Q_0}[1]$
are quasi-isomorphisms
by Lemma \ref{lmNOm}.2
and 3.
Hence 
by the distinguished triangle
$L_{Q_0/S_0}\otimes
_{{\mathcal O}_{Q_0}}{\mathcal O}_X
\to L_{X/S_0}\to
L_{X/Q_0}\to$,
the cotangent complex $L_{X/S_0}$
is of tor-amplitude $[-1,0]$.

4.
Let the notation be
as in the proof of 3.
If $K$ is separable over $K_0$,
we may assume that
the generic fiber $X_{K_0}=X\times_{S_0}
{\rm Spec}\, K_0$
is smooth over $K_0$.
Since $L_{X_K/K_0}\to 
\Omega^1_{X_K/K_0}[0]$
is a quasi-isomorphism
by Lemma \ref{lmNOm}.2,
the morphism
$N_{X/Q_0}\to 
\Omega^1_{Q_0/S_0}
\otimes_{{\mathcal O}_{Q_0}}{\mathcal O}_X$
of locally free ${\mathcal O}_X$-modules
is an injection and
the assertion follows.
\qed

}

\begin{lm}\label{lmcar0}
Let ${\mathcal O}_{K_0}$ be a discrete valuation 
ring and let
${\mathcal O}_L\to {\mathcal O}_{L'}$
be a morphism of ramification index $e$
of discrete valuation rings 
essentially of finite type over ${\mathcal O}_{K_0}$.
Let
\begin{equation}
\begin{CD}
Q'_0@<<< T'&
={\rm Spec}\, {\mathcal O}_{L'}\\
@VVV@VVV&\\
Q_0@<<<T&
={\rm Spec}\, {\mathcal O}_L
\end{CD}
\label{eqcar0}
\end{equation}
be a commutative diagram
of essential immersions
to smooth schemes over $S_0={\rm Spec}\, {\mathcal O}_{K_0}$.

{\rm 1.}
Let $\xi\in Q_0$ and $\xi'\in Q'_0$ 
be the images of the closed points
of $T$ and of $T'$ respectively.
Then, the following conditions
are equivalent:

{\rm (1)}
The diagram {\rm (\ref{eqcar0})}
is essentially cartesian.

{\rm (2)}
The local ring
${\mathcal O}_{Q'_0\times_{Q_0}\xi,\xi'}$
is of length $e$.

{\rm 2.}
Assume that {\rm (\ref{eqcar0})} is 
essentially cartesian.
Then, 
the morphism $Q'_0\to Q_0$
is flat on a neighborhood of 
the image of $T'$
and the canonical morphism
\begin{equation}
N_{T/Q_0}\otimes
_{{\mathcal O}_T}{\mathcal O}_{T'}
\to N_{T'/Q'_0}
\label{eqNTQ0}
\end{equation}
is an isomorphism.

{\rm 3.}
Assume that {\rm (\ref{eqcar0})} is 
essentially cartesian
and that $T'\to T$ is finite.
Then, the diagram {\rm (\ref{eqcar0})} is 
cartesian on a neighborhood of $T'$.
\end{lm}

\proof{
1.
Let $X\subset Q_0$,
$X'\subset Q'_0$
be as in Definition \ref{dfdvr}.3
and let
$\xi \in X$,
$\xi' \in X'$ be
the images of the closed points
of $T,T'$.
Let
$I\subset {\mathcal O}_{Q_0,\xi}$
and
$I'\subset {\mathcal O}_{Q'_0,\xi'}$
be the kernels
of the morphisms
${\mathcal O}_{Q_0,\xi}\to {\mathcal O}_L$
and
${\mathcal O}_{Q'_0,\xi'}\to {\mathcal O}_{L'}$.
Then, we have 
an inclusion
$I{\mathcal O}_{Q'_0,\xi'}\subset I'$ of ideals of ${\mathcal O}_{Q'_0,\xi'}$
and the condition (1)
is equivalent to the equality
$I{\mathcal O}_{Q'_0,\xi'}=I'$. 
By the surjection
${\mathcal O}_{Q'_0\times_{Q_0}\xi,\xi'}
\to 
{\mathcal O}_{X'\times_X\xi,\xi'}
={\mathcal O}_{L'}/{\mathfrak m}_L
{\mathcal O}_{L'}
={\mathcal O}_{L'}/{\mathfrak m}_{L'}^e$,
the condition (2)
is equivalent to the isomorphism
${\mathcal O}_{Q'_0\times_{Q_0}\xi,\xi'}
\to 
{\mathcal O}_{X'\times_X\xi,\xi'}$.

Let $s\in {\mathcal O}_{Q_0,\xi}$
be a lifting of a uniformizer
of ${\mathcal O}_L$
and let $\bar I\subset 
A={\mathcal O}_{Q_0,\xi}/(s)$
and 
$\bar I'
\subset A'={\mathcal O}_{Q'_0,\xi'}/(s)$
be the images of 
$I$ and $I'$. 
Let $E$ and $E'$
be the residue fields of $L$
and of $L'$.
Then, since
$E=
A/\bar I$,
${\mathcal O}_{Q'_0\times_{Q_0}\xi,\xi'}
=A'/\bar IA'$
and
${\mathcal O}_{X'\times_X\xi,\xi'}
=A'/\bar I'$,
the condition (2) is further equivalent to
the equality
$\bar IA'=\bar I'$.

Since $I$ and $I'$
are generated by regular sequences
of ${\mathcal O}_{Q_0,\xi}$
and of ${\mathcal O}_{Q'_0,\xi'}$
that remain to be regular sequences
after joining $s$,
the surjections
$I\otimes_{{\mathcal O}_{Q_0,\xi}} E\to 
\bar I\otimes_A E$
and $I'\otimes_{{\mathcal O}_{Q'_0,\xi'}} E'
\to \bar I'\otimes_{A'} E'$
are isomorphisms.
By Nakayama's lemma,
the equality
$\bar IA'=\bar I'$
is equivalent to
the surjectivity of
$\bar I\otimes_A E'
\to \bar I'\otimes_{A'} E'$
and 
the equality
$I{\mathcal O}_{Q'_0,\xi'}=I'$
is equivalent to
the surjectivity of
$I\otimes_{{\mathcal O}_{Q_0,\xi}} E'\to 
I'\otimes_{{\mathcal O}_{Q'_0,\xi'}} E'$
and the assertion follows.

2.
Let $X\subset Q_0$ and
$X'\subset Q'_0$ be
as in Definition \ref{dfdvr}.3.
Since $T'\to T$ is flat,
after replacing $Q'_0$
by a neighborhood of 
the image of $T'$,
we may assume that $X'\to X$
is flat.
Since (\ref{eqQX}) is cartesian
on a neighborhood of 
the image of $T'$,
the morphism $Q'\to Q$
is also flat
on a neighborhood of 
the image of $T'$
by \cite[Chapitre 0, Proposition (15.1.21) b)$\Rightarrow$ a)]{EGA4}.
Hence the morphism (\ref{eqNTQ0})
is an isomorphism.

3.
By 2.,
we may assume
that $Q'\to Q$ is flat
and quasi-finite.
Since the underlying
set of a quasi-finite scheme
$T\times_QQ'$ over $T$
is finite,
the normal scheme $T'$ is an open
subscheme of $T\times_QQ'$
and the assertion follows.
\qed

}

\begin{lm}\label{lmcar}
Let ${\mathcal O}_K$ be a discrete valuation 
ring essentially of finite type over 
a discrete valuation ring
${\mathcal O}_{K_0}$.
Let ${\mathcal O}_K\to
{\mathcal O}_L$ be a finite morphism
of discrete valuation rings.
We consider a commutative diagram
\begin{equation}
\begin{CD}
Q_0@<<< T&\ ={\rm Spec}\,{\mathcal O}_L\\
@VVV@VVV\\
S_0@<<<S&\ ={\rm Spec}\,{\mathcal O}_K
\end{CD}
\label{eqQ0TS0}
\end{equation}
where $Q_0$
is smooth over $S_0$
and $T\to Q_0$ is
essentially an immersion.

{\rm 1.}
The morphism
$T\to Q=Q_0\times_{S_0}S$
defined by the diagram {\rm (\ref{eqQ0TS0})}
is an immersion.
If $Q_0\to S_0$ is separated,
then $T\to Q$ is a closed immersion.

{\rm 2.}
Let {\rm (\ref{eqcar0})} be
an essentially cartesian diagram
such that $T'\to T$
is defined by a finite morphism
of discrete valuation rings.
Then, the diagram
\begin{equation}
\begin{CD}
Q'&\ =Q'_0\times_{S_0}S@<<< T'\\
@VVV&@VVV\\
Q&\ =Q_0\times_{S_0}S@<<< T
\end{CD}
\label{eqQQ0}
\end{equation}
is cartesian on a neighborhood
of the image of $T'$ and
the canonical morphism
\begin{equation}
N_{T/Q}\otimes
_{{\mathcal O}_T}{\mathcal O}_{T'}
\to N_{T'/Q'}
\label{eqNTQ'0}
\end{equation}
is an isomorphism.
After replacing $Q_0$
and $Q'_0$ by \'etale neighborhoods,
the vertical arrows 
are finite flat.
\end{lm}

\proof{
1.
Since $T$ is finite over $S$
and $T\times_S{\rm Spec}\, F
\to Q$ is an immersion
for the residue field $F$,
the morphism $T\to Q$
is an immersion
by Nakayama's lemma.
If $Q_0\to S_0$ is separated,
then $Q\to S$ is also separated
and the immersion
$T\to Q$
is a closed immersion
since $T$ is finite over $S$.

2.
By Lemma \ref{lmcar0}.3,
we may assume that
(\ref{eqcar0}) is cartesian.
Then,
(\ref{eqQQ0}) is cartesian
and 
(\ref{eqNTQ'0})
is an isomorphism.

Since $T'\to T$ is flat,
the morphism
$Q'_0\to Q_0$
is quasi-finite and flat
on a neighborhood of the image
of $T'\to Q'_0$.
Hence $Q'_0\to Q_0$
is finite and flat on
\'etale neighborhoods.
\qed

}
\medskip

\begin{pr}\label{prbsc}
Let $S_0={\rm Spec}\, {\mathcal O}_{K_0}$
be a discrete valuation ring
with perfect residue field $k$.
Let ${\mathcal O}_L$ be a discrete
valuation ring essentially of finite
type over ${\mathcal O}_{K_0}$.

{\rm 1.}
There exists a divisorial
essential immersion
$T\to Q_0$ over $S_0$,

{\rm 2.}
Let $T\to Q_0$
be a divisorial essential immersion.
Then, on a neighborhood of the image of $T$,
there exists a divisor $D_0\subset Q_0$
smooth over $S_0$
such that
$D_0$ meets $T$
transversally at the image
$\xi\in Q_0$ of the closed point of $T$
and that $\xi$ is the generic point
of $D_0\times_{S_0}k$.

{\rm 3.}
Let
\begin{equation}
\begin{CD}
S'_0@<<< Q'_0@<<< T'\\
@VVV@.@VVV\\
S_0@<<<Q_0@<<< T
\end{CD}
\label{eqbscf}
\end{equation}
be a commutative diagram of
schemes such that
the rectangle
is defined by morphisms of
discrete valuation rings.
Assume that
$T$ and $T'$ are
essentially of finite
type over $S_0$ and $S'_0$
respectively
and that the right horizontal
arrows are
divisorial essential immersions to smooth schemes.

Assume that the ramification index
$e_{K'_0/K_0}$ is $1$.
Then, after replacing $Q'_0$
by an \'etale neighborhood,
there exists a morphism
$Q'_0\to Q_0$ of schemes
such that the diagram
{\rm (\ref{eqbscf})}
with $Q'_0\to Q_0$ inserted is
commutative and its right square
is essentially cartesian.
\end{pr}

\proof{
1.
We may assume that
${\mathcal O}_L={\mathcal O}_{X,\xi}$
is a local ring of a normal flat
scheme of finite type over $S_0$.
Let $t\in {\mathcal O}_L$ be a uniformizer and let $v_1,\ldots,v_n\in {\mathcal O}_L$
be liftings of a transcendental basis of
the residue field
$F={\mathcal O}_L/{\mathfrak m}_L$ over $k$ such that
$F$ is a finite separable extension of
$k(\bar v_1,\ldots,\bar v_n)$.
After shrinking $X$ if necessary,
define a morphism
$X\to {\mathbf A}^{n+1}_{S_0}$
by $t,v_1,\ldots,v_n$.
Then, $X\to
{\mathbf A}^{n+1}_{S_0}$
is unramified at $\xi$.
Hence
by \cite[Chapitre 4, Corollaire (18.4.7)]{EGA4}
after further shrinking $X$ if necessary,
there exists a scheme $Q_0$
\'etale over
${\mathbf A}^{n+1}_{S_0}$
containing $X$ as
an effective Cartier divisor.

2.
Let $\xi\in X\subset Q_0$
be as in Definition \ref{dfdvr}.2.
The local ring
${\mathcal O}_{Q_0,\xi}$
is a regular local ring of dimension $2$.
Let $f\in {\mathcal O}_{Q_0,\xi}$
be a generator of the divisor
$X$ and $t_0\in {\mathcal O}_{K_0}$
be a uniformizer.
Then, the classes of
$f$ and of $t_0$
in ${\mathfrak m}_\xi/
{\mathfrak m}_\xi^2$ are non-zero.
Hence, there exists an
element $t\in {\mathfrak m}_\xi$
such that
${\mathfrak m}_\xi
=(f,t)=(t_0,t)$.

Let $D_0\subset Q_0$ be
the divisor defined by $t$.
Then, by ${\mathfrak m}_\xi=(t_0,t)$,
the divisor
$D_0\times_{S_0}k\subset
Q_0\times_{S_0}k$ is regular
and $\xi$ is the generic point
of $D_0\times_{S_0}k$,
on a neighborhood of $\xi$.
Since $k$ is perfect,
$D_0\times_{S_0}k$
is smooth over $k$
and $D_0\subset Q_0$
is smooth over $S_0$
on a neighborhood of $\xi$.
Further, by ${\mathfrak m}_\xi
=(f,t)$,
the divisor $D_0$ meets $T$
transversally at the closed point $\xi$.

3.
Let $\xi\in X\subset Q_0$
and $\xi'\in X'\subset Q'_0$
be as in Definition \ref{dfdvr}.2.
Let $t_0\in {\mathcal O}_{K_0}\subset {\mathcal O}_{K'_0}$
be a common uniformizer
and let $f\in {\mathfrak m}_{\xi}\subset
{\mathcal O}_{Q_0,\xi}$
and $f'\in {\mathfrak m}_{\xi'}\subset
{\mathcal O}_{Q'_0,\xi'}$ be elements
defining $X\subset Q_0$ and 
$X'\subset Q'_0$ on neighborhoods of
$\xi$ and $\xi'$ respectively.
Let $s\in {\mathfrak m}_{\xi}$
and $s'\in {\mathfrak m}_{\xi'}$ be
elements such that $(f,s)=(t_0,s)
={\mathfrak m}_{\xi}$
and
$(f',s')=(t_0,s')
={\mathfrak m}_{\xi'}$
as in the proof of 2.

Let $e=e_{L'/L}$ be the ramification index.
Since $s\in {\mathfrak m}_{\xi}$
and $s'\in {\mathfrak m}_{\xi'}$
are liftings of uniformizers of
${\mathcal O}_L$
and of
${\mathcal O}_{L'}$,
there exists a unit
$u\in {\mathcal O}_{Q'_0,\xi'}^\times$
such that 
$s=us^{\prime e}$ in ${\mathcal O}_{L'}$.
Then, since $Q_0\to S_0$
is smooth and since $T'\to Q'_0$
is essentially an immersion,
after replacing $Q'_0$ by an \'etale neighborhood
if necessary,
there exists a morphism 
$Q'_0\to Q_0$ that makes the diagram
(\ref{eqbscf}) commutative and
satisfies
$s\mapsto us^{\prime e}$.
Then since
${\mathcal O}_{Q'_0\times_{Q_0}\xi,\xi'}
={\mathcal O}_{Q'_0,\xi'}/
(f,s)
={\mathcal O}_{Q'_0,\xi'}/
(f',s^{\prime e})
={\mathcal O}_{L'}/(s^{\prime e})$
is of length $e$,
the right square is essentially cartesian
by Lemma \ref{lmcar0}.1
(2)$\Rightarrow$(1).
\qed

}

\begin{pr}\label{prImq}
Let ${\mathcal O}_{K_0}$ be a discrete
valuation ring with perfect residue field
$k$ and
let 
$$\begin{CD}
{\mathcal O}_{K'}@>>> {\mathcal O}_{L'}
\\
@AAA @AAA\\
{\mathcal O}_K@>>> {\mathcal O}_L
\end{CD}
$$ be a commutative diagram
of morphisms of
discrete valuation rings
essentially of finite type over ${\mathcal O}_{K_0}$
such that the horizontal arrows
are finite.
Let
\begin{equation}
\begin{CD}
Q'_0@<<< T'&={\rm Spec}\, {\mathcal O}_{L'}\\
@VVV@VVV\\
Q_0@<<< T&={\rm Spec}\, {\mathcal O}_L
\end{CD}
\label{eqQQTT0}
\end{equation}
be an essentially cartesian diagram
of essential immersions
to smooth separated
schemes over $S_0={\rm Spec}\, {\mathcal O}_{K_0}$.
Set $S={\rm Spec}\, {\mathcal O}_K$,
$S'={\rm Spec}\, {\mathcal O}_{K'}$
and let
$$\begin{CD}
Q'&=Q'_0\times_{S_0}S'@<<< T'\\
@VVV&@VVV\\
Q&=Q_0\times_{S_0}S@<<< T
\end{CD}$$
be the induced commutative diagram
of immersions.

{\rm 1.}
The morphism
$T'\to T\times_{Q_0}Q'_0$
defined by {\rm (\ref{eqQQTT0})}
is formally \'etale
and induces 
a formally \'etale morphism
$T'\times_SS'\to T\times_QQ'
=T\times_{Q_0}Q'_0\times_SS'$.

{\rm 2.}
Assume that
$T\to Q_0$ and
$T'\to Q'_0$ are separated
and let
${\mathcal J}\subset {\mathcal O}_{S'\times_SS'}$
be the ideal defining the diagonal
$S'\to S'\times_SS'$
and let ${\mathcal I}_T\subset {\mathcal O}_Q$
and ${\mathcal I}_{T'}\subset {\mathcal O}_{Q'}$
be the ideals defining the closed subschemes
$T\subset Q$ and $T'\subset Q'$.
Let $m\geqq 0$ and $q\geqq 2$ be
integers and set
$d'=\dim \Omega^1_{F'/k}$.
Let 
$A$ denote the quotient ring
${\mathcal O}_{K'}/{\mathfrak m}_{K'}^m
[X_0,\ldots,X_{d'}]/(X_0,\ldots,X_{d'})^q$.
Then, the following conditions
are equivalent:

{\rm (1)}
We have ${\mathcal I}_T{\mathcal O}_{Q'}
\subset {\mathfrak m}_{K'}^m
{\mathcal I}_{T'}+
{\mathcal I}_{T'}^q$.

{\rm (2)}
The closed subscheme
$Y\subset S'\times_SS'$
defined by the ideal
${\mathcal I}_Y={\mathfrak m}_{K'}^m
{\mathcal O}_{S'\times_SS'}+{\mathcal J}^q
\subset
{\mathcal O}_{S'\times_SS'}$
is isomorphic to
${\rm Spec}\, A$.
\end{pr}
\medskip

Note that the condition (2)
does not depend on $T$ or $S_0$.

\proof{
1.
Let $X\subset Q_0$ and
$X'\subset Q'_0$ be
as in Definition \ref{dfdvr}.3.
We consider the commutative diagram
$$\begin{CD}
Q'_0
@<<<
X'
@<<<
T'\\
@VVV@VVV@VVV\\
Q_0
@<<<
X
@<<<
T.
\end{CD}$$
Since (\ref{eqQQTT0})
is assumed essentially cartesian,
we may assume that
the left square is cartesian.
Since the right horizontal arrows
are formally \'etale,
the morphism
$T'\to T\times_{Q_0}Q'_0$
is formally \'etale.
Since $T\times_QQ'
=T\times_{Q_0}Q'_0\times_SS'$,
the morphism
$T'\times_SS'\to T\times_QQ'$
is formally \'etale.

2.
Since $T\to Q_0$ and
$T'\to Q'_0$ are separated,
the immersions $T\to Q$
and $T'\to Q'$ are closed immersions
by Lemma \ref{lmcar}.1.
Let $N$ be the length of
the ring $A$
and $e=e_{L'/K'}$ be 
the ramification index.
Let $Z\subset Q'$ be the
closed subscheme
defined by the ideal
${\mathcal I}_Z=
{\mathfrak m}_{K'}^m{\mathcal O}_{Q'}+
{\mathcal I}_{T'}^q
\subset {\mathcal O}_{Q'}$.
Since $Q'\to S'$ is smooth of relative dimension
$d'+1$ and since $T'\to Q'$ is a regular immersion,
the length of ${\mathcal O}_Z$ is $Ne$.
We show that,
for the fiber product
$T'\times_{S'}Y$
with respect to the restriction of
the first projection,
the both conditions
(1) and (2) are equivalent to 
the condition
that ${\mathcal O}_{T'\times_{S'}Y}$
is of length $Ne$.

Since
${\mathcal I}_T{\mathcal O}_{Q'}
\subset
{\mathcal I}_{T'}$
and $
{\mathfrak m}_{K'}^m{\mathcal I}_{T'}+
{\mathcal I}_{T'}^q
=
({\mathfrak m}_{K'}^m{\mathcal O}_{Q'}+
{\mathcal I}_{T'}^q)
\cap
{\mathcal I}_{T'}$,
the condition (1) is equivalent to
the inclusion
${\mathcal I}_T{\mathcal O}_{Q'}
\subset
{\mathfrak m}_{K'}^m{\mathcal O}_{Q'}+
{\mathcal I}_{T'}^q={\mathcal I}_Z$.
This means that
the closed immersion
$T\times_QZ\to Z$ is an isomorphism.
Since ${\mathcal O}_Z$ is of length $Ne$,
the condition (1)
is equivalent to the condition
that ${\mathcal O}_{T\times_QZ}$
is of length $Ne$.
To show that this
is equivalent to the condition
that ${\mathcal O}_{T'\times_{S'}Y}$
is of length $Ne$,
it suffices to show that
the formally \'etale morphism
$T'\times_SS'\to T\times_QQ'$
induces an isomorphism
$T'\times_{S'}Y\to T\times_QZ$
of closed subschemes.

The closed subscheme
$T'$ of
$T'\times_SS'$ and
of $T\times_QQ'$
is defined by
the pull-backs of ${\mathcal J}
\subset {\mathcal O}_{S'\times_SS'}$ and
of ${\mathcal I}_{T'}\subset {\mathcal O}_{Q'}$ respectively.
Since the morphism
$T'\times_SS'\to T\times_QQ'$
is formally \'etale by 1,
this induces an isomorphism
on the closed subschemes
defined by the pull-backs of
${\mathcal J}^q$ and
of ${\mathcal I}^q_{T'}$ respectively.
Hence by 
${\mathcal I}_Y={\mathfrak m}_{K'}^m
{\mathcal O}_{S'\times_SS'}+{\mathcal J}^q
\subset
{\mathcal O}_{S'\times_SS'}$
and ${\mathcal I}_Z=
{\mathfrak m}_{K'}^m{\mathcal O}_{Q'}+
{\mathcal I}_{T'}^q
\subset{\mathcal O}_{Q'}$,
the morphism
$T'\times_SS'\to T\times_QQ'$ induces an isomorphism
$T'\times_{S'}Y=
(T'\times_SS')
\times_{S'\times_SS'}Y
\to
T\times_QZ
=(T\times_QQ')\times_{Q'}Z$.
Thus the condition (1)
is equivalent to the condition
that ${\mathcal O}_{T'\times_{S'}Y}$
is of length $Ne$
as required.

Since
$\Omega^1_{S'/S}
=N_{S'/S'\times_SS'}$
is generated by at most
$d'+1$ elements,
there exists a closed immersion
$Y\to {\rm Spec}\, A$.
Hence the condition (2)
is equivalent to the condition
that the local ring ${\mathcal O}_Y$
is of length $N$.
Since $T'\to S'$ is faithfully flat,
this is further equivalent to the condition
that ${\mathcal O}_{T'\times_{S'}Y}$
is of length $Ne$.
\qed}

\begin{ex}\label{eg129}
{\rm Let ${\mathcal O}_K$
be a henselian discrete valuation
ring with residue field $F$
of characteristic $p>0$.
Let $t\in {\mathcal O}_K$
be a uniformizer
and $u_1,\ldots,u_n\in {\mathcal O}_K$
be a lifting of a $p$-basis of $F$.
Let $m=qm'\geqq 1$ be an
integer divisible by a power $q>1$
of $p$ and 
${\mathcal O}_{K'}$
be the henselization of
$${\mathcal O}_K[v_0,\ldots,v_n,
w_1,\ldots,w_n,t']/
(t-v_0t^{\prime m},
u_i-(w_i^m+v_it^{\prime m});
i=1,\ldots,n)$$
at the prime ideal $(t')$.

Then, $m$ and $q$ satisfies
the condition (2) in Proposition \ref{prImq}.
In fact, 
set $R={\mathcal O}_K
[T',V_0,\ldots,V_n,
W_1,\ldots,W_n]$ and
define a formally unramified morphism
$S'={\rm Spec}\, {\mathcal O}_{K'}
\to {\mathbf A}^{2n+2}_S
={\rm Spec}\, R$
by $t',v_0,\ldots,v_n,w_1,\ldots,w_n
\in {\mathcal O}_{K'}$.
Then its base change
$S'\times_SS'\to {\mathbf A}^{2n+2}_{S'}$ 
defines a formally \'etale morphism
to the closed subscheme defined by
the ideal 
$$I= (v_0t^{\prime m}\otimes 1-1\otimes V_0T^{\prime m},
(w_i^m+v_it^{\prime m})\otimes 1-
1\otimes (W_i^m+V_iT^{\prime m}); i=1,\ldots, n)
\subset R.$$
Let $$J
=(v_i-V_i; i=0,\ldots,n,
w_i-W_i; i=1,\ldots,n, t'-T')\subset R$$ 
denote the ideal
defining the diagonal section
$S'\subset S'\times_SS'\to  {\mathbf A}^{2n+2}_{S'}$.
Then, by the conguence $a^m-b^m
\equiv (a^{m'}-b^{m'})^q
\bmod p$
and $p,t^{\prime m}
\in {\mathfrak m}_K{\mathcal O}_{K'}
=
{\mathfrak m}_{K'}^m$,
we have 
$I\subset J^q+{\mathfrak m}_{K'}^mR$
and the condition (2) is satisfied.
\qed}
\end{ex}

\subsection{Differentials
and conormal modules}\label{ssTor}

\begin{lm}\label{lmTorH}
Let ${\mathcal O}_K\to {\mathcal O}_L$
be a finite morphism of discrete
valuation rings such that $L$ is
a finite separable extension of $K$
and let
$E$ be the residue field of $L$.
Set $S={\rm Spec}\, {\mathcal O}_K$ 
and $T={\rm Spec}\, {\mathcal O}_L$.

{\rm 1.}
The morphisms ${\rm Spec}\, E
\to T\to S$
define an injection
\begin{equation}
{\rm Tor}_1^{{\mathcal O}_T}
(\Omega^1_{T/S},E)
\to H_1(L_{E/S})
\label{eqTorL}
\end{equation}
of $E$-vector spaces.

{\rm 2.}
Let $T\to Q$ be an immersion
to a smooth scheme $Q$ over $S$.
Then, the distinguished triangle
$L_{Q/S}\otimes^L_{{\mathcal O}_Q}{\mathcal O}_T
\to L_{T/S}\to L_{T/Q}\to $
defines an injection
\begin{equation}
{\rm Tor}_1^{{\mathcal O}_T}
(\Omega^1_{T/S},E)
\to N_{T/Q}\otimes_{{\mathcal O}_T}E
\label{eqTorEN}
\end{equation}
of $E$-vector spaces.
If $T\to Q$ is minimal,
then
{\rm (\ref{eqTorEN})}
is an isomorphism.

{\rm 3}.
Let $K\subset M\subset L$
be an intermediate extension
and let $T\to X={\rm Spec}\, {\mathcal O}_M
\to S$ be the corresponding
morphisms.
Then, 
the associated morphism
\begin{equation}
{\rm Tor}_1^{{\mathcal O}_M}
(\Omega^1_{{\mathcal O}_M/{\mathcal O}_K},E)
\to
{\rm Tor}_1^{{\mathcal O}_L}
(\Omega^1_{{\mathcal O}_L/{\mathcal O}_K},E)
\label{eqML}
\end{equation}
is an injection.
\end{lm}

\proof{
1.
The morphisms ${\rm Spec}\, E
\to T\to S$ define 
a distinguished triangle 
\begin{equation}
\to L_{T/S}\otimes^L_{{\mathcal O}_T}E
\to L_{E/S}
\to L_{E/T}\to 
\end{equation}
of cotangent complexes.
Since $L$ is assumed to be
a separable extension of $K$,
the canonical morphism
$L_{T/S}\to \Omega^1_{T/S}[0]$
is a quasi-isomorphism by Lemma \ref{lmLSS0}.4.
By this and
the quasi-isomorphism
$L_{E/T}\to N_{E/T}[1]$,
we obtain an exact sequence
\begin{equation}
0\to {\rm Tor}_1^{{\mathcal O}_T}
(\Omega^1_{T/S},E)
\to
H_1(L_{E/S})
\to N_{E/T}\to
\Omega^1_{T/S}\otimes_{{\mathcal O}_T}E
\to \Omega^1_{E/F}\to 0
\end{equation}
of $E$-vector spaces
and the assertion follows.

2.
By Lemma \ref{lmLSS0}.4,
the canonical morphism
$L_{T/S}\to \Omega^1_{T/S}[0]$
is a quasi-isomorphism.
Hence the distinguished triangle
$L_{Q/S}\otimes^L_{{\mathcal O}_Q}{\mathcal O}_T
\to L_{T/S}\to L_{T/Q}\to $
defines an exact sequence
\begin{equation}
0\to N_{T/Q}\to \Omega^1_{Q/S}
\otimes_{{\mathcal O}_Q}{\mathcal O}_T\to 
\Omega^1_{T/S}
\to 0
\label{eqTorN}
\end{equation}
of ${\mathcal O}_T$-modules.
This defines
an exact sequence
\begin{equation}
0\to 
{\rm Tor}_1^{{\mathcal O}_T}
(\Omega^1_{T/S},E)
\to N_{T/Q}\otimes_{{\mathcal O}_T}E\to \Omega^1_{Q/S}
\otimes_{{\mathcal O}_T}E\to 
\Omega^1_{T/S}
\otimes_{{\mathcal O}_T}E\to 
0
\label{eqTor}
\end{equation}
of $E$-vector spaces.

If $T\to Q$ is minimal,
the right morphism
$\Omega^1_{Q/S}
\otimes_{{\mathcal O}_T}E\to 
\Omega^1_{T/S}
\otimes_{{\mathcal O}_T}E$
is an isomorphism
and hence the left morphism
$
{\rm Tor}_1^{{\mathcal O}_T}
(\Omega^1_{T/S},E)
\to N_{T/Q}\otimes_{{\mathcal O}_T}E$
is an isomorphism.

3.
By the distinguished triangle
$L_{X/S}\otimes_{{\mathcal O}_X}^L
{\mathcal O}_T
\to L_{T/S}\to L_{T/X}\to$
and Lemma \ref{lmLSS0}.3
applied to $L_{T/X}$,
we obtain an injection (\ref{eqML}).
\qed

}

\medskip
We make Lemma \ref{lmTorH}
explicit in the case where
${\mathcal O}_L$
is generated by one element
over ${\mathcal O}_K$.
Let ${\mathcal O}_K\to {\mathcal O}_L$
be a finite morphism of discrete valuation rings
and assume that $L\supsetneqq K$
is a separable extension.
Assume that the residue field $E$
of $L$ is a purely inseparable extension of
degree $q$ of
the residue field $F$ of $K$
and let $e=e_{L/K}$ be the ramification
index.
Assume that ${\mathcal O}_L$
is generated by one element $\alpha$
over ${\mathcal O}_K$
and let $f\in {\mathcal O}_K[X]$
be the minimal polynomial of $\alpha$.
Since $\bar f\in F[X]$
is a power of the minimal polynomial
of $\bar \alpha \in E=F(\bar \alpha)$,
there exists $\bar u\in F$
such that
$\bar f\equiv (X^q-\bar u)^e$
and that $X^q-\bar u
\in F[X]$ is irreducible.
If $q=1$ and $\bar u=0$,
then $f$ is an Eisenstein polynomial.

Since the immersion
$T={\rm Spec}\, {\mathcal O}_L
\to
Q={\rm Spec}\, {\mathcal O}_K[X]$
defined by $\alpha$ is minimal,
the canonical morphism
${\rm Tor}_1^{{\mathcal O}_T}
(\Omega^1_{T/S},E)
\to N_{T/Q}\otimes_{{\mathcal O}_T}E$
(\ref{eqTorEN})
is an isomorphism by Lemma \ref{lmTorH}.2.

\begin{lm}\label{lmNH}
Let the notation be as above
and let 
\begin{equation}
N_{T/Q}\otimes_{{\mathcal O}_T}E
\to
H_1(L_{E/S})
\label{eqNL}
\end{equation}
be the composition of the inverse
of the isomorphism {\rm (\ref{eqTorEN})}
with the injection
${\rm Tor}^{{\mathcal O}_L}_1
(\Omega^1_{{\mathcal O}_L/{\mathcal O}_K},
E)
\to H_1(L_{E/S})$
{\rm (\ref{eqTorL})}.

{\rm 1.}
Assume that $f
\in {\mathcal O}_K[X]$
is an Eisenstein polynomial
and
let $\pi=f(0)$
be the constant term.
Then the image of 
the basis $f\in 
N_{T/Q}\otimes_{{\mathcal O}_T}E$
by {\rm (\ref{eqNL})} is
$d\pi$.

{\rm 2.}
Assume that $q>1$
and set $f=(X^q-u)^e+t\cdot h$
for a uniformizer $t\in {\mathcal O}_K$,
a lifting $u\in {\mathcal O}_K$ of $\bar u$
and $h\in {\mathcal O}_K[X]$
such that $\deg h<qe$.
Then,
the image of $f$
by {\rm (\ref{eqNL})} is
$-\widetilde d u+\bar h\cdot
dt$ if $e=1$
and
$\bar h\cdot
dt$ if $e>1$.
\end{lm}

\proof{
1.
This follows from $f\equiv \pi
\bmod (X,\pi)^2$.

2.
By the commutative diagram
\begin{equation}
\begin{CD}
H_1(L_{E/S})
@>>> N_{E/Q}\\
@AAA@AAA\\
{\rm Tor}^{{\mathcal O}_L}_1
(\Omega^1_{{\mathcal O}_L/{\mathcal O}_K},
E)
@>{\simeq}>> N_{T/Q}\otimes
_{{\mathcal O}_T}E
\end{CD}
\label{eqcdTL}
\end{equation}
of injections,
it suffices to show that
$X^q-u\in N_{E/Q}$
is the image of
$-\widetilde d u\in H_1(L_{E/S})$.
We consider the commutative diagram
\begin{equation}
\begin{CD}
{\rm Spec}\, E@>>> Q&\, ={\rm Spec}\,
{\mathcal O}_K[X]\\
@VVV@VVV\\
{\rm Spec}\, F(u^{1/p})
@>>> {\mathbf A}^1_S&\,
={\rm Spec}\, {\mathcal O}_K[V]
\end{CD}
\label{eqEFup}
\end{equation}
where the lower line is defined
by $V\mapsto u^{1/p}$
and the right vertical arrow
is defined by
$V\mapsto X^{q/p}$.
Since $X^q-u$ is the image
of the basis $V^p-u$
of $N_{F(u^{1/p})/
{\mathbf A}^1_S}$,
the assertion follows.
\qed

}

\begin{ex}{\rm 
Let $K=k(u)((t))$ be the field
of formal power series of characteristic $2$
and let $L$ be the Artin-Schreier
extension of $K$ defined by
$X^2-X=u/t^2$.
Then the integer ring
${\mathcal O}_L$ is generated by
$\alpha=tX$ over ${\mathcal O}_K$
and its minimal polynomial is
$f=T^2-tT-u$.
Hence the class of
$f$ is $-\widetilde d u-\alpha\cdot 
\widetilde dt$.}
\end{ex}

\medskip

Let ${\mathcal O}_K\to {\mathcal O}_L$
be a finite morphism
of discrete valuation rings
essentially of finite type over
${\mathcal O}_{K_0}$
such that the
extensions of discrete valuation fields
are separable.
We consider a commutative diagram
\begin{equation}
\begin{CD}
Q_0@<<< T&\ ={\rm Spec}\,{\mathcal O}_L\\
@VVV@VVV\\
S_0@<<<S&\ ={\rm Spec}\,{\mathcal O}_K
\end{CD}
\label{eqQ0TS}
\end{equation}
where $Q_0$
is smooth over $S_0$
and $T\to Q_0$ is
essentially an immersion.
We consider the distinguished triangle
$L_{Q/Q_0}\otimes^L
_{{\mathcal O}_Q}{\mathcal O}_T
\to L_{T/Q_0}\to L_{T/Q}\to $
of cotangent complexes
associated to morphisms
$T\to Q=Q_0\times_{S_0}S\to Q_0$.
This defines an exact sequence
\begin{equation}
\begin{CD}
0\to  N_{T/Q_0}\to 
N_{T/Q}\to \Omega^1_{S/S_0}
\otimes_{{\mathcal O}_S}{\mathcal O}_T
\to 0
\end{CD}
\label{eqNTQ}
\end{equation}
of ${\mathcal O}_T$-modules.
In fact, $L_{S/S_0}
\to \Omega^1_{S/S_0}[0]$ 
is a quasi-isomorphism
by the assumption
that $K$ is separable over $K_0$
and by Lemma \ref{lmLSS0}.4
and 
$L_{S/Q_0}
\to N_{S/Q_0}[1]$ 
is a quasi-isomorphism
by Lemma \ref{lmLSS0}.2.

If $T=S$,
then the immersion $S\to Q=Q_0\times_{S_0}S$
is a section of a smooth morphism
$Q=Q_0\times_{S_0}S\to S$
and 
the exact sequence (\ref{eqNTQ})
is identified with
\begin{equation}
\begin{CD}
0\to  N_{S/Q_0}\to 
\Omega^1_{Q_0/S_0}
\otimes_{{\mathcal O}_{Q_0}}{\mathcal O}_S
\to 
\Omega^1_{S/S_0}
\to 0
\end{CD}
\label{eqNSQ}
\end{equation}
by the isomorphism
$d\colon N_{S/Q}\to\Omega^1_{Q_0/S_0}
\otimes_{{\mathcal O}_{Q_0}}{\mathcal O}_S$.

\begin{pr}\label{prTOm}
Let ${\mathcal O}_K
\to {\mathcal O}_L$ be a finite
morphism of discrete valuation rings
essentially of finite type
over a discrete valuation
ring ${\mathcal O}_{K_0}$
with perfect residue field $k$
of characteristic $p>0$.
We assume that the
extensions of discrete valuation fields
are separable.
Set $S={\rm Spec}\, {\mathcal O}_K$,
$T={\rm Spec}\, {\mathcal O}_L$ and
$S_0={\rm Spec}\, {\mathcal O}_{K_0}$.
Let $Q_0$ be a smooth scheme
over $S_0$
and $T\to Q_0$
be essentially an immersion.
Let $T\to Q=Q_0\times_{S_0}S$
be the induced immersion.
Then we have a commutative
diagram
\begin{equation}
\xymatrix{
{\rm Tor}^{{\mathcal O}_L}_1
(\Omega^1_{{\mathcal O}_L/{\mathcal O}_K}
,E)
\ar[rrd]\ar[d]_{\rm (\ref{eqTorEN})}
\ar[rr]^{\rm (\ref{eqTorL})}&&
H_1(L_{E/S})\ar[d]^
{\rm (\ref{eqH1LOm})}\\
N_{T/Q}
\otimes_{{\mathcal O}_L}E
\ar[rr]^{\rm (\ref{eqNTQ})}&&
\Omega^1_{S/S_0}
\otimes_{{\mathcal O}_K}E
}
\label{eqTOm}
\end{equation}
where
the slant arrow is
minus of the boundary mapping
defined by 
the exact sequence
$0\to \Omega^1_{S/S_0}
\otimes_{{\mathcal O}_S}{\mathcal O}_T
\to \Omega^1_{T/S_0}
\to \Omega^1_{T/S}\to 0$.
\end{pr}

\proof{
The morphisms 
${\rm Spec}\, E\to T\to S\to S_0$ of schemes
define a commutative diagram
of distinguished triangles
\begin{equation}
\begin{CD}
@>>> L_{S/S_0}\otimes^L_{{\mathcal O}_S}E
@>>> L_{T/S_0}\otimes^L_{{\mathcal O}_S}E
@>>> L_{T/S}\otimes^L_{{\mathcal O}_T}E@>>>\\
@.@|@VVV@VVV@.\\
@>>> L_{S/S_0}\otimes^L_{{\mathcal O}_S}E
@>>> L_{E/S_0}
@>>> L_{E/S}@>>>
\end{CD}
\end{equation}
and the upper right triangle 
in (\ref{eqTOm}) is commutative.

The commutative diagram of schemes
$$
\xymatrix{
T\ar[r]\ar[rd]& S\ar[r]&S_0\\
& Q\ar[r]\ar[u]&Q_0\ar[u]
}
$$
defines a commutative diagram
\begin{equation}
\xymatrix{
L_{S/S_0}\otimes_{{\mathcal O}_S}
{\mathcal O}_T\ar[r]\ar[d]^{\simeq}
& 
L_{T/S_0}
\ar[r]\ar[d]
&
L_{T/S}\ar[d]
\ar[r]&
\\
L_{Q/Q_0}\otimes^L_{{\mathcal O}_Q}
{\mathcal O}_T\ar[r]
& L_{T/Q_0}
\ar[r]&
L_{T/Q}\ar[r]&
}
\label{eqLLT}
\end{equation}
of distinguished triangles
of complexes
of ${\mathcal O}_T$-modules.
The canonical morphisms
$L_{T/S}\to \Omega^1_{T/S}[0],
L_{S/S_0}\to \Omega^1_{S/S_0}[0]$
are quasi-isomorphisms
by Lemma \ref{lmLSS0}.4
and 
$L_{T/Q}\to N_{T/Q}[1]$
is a quasi-isomorphism
by Lemma \ref{lmNOm}.3.
Hence the diagram
(\ref{eqLLT})
with $\otimes_{{\mathcal O}_T}^LE$
shows that 
the lower left triangle in (\ref{eqTOm}) 
is commutative
by Lemma \ref{lmNOm}.1.
\qed

}

\medskip

The exact sequence (\ref{eqNTQ})
has the following functoriality.
Let
\begin{equation}
\begin{CD}
T'@>>> S'@>>> S'_0\\
@VVV@VVV@VVV\\
T@>>> S@>>> S_0,
\end{CD}
\qquad
\begin{CD}
T'@>>> Q'_0\\
@VVV@VVV\\
T@>>> Q_0
\end{CD}
\label{eqQTS}
\end{equation}
be commutative diagrams
of schemes 
such that the left one
is defined
by morphisms of
discrete valuation rings
and the horizontal arrows in the right one
are essentially regular immersions
to smooth schemes.
We assume that the
extensions of discrete valuation fields
are separable.
Then by the commutative diagram
\begin{equation}
\begin{CD}
L_{S'/S'_0}
\otimes_{{\mathcal O}_{S'}}{\mathcal O}_{T'}
@>>>L_{T'/Q'_0}@>>>
L_{T'/Q'}@>>>
 \\
@AAA@AAA@AAA@.\\
L_{S/S_0}
\otimes_{{\mathcal O}_S}{\mathcal O}_{T'}
@>>> 
L_{T/Q_0}
\otimes_{{\mathcal O}_T}{\mathcal O}_{T'}
@>>>
L_{T/Q}
\otimes_{{\mathcal O}_T}{\mathcal O}_{T'}
@>>>
\end{CD}
\end{equation}
of distinguished triangles
and Lemma \ref{lmLSS0}.4,
we have a commutative diagram
\begin{equation}
\begin{CD}
0@>>>  N_{T'/Q'_0}@>>>
N_{T'/Q'}@>d>>
 \Omega^1_{S'/S'_0}
\otimes_{{\mathcal O}_{S'}}{\mathcal O}_{T'}
@>>> 0\\
@.@AAA@AAA@AAA@.\\
0@>>>  N_{T/Q_0}
\otimes_{{\mathcal O}_T}{\mathcal O}_{T'}
@>>>
N_{T/Q}
\otimes_{{\mathcal O}_T}{\mathcal O}_{T'}
@>d>>
 \Omega^1_{S/S_0}
\otimes_{{\mathcal O}_S}{\mathcal O}_{T'}
@>>> 0
\end{CD}
\label{eqNTTQ}
\end{equation}
of ${\mathcal O}_{T'}$-modules.
If $S_0=S'_0, S=S'$ and if
the right diagram in
(\ref{eqQTS})
is essentially cartesian,
the vertical arrows are isomorphisms
by (\ref{eqNTQ0}) and
(\ref{eqNTQ'0}).
Further if $T=S$,
the diagram (\ref{eqNTTQ}) 
and the exact sequence
(\ref{eqNSQ}) give
an isomorphism
\begin{equation}
\begin{CD}
0@>>>  N_{T'/Q'_0}@>>>
N_{T'/Q'}@>>>
 \Omega^1_{S/S_0}
\otimes_{{\mathcal O}_{S}}{\mathcal O}_{T'}
@>>> 0\\
@.@AAA@AAA@|@.\\
0@>>> N_{S/Q_0}
\otimes_{{\mathcal O}_S}
{\mathcal O}_{T'}
@>>> \Omega^1_{Q_0/S_0}
\otimes_{{\mathcal O}_{Q_0}}
{\mathcal O}_{T'}
@>>>
\Omega^1_{S/S_0}\otimes_{{\mathcal O}_S}
{\mathcal O}_{T'}@>>>0
\end{CD}
\label{eqNTS}
\end{equation}
of exact sequences.

\begin{pr}\label{prTorTS}
Let ${\mathcal O}_{K_0}$ be a discrete valuation
ring with perfect residue field
and let 
$$\begin{CD}
{\mathcal O}_{K'}@>>> {\mathcal O}_{L'}\\
@AAA@AAA\\
{\mathcal O}_K@>>> {\mathcal O}_L
\end{CD}$$
be a commutative diagram of
morphisms of discrete valuation rings
essentially of finite type over ${\mathcal O}_{K_0}$.
Assume that the horizontal arrows are finite
and that the extensions of
discrete valuation fields are separable.
Let $S={\rm Spec}\, {\mathcal O}_K$,
$S'={\rm Spec}\, {\mathcal O}_{K'}$
and
$S_0={\rm Spec}\, {\mathcal O}_{K_0}$.
Let
\begin{equation}
\begin{CD}
Q'_0@<<< T'\,&={\rm Spec}\, {\mathcal O}_{L'}\\
@VVV@VVV&\\
Q_0@<<< T\,&={\rm Spec}\, {\mathcal O}_L
\end{CD}
\label{eqbscN1}
\end{equation}
be an essentially cartesian diagram
of divisorial essential immersions
and let $T\to Q=Q_0\times_{S_0}S$
and $T'\to Q'=Q'_0\times_{S_0}S'$
be the induced immersions.
Let $F,F',E$ and $E'$
be residue fields of $K,K',L$
and $L'$ and set $d=\dim_F \Omega^1_{F/k}$,
$d'=\dim_{F'}\Omega^1_{F'/k}$.

{\rm 1.}
The commutative diagram
{\rm (\ref{eqNTTQ})} defines
an exact sequence
\begin{equation}
\begin{CD}
0@>>>
N_{T/Q}\otimes_{{\mathcal O}_T}{\mathcal O}_{T'}
@>>>
N_{T'/Q'}
@>>>
\Omega^1_{S'/S}
\otimes_{{\mathcal O}_{S'}}{\mathcal O}_{T'}
@>>>0.
\end{CD}
\label{eqTorTS}
\end{equation}

{\rm 2.}
Assume 
$\dim_{F'}
\Omega^1_{S'/S}
\otimes_{{\mathcal O}_{S'}}F'=d'+1$.
Then,
$\dim_{E'}
{\rm Tor}_1^{{\mathcal O}_{S'}}
(\Omega^1_{S'/S},E')=d+1$.
The morphism
\begin{equation}
{\rm Tor}_1^{{\mathcal O}_{S'}}
(\Omega^1_{S'/S},E')
\to H_1(L_{E'/S})
\label{eqTorSE}
\end{equation}
defined by the distinguished triangle
$L_{S'/S}\otimes^L_{{\mathcal O}_{S'}}E'
\to
L_{E'/S}\to L_{E'/S'}\to$
is an isomorphism
and the morphism
\begin{equation}
H_1(L_{E'/S'})
\to
\Omega^1_{S'/S}\otimes_{{\mathcal O}_{S'}}
E'
\label{eqTorS'E}
\end{equation}
is an injection.
The injection
\begin{equation}
{\rm Tor}_1^{{\mathcal O}_{S'}}
(\Omega^1_{S'/S},E')
\to N_{T/Q}\otimes_{{\mathcal O}_T}E'
\label{eqTorSN}
\end{equation}
defined by {\rm (\ref{eqTorTS})}
is an isomorphism.

{\rm 3.}
We keep assuming
$\dim_{F'}
\Omega^1_{S'/S}
\otimes_{{\mathcal O}_{S'}}F'=d'+1$.
Then, the diagram
\begin{equation}
\begin{CD}
{\rm Tor}_1^{{\mathcal O}_{S'}}
(\Omega^1_{S'/S},E')
@>{\simeq}>{\rm (\ref{eqTorSE})}>
H_1(L_{E'/S})\\
@V{\rm (\ref{eqTorSN})}V{\simeq}V
@AA{\rm (\ref{eqTorL})}A\\
N_{T/Q}\otimes_{{\mathcal O}_T}E'
@<{\rm (\ref{eqTorEN})}<<
{\rm Tor}_1^{{\mathcal O}_T}
(\Omega^1_{T/S},E')
\end{CD}
\label{eqTorTSE}
\end{equation}
is commutative.

{\rm 4}.
For an integer
$m\geqq 0$,
the following conditions are equivalent:

{\rm (1)}
The image of the canonical morphism
$N_{T/Q}\to N_{T'/Q'}$
is a submodule of
${\mathfrak m}_{K'}^m
N_{T'/Q'}$
and 
the morphism
$N_{T/Q}\otimes_{{\mathcal O}_L}E'
\to {\mathfrak m}_{K'}^mN_{T'/Q'}\otimes_{{\mathcal O}_L'}E'$
induced by
$N_{T/Q}\to {\mathfrak m}_{K'}^mN_{T'/Q'}$
is an injection.

{\rm (2)}
The ${\mathcal O}_{K'}$-module
$\Omega^1_{S'/S}$
is isomorphic to
${\mathcal O}_{K'}^{\oplus (d'-d)}
\oplus
({\mathcal O}_{K'}/{\mathfrak m}_{K'}^m)^{\oplus (d+1)}$.

\noindent
If these equivalent conditions
are satisfied,
we have a canonical isomorphism
\begin{equation}
{\rm Tor}_1^{{\mathcal O}_{K'}}
(\Omega^1_{S'/S},F')
\to 
\Omega^1_{S'/S,{\rm tors}}
\otimes_{{\mathcal O}_{K'}}
{\mathfrak m}_{K'}^m/
{\mathfrak m}_{K'}^{m+1}.
\label{eqTorm}
\end{equation}
\end{pr}

Note that the condition (2) in 4.\
is independent of $T$.

\proof{
1.
Since the diagram (\ref{eqbscN1})
is assumed to be essentially cartesian,
the left vertical arrow in
(\ref{eqNTTQ}) is an isomorphism
by Lemma \ref{lmcar0}.2.
Hence (\ref{eqNTTQ})
induces an isomorphism
\begin{align*}
{\rm Coker}
(N_{T/Q}
\otimes_{{\mathcal O}_T}{\mathcal O}_{T'}
\to
N_{T'/Q'})
\to
&\,
{\rm Coker}
(\Omega^1_{S/S_0}
\otimes_{{\mathcal O}_S}{\mathcal O}_{T'}
\to 
\Omega^1_{S'/S'_0}
\otimes_{{\mathcal O}_{S'}}{\mathcal O}_{T'})
\\&
=
\Omega^1_{S'/S}
\otimes_{{\mathcal O}_{S'}}{\mathcal O}_{T'}
\end{align*}
and we obtain the exact sequence (\ref{eqTorTS}).

2.
Let $S\to X$ and $S'\to X'$ be 
morphisms to irreducible 
regular schemes of finite type over $S_0$
inducing isomorphisms on the henselizations
and let $X'\to X$ be a morphism
compatible with $S'\to S$.
We may assume that
the generic fibers
$X_{K_0}$ and $X'_{K_0}$
are smooth over $K_0$ and
that the morphism
$X'_{K_0}\to X_{K_0}$
is smooth.
Then, we have $\dim X=d+1$,
$\dim X'=d'+1$
and hence ${\rm dim}_{K'}
\Omega^1_{{\mathcal O}_{K'}/
{\mathcal O}_K}\otimes_{{\mathcal O}_{K'}}
K'=d'-d$.
Thus, the assumption $\dim_{F'}
\Omega^1_{S'/S}
\otimes_{{\mathcal O}_{S'}}F'=d'+1$
implies
$\dim_{E'}
{\rm Tor}_1^{{\mathcal O}_{S'}}
(\Omega^1_{S'/S},E')=d+1$.

Since $L_{E'/S'}$ is acyclic except
at degree $[-1,0]$,
the distinguished triangle
$L_{S'/S}\otimes^L_{{\mathcal O}_{S'}}E'
\to
L_{E'/S}\to L_{E'/S'}\to$ defines
an exact sequence
$0\to {\rm Tor}_1^{{\mathcal O}_{S'}}
(\Omega^1_{S'/S},E')
\to H_1(L_{E'/S})\to 
H_1(L_{E'/S'})
\to
\Omega^1_{S'/S}\otimes_{{\mathcal O}_{S'}}
E'$.
By the exact sequence (\ref{eqH1L}),
we have
$\dim_{E'}H_1(L_{E'/S})
\leqq d+1$.
Hence
the injection (\ref{eqTorSE}) is 
an isomorphism and the morphism
(\ref{eqTorS'E}) is an injection.
Since $Q_0\to S_0$
is of relative dimension $d+1$,
the free ${\mathcal O}_T$-module
$N_{T/Q}$ is of rank $d+1$.
Hence the injection 
(\ref{eqTorSN}) induced by (\ref{eqTorTS})
is an isomorphism.

3.
By the functoriality of (\ref{eqTorL}),
the diagram
\begin{equation}
\xymatrix{
{\rm Tor}_1^{{\mathcal O}_{S'}}
(\Omega^1_{S'/S},E')
\ar[r]^-\simeq_-{\rm (\ref{eqTorSE})}\ar[d]
&
H_1(L_{E'/S})\\
{\rm Tor}_1^{{\mathcal O}_{T'}}
(\Omega^1_{T'/S},E')\ar[ru]
&
{\rm Tor}_1^{{\mathcal O}_T}
(\Omega^1_{T/S},E')
\ar[l]\ar[u]_{\rm (\ref{eqTorL})}
}
\label{eqTorA}
\end{equation}
is commutative.
Since 
the slant arrow
${\rm Tor}_1^{{\mathcal O}_{T'}}
(\Omega^1_{T'/S},E')\to
H_1(L_{E'/S})$ is an injection
by Lemma \ref{lmTorH}.1,
the left vertical arrow is an isomorphism.

Since the construction of (\ref{eqTorTS})
is functorial in $S'$,
we obtain a commutative diagram
\begin{equation}
\begin{CD}
0@>>>
N_{T/Q}\otimes_{{\mathcal O}_T}{\mathcal O}_{T'}
@>>>
N_{T'/Q'}
@>>>
\Omega^1_{S'/S}
\otimes_{{\mathcal O}_{S'}}{\mathcal O}_{T'}
@>>>0
\\
@.@|@VVV@VVV@.\\
0@>>>
N_{T/Q}\otimes_{{\mathcal O}_T}{\mathcal O}_{T'}
@>>>
N_{T'/(Q'_0\otimes_{S_0}T')}
@>>>
\Omega^1_{T'/S}
@>>>0
\\
@.@|@AAA@AAA@.\\
0@>>>  
N_{T/Q}\otimes_{{\mathcal O}_T}{\mathcal O}_{T'}
@>>>
N_{T/(Q_0\times_{S_0}T)}
\otimes_{{\mathcal O}_T}{\mathcal O}_{T'}
@>>>
 \Omega^1_{T/S}
\otimes_{{\mathcal O}_T}{\mathcal O}_{T'}
@>>> 0
\end{CD}
\end{equation}
of exact sequences of
${\mathcal O}_{T'}$-modules.
Since the middle term in the lower line
is canonically identified with
$\Omega^1_{Q/S}\otimes_{{\mathcal O}_Q}
{\mathcal O}_{T'}$,
the diagram
\begin{equation}
\xymatrix{
{\rm Tor}_1^{{\mathcal O}_{S'}}
(\Omega^1_{S'/S},E')
\ar[r]^-\simeq
\ar[d]^\simeq_{\rm (\ref{eqTorSN})}
&
{\rm Tor}_1^{{\mathcal O}_{T'}}
(\Omega^1_{T'/S},E')\ar[ld]
\\
N_{T/Q}\otimes_{{\mathcal O}_T}E'
&
{\rm Tor}_1^{{\mathcal O}_T}
(\Omega^1_{T/S},E')
\ar[u]\ar[l]_-{\rm (\ref{eqTorEN})}
}
\label{eqTorB}
\end{equation}
is commutative.
By combining (\ref{eqTorA})
and (\ref{eqTorB}), we obtain
(\ref{eqTorTSE}).

4.
The exact sequence
(\ref{eqTorTS}) induces an isomorphism
\begin{equation}
M=
{\rm Coker}(
N_{T/Q}
\otimes_{{\mathcal O}_T}
{\mathcal O}_{T'}
\to
N_{T'/Q'})
\to 
\Omega^1_{S'/S}
\otimes_{{\mathcal O}_{S'}}
{\mathcal O}_{T'}.
\label{eqNOmt}
\end{equation}
The condition (1) is equivalent to
that $M/{\mathfrak m}_{K'}^mM$
is a free ${\mathcal O}_{L'}/{\mathfrak m}_{K'}^m{\mathcal O}_{L'}$-module of rank $d'+1$
and that
${\mathfrak m}_{K'}^mM$
is a free ${\mathcal O}_{L'}$-module of rank $d'-d$.
Hence the conditions (1) and (2)
are equivalent to each other.

Since $M_{\rm tors}=
N_{T/Q}
\otimes_{{\mathcal O}_T}
({\mathfrak m}_{K'}^{-m}/{\mathcal O}_{T'})$,
the isomorphism
(\ref{eqNOmt}) induces an isomorphism
(\ref{eqTorm}).
\qed

}

\begin{ex}\label{eg135}
{\rm Let ${\mathcal O}_K
\to {\mathcal O}_{K'}$ be a
morphism of henselian discrete valuation rings
in Example {\rm \ref{eg129}}.
Then $\dim \Omega^1_{F/k}=
n$, $\Omega^1_{F'/k}=2n+1$
and
the ${\mathcal O}_{K'}$-module
$\Omega^1_{{\mathcal O}_{K'}/{\mathcal O}_K}$
is isomorphic to
${\mathcal O}_{K'}^{\oplus n+1}
\oplus
({\mathcal O}_{K'}/{\mathfrak m}_K{\mathcal O}_{K'})^{\oplus n+1}$.
\qed}
\end{ex}


\section{Smooth group schemes}\label{sGH}

We first study extensions of
vector spaces by finite groups
in Section \ref{ssExt}.
The main goal here is a criterion
(Proposition \ref{prlin})
for a finite \'etale covering
of a vector space to be an 
extension by a finite group.
This will be used in the reduction
to the perfect residue field case
in the first proof of Theorem \ref{thmgr}.

In Section \ref{ssGH},
we give a criterion 
(Proposition \ref{prisog})
for a finite \'etale covering of
a smooth group scheme
to be a morphism
of group schemes.
This will be used to
prove the crucial case
where the index $r>1$
is an integer 
in the second proof of Theorem \ref{thmgr}.

\subsection{Additive torsors over vector spaces}\label{ssExt}

Let $k$ be a field of characteristic $p>0$
and let $E$ be a $k$-vector space
of finite dimension. 
We consider $E={\rm Spec}\,
S(E^\vee)$ as a smooth algebraic group over $k$
defined by the symmetric algebra over $k$
of the dual space $E^\vee={\rm Hom}_k(E,k)$.

Let ${\rm Ext}(E,{\mathbf F}_p)$
denote the group of isomorphism
classes of extensions
$0\to {\mathbf F}_p
\to H\to E\to 0$
of smooth commutative group schemes
over $k$.
We define a morphism
\begin{equation}
E^\vee
\to
{\rm Ext}(E,{\mathbf F}_p)
\label{eqExt2}
\end{equation}
by sending a linear form
$f\colon E\to {\mathbf G}_a$
to the class $[H]$ of the extension 
defined by the commutative diagram
\begin{equation}
\begin{CD}
0@>>>
{\mathbf F}_p
@>>>
H@>>> E@>>>0\\
@.@|@VVV@VVfV@.\\
0@>>>
{\mathbf F}_p
@>>>
{\mathbf G}_a
@>{F-1}>>
{\mathbf G}_a
@>>>0
\end{CD}
\label{eqExtAS}
\end{equation}
of extensions,
where $F\colon {\mathbf G}_a
\to {\mathbf G}_a$
denotes the Frobenius
endomorphism.

Since an extension of
$E$ by ${\mathbf F}_p$
is an ${\mathbf F}_p$-torsor
over $E$,
by forgetting the group structure,
we obtain a canonical morphism
\begin{equation}
{\rm Ext}(E,{\mathbf F}_p)
\to
H^1(E,{\mathbf F}_p).
\label{eqExt}
\end{equation}
Let $F\colon S(E^\vee)\to S(E^\vee)$
denote the absolute
Frobenius endomorphism of the symmetric algebra
and identify 
\begin{equation}
H^1(E,{\mathbf F}_p)
={\rm Coker}(F-1\colon 
S(E^\vee)\to S(E^\vee))
\label{eqH1E}
\end{equation}
by the Artin-Schreier sequence
$0\to {\mathbf F}_p\to {\mathbf G}_a
\overset{F-1}\longrightarrow {\mathbf G}_a\to 0$.
For an integer $m\geqq 1$
prime to $p$,
we set
\begin{equation}
S^{mp^\infty}(E^\vee)
=
\varinjlim_{e\in {\mathbf N}}
S^{mp^e}(E^\vee)
\label{eqSinfE}
\end{equation}
with respect to $F$.

\begin{pr}\label{prExt}
Let $k$ be a field
of characteristic $p>0$
and let $E$ be a $k$-vector space
of finite dimension.

{\rm 1.}
For an integer $n\geqq 1$,
the mapping
$F\colon S^n(E^\vee)\to S^{pn}(E^\vee)$
is an injection.

{\rm 2.}
The morphism
$S(E^\vee)\to H^1(E,{\mathbf F}_p)$
{\rm (\ref{eqH1E})}
induces an isomorphism
\begin{equation}
H^1(k,{\mathbf F}_p)
\oplus
\bigoplus_{p\nmid m}
S^{mp^\infty}(E^\vee)
\to
H^1(E,{\mathbf F}_p).
\label{eqHES}
\end{equation}

{\rm 3.}
The morphism
{\rm (\ref{eqExt})} is an injection
and the morphism
{\rm (\ref{eqExt2})}
is an isomorphism.
The diagram
\begin{equation}
\begin{CD}
S(E^\vee)
@>{\rm (\ref{eqH1E})}>> H^1(E,{\mathbf F}_p)\\
@A{\cup}AA@AA{\rm (\ref{eqExt})}A\\
E^\vee
@>{\rm (\ref{eqExt2})}>>
{\rm Ext}(E,{\mathbf F}_p)
\end{CD}
\label{eqExtAS2}
\end{equation}
is commutative.
\end{pr}

\proof{
1.
Since $S(E^\vee)$ is an integral
domain,
the endomorphism
$F\colon S(E^\vee)\to S(E^\vee)$
is an injection.

2.
We identify
$H^1(E,{\mathbf F}_p)
=
{\rm Coker}(F-1\colon S(E^\vee)\to S(E^\vee))$.
Then the right hand side
is the direct sum of
$H^1(k,{\mathbf F}_p)=
{\rm Coker}(F-1\colon k\to k)$
and 
${\rm Coker}(F-1\colon \bigoplus_{e=0}^\infty S^{p^em}(E^\vee)\to 
\bigoplus_{e=0}^\infty S^{p^em}(E^\vee))
=S^{mp^\infty}(E^\vee)$
for $m\geqq 1$ prime to $p$.

3.
We give a proof for convenience.
By the definitions of
(\ref{eqH1E})
and (\ref{eqExt2}),
the diagram is commutative.
Since the composition through
the upper left is an injection
by 1.\ and 2.,
it suffices to show the surjectivity of
{\rm (\ref{eqExt2})}.
Since $E$ is a direct sum of
copies of ${\mathbf G}_a$,
we may assume $E={\mathbf G}_a$.

Since an automorphism of hyperbolic curve
has no non-trivial deformation,
a connected affine smooth group
scheme of dimension $1$ over $k$ is
isomorphic to ${\mathbf G}_a$ or
a form of 
${\mathbf G}_m$.
Since there is no non-trivial morphism
${\mathbf G}_m\to {\mathbf G}_a$
of group schemes,
a non-trivial extension of ${\mathbf G}_a$
by ${\mathbf F}_p$ is isomorphic to
${\mathbf G}_a$.
Since a finite \'etale morphism
${\mathbf G}_a\to {\mathbf G}_a$
of group schemes of degree $p$
with split kernel
is $a(F-b^{p-1})$ for $a,b\in k^\times$,
the morphism {\rm (\ref{eqExt2})}
for $E={\mathbf G}_a$
is a surjection.
\qed
}

\medskip
By the isomorphism
{\rm (\ref{eqExt2})},
we identify 
${\rm Ext}(E,{\mathbf F}_p)=E^\vee$.
By the isomorphism
{\rm (\ref{eqHES})},
we identify
$H^1(E,{\mathbf F}_p)
=
H^1(k,{\mathbf F}_p)
\oplus
\bigoplus_{p\nmid m}
S^{mp^\infty}(E^\vee)$.
Then, the injection
${\rm Ext}(E,{\mathbf F}_p)
\to
H^1(E,{\mathbf F}_p)$
{\rm (\ref{eqExt})}
is induced by
$E^\vee\to S^{p^\infty}(E^\vee)$.
The injection
$E^\vee\to S^{p^\infty}(E^\vee)$
is an isomorphism if and only 
if $E$ is of dimension 1.

Assume that $k$ is separably
closed and 
let $\pi_1(E,0)_{\text{\rm pro-}p}$ denote
the maximum pro-$p$ quotient
of the \'etale fundamental group
defined by the base point at the origin $0$.

\begin{cor}\label{corpiE}
Let $k\to k'$ be a morphism of
separably closed fields of characteristic $p>0$.
Let $E$ and $E'$ be a $k$-vector space
and a $k'$-vector space of finite dimension
respectively and
let $E'\to E_{k'}=E\otimes_kk'$ be a $k'$-linear morphism.

{\rm 1.}
The following conditions are equivalent:

{\rm (1)}
The morphism of schemes $E'\to E$
is dominant.

{\rm (2)}
The morphism
$H^1(E,{\mathbf F}_p)\to
H^1(E',{\mathbf F}_p)$
of \'etale cohomology is injective.

{\rm (3)}
The morphism of pro-$p$ fundamental groups
$\pi_1(E',0)_{\text{\rm pro-}p}\to
\pi_1(E,0)_{\text{\rm pro-}p}$ is surjective.

{\rm 2.}
Assume that 
$k$ is algebraically closed
and that
$E'\to E_{k'}=E\otimes_kk'$ is
a surjection.
Then, the diagram
\begin{equation}
\begin{CD}
H^1(E,{\mathbf F}_p)
@>>>
H^1(E',{\mathbf F}_p)
\\
@AAA@AAA\\
{\rm Ext}(E,{\mathbf F}_p)
@>>>
{\rm Ext}(E',{\mathbf F}_p)
\end{CD}
\label{eqpilcoc}
\end{equation}
is cartesian
and the horizontal arrows are injections.
\end{cor}

\proof{
1.
(1)$\Leftrightarrow$(2):
The condition (1) is equivalent
to the injectivity of
the morphism of symmetric algebras
$S(E^\vee)\to S(E^{\prime \vee})$.
Since $k\to k'$ is an injection,
this is equivalent
to the injectivity of
$S^n(E^\vee)\to S^n(E^{\prime \vee})$
for every $n\geqq 1$.
By Proposition \ref{prExt}.1,
this is equivalent
to the injectivity of
$S^{mp^\infty}(E^\vee)\to S^{mp^\infty}(E^{\prime \vee})$
for every $m\geqq 1$
prime to $p$.
Since $k$ and $k'$
are separably closed,
by Proposition \ref{prExt}.2,
this is equivalent (2).

(2)$\Leftrightarrow$(3):
Since $H^1(E,{\mathbf F}_p)$
is identified with the character group
${\rm Hom}(\pi_1(E,0)_{\text{\rm pro-}p},{\mathbf F}_p)$
it suffices to apply the following Lemma \ref{lmGp}.

2.
By taking a splitting, we may assume
$E'=E_{k'}$.
Taking a basis, we may assume $E=k^n$
and identify $S(E^\vee)=k[X_1,\ldots,X_n]$.
Then, by Proposition \ref{prExt}.2,
$H^1(E,{\mathbf F}_p)$
is identified with $\bigoplus_mkX^m$
where $m=(m_1,\ldots,m_n)$ runs multiindices
such that $p\nmid m_i$ for some $i$
and ${\rm Ext}(E,{\mathbf F}_p)$
is identified with $\bigoplus_{i=1}^nkX_i$.
Hence the assertion follows.
\qed}

\begin{lm}\label{lmGp}
Let $G\to G'$ be a morphism
of pro-$p$ groups.
The following conditions are equivalent:

{\rm (1)}
$G\to G'$ is a surjection.

{\rm (2)}
%
The induced morphism
${\rm Hom}(G',{\mathbf F}_p)\to
{\rm Hom}(G,{\mathbf F}_p)$ is an injection.
%
\end{lm}

\proof{
We may assume that $G$ and $G'$ are
finite $p$-groups.
Then, since $p$-groups are nilpotent,
the condition (1) is equivalent to
the surjectivity of the morphism
$G_{\rm ab}\to G'_{\rm ab}$ 
of abelianizations
by \cite[Chapter 1, \S3 Proposition 8, Corollary 4]{NB}.
By Nakayama's lemma,
this is equivalent to the surjectivity
of $G_{\rm ab}/G^p_{\rm ab}\to 
G'_{\rm ab}/G_{\rm ab}^{\prime p}$.
This is equivalent to (2).
%
\qed}

\begin{df}\label{dflin}
Let $k$ be a field of characteristic $p>0$
and let $E$ be a $k$-vector space
of finite dimension.
Let $G$ be a finite group
and let $f\colon H\to E$ be
a $G$-torsor on $E$.

{\rm 1.}
We say that the $G$-torsor $H$ on $E$
is {\em additive} if the following conditions
are satisfied:
There exists an isomorphism 
$H \to {\mathbf G}_a^n$ for $n=\dim_kE$
such that
$f\colon H\to E$ is an \'etale morphism
of group schemes over $k$
and that the action of
$G$ on $H$ is the same
as the translation by 
$G={\rm Ker}\, f$.

{\rm 2.}
We define a morphism
\begin{equation}
[H]\colon
{\rm Hom}(G,{\mathbf F}_p)
\to 
H^1(E,{\mathbf F}_p)
\label{eqAv}
\end{equation}
by sending a character
$\chi\colon G\to {\mathbf F}_p$
to the image of the
class
$[H]\in H^1(E,G)$
by $\chi_*\colon H^1(E,G)
\to H^1(E,{\mathbf F}_p)$.
\end{df}

If $H$ is additive, 
the morphism
{\rm (\ref{eqAv})}
induces a morphism
\begin{equation}
[H]\colon
{\rm Hom}(G,{\mathbf F}_p)
\to 
E^\vee
={\rm Ext}(E,{\mathbf F}_p)
\label{eqAvl}
\end{equation}
sending a character
$\chi\colon G\to {\mathbf F}_p$
to the linear form
$f\in E^\vee$
such that there exists
a commutative diagram
\begin{equation}
\begin{CD}
0@>>> G@>>> H@>>> E@>>> 0\\
@.@V{\chi}VV @VVV@VVfV@.\\
0@>>> {\mathbf F}_p@>>> 
{\mathbf G}_a@>{F-1}>> 
{\mathbf G}_a@>>> 0
\end{CD}
\end{equation}
of exact sequences of smooth group schemes over $k$.

An additive torsor over 
${\mathbf G}_a$ is described as follows.

\begin{lm}\label{lmapol}
Let $k$ be a field of characteristic $p>0$
and let $G\subset k$
be a finite subgroup of the additive group.

{\rm 1.}
The separable polynomial
$a=\prod_{\sigma\in G}(X-\sigma)
\in k[X]$ is an additive polynomial and
the sequence
\begin{equation}
\begin{CD}
0@>>> G
@>>>
{\mathbf G}_a
@>a>>
{\mathbf G}_a
@>>>0
\end{CD}
\label{eqapol}
\end{equation}
is exact.

{\rm 2.}
Let $\chi\colon G\to {\mathbf F}_p$
be a character and let
$\tau\in G$ be an element
such that $\chi(\tau)=1$.
Let $b=\prod_{\sigma\in 
{\rm Ker}\, \chi}(X-\sigma)
\in k[X]$ be the additive polynomial.
Then the morphism
$G^\vee={\rm Hom}(G,{\mathbf F}_p)
\to k={\rm Ext}({\mathbf G}_a,
{\mathbf F}_p)$
{\rm (\ref{eqAvl})}
defined by
{\rm (\ref{eqapol})}
maps $\chi$ to
$1/b(\tau)^p$.
\end{lm}

\proof{
1.
Since $G=a^{-1}(0)$,
it suffices to show that $a$
is an additive polynomial.
We show this by induction on
$\dim_{{\mathbf F}_p}G$.
If $G=0$, then $a=X$ is an additive polynomial.
Let $\chi\colon G\to {\mathbf F}_p$
be a character and let
$\tau\in G$ be an element
such that $\chi(\tau)=1$
as in 2.
The polynomial $b\in k[X]$
defined in 2
is an additive polynomial by
the induction hypothesis.
Hence
\begin{equation}
a=\prod_{i\in {\mathbf F}_p}
b(X-i\tau)
=\prod_{i\in {\mathbf F}_p}
(b(X)-ib(\tau))
=b(X)^p-b(\tau)^{p-1}b(X)
\label{eqab}
\end{equation}
is an additive polynomial.

2.
By (\ref{eqab}),
for $b_1=b/b(\tau)
\in k[X]$,
we have a commutative diagram
$$\begin{CD}
0@>>> G
@>>>
{\mathbf G}_a
@>a>>
{\mathbf G}_a
@>>>0\\
@.@V{\chi}VV
@VV{b_1}V
@VV{\cdot \frac
1{b(\tau)^p}}V@.\\
0@>>> {\mathbf F}_p
@>>>
{\mathbf G}_a
@>{F-1}>>
{\mathbf G}_a
@>>>0
\end{CD}$$
and the assertion follows.
\qed

}

\begin{pr}\label{prlin}
Let $k$ be a field of characteristic $p>0$
and let $E$ be a $k$-vector space
of finite dimension.
Let $G$ be a finite group
and let $f\colon H\to E$ be
a $G$-torsor on $E$.
Then, 
the following conditions {\rm (1)}--{\rm (3)}
are equivalent to each other:

{\rm (1)}
$H$ is additive.

{\rm (2)}
There exists a
structure of group scheme over $k$
on $H$ such that
$f\colon H\to E$ is a morphism
of group schemes over $k$
and $H$ is connected.

{\rm (3)}
The group $G$
is an ${\mathbf F}_p$-vector space
and the morphism
$[H]\colon 
{\rm Hom}(G,{\mathbf F}_p)
\to
H^1(E,{\mathbf F}_p)$
is an injection 
to ${\rm Ext}(E,{\mathbf F}_p)=E^\vee$.
\end{pr}

We use the following to prove the implication
(2)$\Rightarrow$(1).

\begin{lm}\label{lmGH}
Let $f\colon H\to G$
be a finite \'etale morphism of smooth
connected group schemes over 
a field $k$.

{\rm 1.}
$H$ is a central extension of
$G$ by the kernel $N$.

{\rm 2.}
If $G$ is abelian,
then $H$ is abelian.
Further if $G$ is killed by an integer $n\geqq 1$,
then $H$ is also killed by $n$.

{\rm 3.}
If $G$ is a vector space over $k$,
then $H$ is also isomorphic
to a vector space over $k$.
\end{lm}

\proof{
1.
By extending $k$ if necessary,
we may assume that $N$ is split.
For $a\in N$,
the composition of
the morphism $[\, ,a]\colon H\to H$
defined by the commutator
$[\, ,a]$ with $f\colon H\to G$
is the constant morphism
and induces a morphism
$H\to N$.
Since $H$ is connected,
this implies that the morphism
$[\, ,a]\colon H\to H$ is also constant.
Hence $H$ is a central extension of
$G$ by $N$.

2.
Since $H$ is an $N$-torsor over $G$,
the morphism
$H\times N\to H\times_GH$
defined by $(x,y)\mapsto (x,xy)$
is an isomorphism.
Since $G$ is commutative,
the pair of morphisms
$H\times H\to H$ defined by
$(x,y)\mapsto xy$ and 
$(x,y)\mapsto yx$ defines 
a morphism $H\times H\to H\times_GH$.
Since $H$ is connected,
the composition with the inverse of
the isomorphism 
$H\times N\to H\times_GH$
and the projection $H\times N\to N$
is the constant morphism
$H\times H\to N$
defined by the unit $e\in N$. 
Hence $H$ is commutative.

Similarly, the morphism
$H\to H$ defined by $x\mapsto x^n$
induces the constant morphism
$H\to N$
defined by the unit $e\in N$.

3.
If $k$ is of characteristic $0$,
the morphism $f\colon
H\to G$ is an isomorphism.
Hence, we may assume
that $k$ is of characteristic $p>0$. 
Assume that $G$ is a vector space.
We show that $H$ is also isomorphic
to a vector
space by induction on 
the dimension of $N$
as an ${\mathbf F}_p$-vector space.
We may assume that
$N={\mathbf F}_p$.
By the isomorphism (\ref{eqExt2}), 
$H$ is the pull-back of 
the Artin-Schreier extension
$0\to {\mathbf F}_p
\to {\mathbf G}_a
\to {\mathbf G}_a\to 0$
by a projection $G\to {\mathbf G}_a$
and the assertion follows.
\qed
}

\proof[Proof of Proposition {\rm \ref{prlin}}]
{
(1)$\Rightarrow$(3):
Since $G$ is a subgroup
of a vector space $H$,
it is an ${\mathbf F}_p$-vector space
and the image of the morphism
$[H]\colon 
{\rm Hom}(G,{\mathbf F}_p)
\to
H^1(E,{\mathbf F}_p)$
is a subgroup of
${\rm Ext}(E,{\mathbf F}_p)$.
For any non-trivial character $\chi\colon G\to {\mathbf F}_p$,
the extension $\chi_*H$ is connected
and its class is non-trivial.

(3)$\Rightarrow$(2):
Since $H$ is an extension of $E$
by $G$,
it is a group scheme over $k$.
Let $H^\circ\subset H$
be the connected component
containing $0$.
If $\chi\colon G\to {\mathbf F}_p$
is a character trivial on
the intersection
$G^\circ =G\cap H^\circ$,
then $\chi_*H$
is a trivial extension of $E$
by ${\mathbf F}_p$.
Hence we obtain $\chi=0$,
$G=G^\circ$ and $H=H^\circ$.

(2)$\Rightarrow$(1):
Since $H$ is connected,
by Lemma \ref{lmGH},
$H$ is isomorphic to a vector space
and is a central extension of
$E$ by the kernel $N={\rm Ker}\, f$.
Further since $H$ is connected,
the mapping
$G={\rm Aut}(H/E)
\to N$
sending $g$ to $g(0)$
is an isomorphism
of groups.
Hence, $G$ is an ${\mathbf F}_p$-vector
space and the action
of $G$ on $H$
is the same as the
translation by $N=G$.
%
%
\qed

}
\medskip

Let $k\to k'$ be a morphism
of 
fields
of characteristic $p>0$
and let $E$ and $E'$
be a $k$-vector space
and a $k'$-vector space respectively.
Let $G'\to G$ be a morphism 
of finite groups
and $H$ be a $G$-torsor over $E$
and $H'$ be a $G'$-torsor over $E'$.
Let $E'\to E$ be a morphism
of schemes induced by
a $k'$-linear morphism
$E'\to E\otimes_kk'$
and 
$$\begin{CD}
H'@>>> E'\\
@VVV@VVV\\
H@>>>E
\end{CD}$$ be a 
commutative diagram
of schemes
compatible with
$G'\to G$.
Then, since
the $G$-torsor
$H\times_EE'\to E'$ 
is isomorphic to
the direct image of
the $G'$-torsor $H'\to E'$
by $G'\to G$,
the diagram
\begin{equation}
\begin{CD}
{\rm Hom}(G',{\mathbf F}_p)
@>{[H']}>>
H^1(E',{\mathbf F}_p)
\\
@AAA@AAA&\\
{\rm Hom}(G,{\mathbf F}_p)
@>{[H]}>>
H^1(E,{\mathbf F}_p)
\end{CD}
\label{eqFH1}
\end{equation}
is commutative.

%
%

\begin{cor}\label{corlinb}
Let $k\to k'$ be a
morphism of 
fields
of characteristic $p>0$
and $G'\to G$ be an injection of 
finite $p$-groups.
Let $E$ (resp.\ $E'$)
be a $k$-vector space
(resp.\ $k'$-vector space)
of finite dimension and
$H$ (resp.\ $H'$) be a connected $G$-torsor 
(resp.\ connected $G'$-torsor) over $E$ 
(resp.\ over $E'$).
Let 
\begin{equation}
\begin{CD}
H'@>>> E'\\
@VVV@VVV\\
H@>>> E
\end{CD}
\label{eqFE}
\end{equation}
be a commutative diagram
of schemes such that
the left vertical arrow is
compatible with $G'\to G$
and that the right vertical arrow
is induced by a $k'$-linear
mapping $E'\to E\otimes_kk'$.

{\rm 1.}
Assume that $k$ is
separably closed and that
the morphism $E'\to E$ of schemes is dominant.
Then, the injection $G'\to G$
 is an isomorphism.

{\rm 2.}
If the $G$-torsor $H$ is additive,
then the $G'$-torsor $H'$ is also additive.

{\rm 3.}
Assume that 
the $k'$-linear
mapping $E'\to E\otimes_kk'$
is a surjection
and that $H$ is geometrically connected.
Then $G'\to G$ is an isomorphism. 
The $G$-torsor $H$ is additive
if and only if the $G'$-torsor $H'$ is additive.
\end{cor}

\proof{
1.
Let $\bar k'$ be a separable closure of $k'$.
By replacing $k', H'$ and $G'$ by
$\bar k'$, a connected component of
the base change $H'_{k'}$
and its stabilizer,
we may assume that $k'$ is also
separably closed.
Since $H$ and $H'$ are connected,
the commutative diagram
(\ref{eqFE}) induces a commutative diagram
$$\begin{CD}
\pi_1(E',0)_{\text{\rm pro-}p}@>>> G'\\
@VVV@VVV\\
\pi_1(E,0)_{\text{\rm pro-}p}@>>> G
\end{CD}$$
where the horizontal arrows are surjections.
The left vertical arrow is
also a surjection by the assumption
that $k$ and $k'$ are
separably closed
and that $E'\to E$ is dominant by
Corollary \ref{corpiE}.1.
Hence $G'\to G$ is also a surjection.

2.
Since $G$ is an ${\mathbf F}_p$-vector
space and since $G'\to G$
is an injection,
$G'$ is also an ${\mathbf F}_p$-vector
space.
By the commutative diagram
(\ref{eqFH1}),
if the image of
$[H]$ is a subgroup of $E^\vee$,
the image of
$[H']$ is a subgroup of $E^{\prime \vee}$.
Since $H'$ is connected,
the morphism $[H']$ is an injection.
Hence $H'$ is additive by Proposition \ref{prlin}
(3)$\Leftrightarrow$(1).

3.
We show that
$G'\to G$ is an isomorphism.
Since $H$ is geometrically connected,
by replacing $k$ and $k'$
by algebraic closures,
we may assume that
$k$ and $k'$
are algebraically closed,
similarly as in the beginning 
of the proof of 1.
If $E'\to E\otimes_kk'$ is surjective,
the composition 
$E'\to E\otimes_kk'\to E$
is also surjective and
$G'\to G$ is an isomorphism by 1.

Since the if part is proved in 2.,
it suffices to show the only if part.
Since $G'$ is an ${\mathbf F}_p$-vector
space and since $G'\to G$
is an isomorphism,
$G$ is also an ${\mathbf F}_p$-vector space.
In the commutative diagram (\ref{eqFH1}),
the assumption implies that
the vertical arrows are injections
and hence $[H]$ is an injection
if $[H']$ is an injection.
By Corollary \ref{corpiE}.2,
the image of
$[H]$ is a subgroup of $E^\vee$
if the image of
$[H']$ is a subgroup of $E^{\prime \vee}$.
Hence $H$ is additive by Proposition \ref{prlin}
(3)$\Leftrightarrow$(1).
\qed}

\subsection{\'Etale isogenies of smooth group schemes}\label{ssGH}

\begin{lm}\label{lmXY}
Let $S$ be a scheme
and $X$ and $Y$ be
finite \'etale schemes over $S$.
Then, the functor
${\mathcal I}som_S(X,Y)$
sending a scheme $T$ over $S$
to the set
${\rm Isom}_T(X_T,Y_T)$
of isomorphisms
$X_T\to Y_T$ over $T$
is representable by
a finite \'etale scheme 
$I\!som_S(X,Y)$
over $S$.
\end{lm}

\proof{
The functor ${\mathcal I}som_S(X,Y)$
is representable by
a closed subscheme
$I\!som_S(X,Y)$
of the vector bundle
associated to the locally free
${\mathcal O}_S$-module
${\mathcal H}om_{{\mathcal O}_S}({\mathcal O}_Y,{\mathcal O}_X)$.

We show that
the scheme $I\!som_S(X,Y)$
is finite \'etale over $S$.
Since the assertion is \'etale local on
$S$
and since the sheaf
${\mathcal I}som_S(X,Y)$ is \'etale locally
constant,
the assertion follows.
\qed}

\begin{lm}\label{lmetale}
Let $X\to Y$ be a finite \'etale morphism
of schemes and
let $S$ be a normal integral scheme
over $Y$.
Let $K$ be the residue field
of the generic point of $S$
and let $K'$ be an extension of $K$.
Then, the mapping
\begin{equation}
\begin{CD}
{\rm Mor}_Y(S,X)
@>>>
{\rm Mor}_Y({\rm Spec}\, K',X)
\end{CD}
\label{eqet}
\end{equation}
sending a morphism
$f\colon S\to X$ over $Y$
to the composition
${\rm Spec}\, K'\to S\to X$
is an injection.
The image of {\rm (\ref{eqet})}
is the subset
defined by the condition
that the image of the induced morphism
${\rm Spec}\, (K'\otimes_KK')
\to X\times_YX$
is a subset of the diagonal
$X\subset X\times_YX$.
\end{lm}

\proof{
By replacing $X\to Y$
by the base change
$X\times_YS\to S$,
we may assume $Y=S$.
Since $S$ is normal,
the \'etale scheme $X$ over $S$
is also normal.
Further replacing $S$ by
${\rm Spec}\, K$,
we may assume $S={\rm Spec}\, K$.
Since $K\to K'$ is injective,
the mapping (\ref{eqet}) is an injection.

Assume that the image of 
${\rm Spec}\, (K'\otimes_KK')
\to X\times_YX$
is a subset of the diagonal
$X\subset X\times_YX$.
By replacing $X$ by the
image of ${\rm Spec}\, K'\to X$
and further $K'$ by
$\Gamma(X,{\mathcal O}_X)$,
we may assume 
${\rm Spec}\, K'= X$.
Then, the assumption means
that the diagonal
$X\to X\times_YX$ is
an isomorphism.
Hence, 
the morphism $X\to Y$ is 
of degree 1 and
is an isomorphism.
\qed}

\begin{pr}\label{prET}
Let $G$ be a smooth connected group scheme over a field $k$
with the group structure
$\mu\colon G\times G\to G$
and the unit $e\in G(k)$.

Let $f\colon H\to G$
be an \'etale morphism of schemes over $k$. 
Let 
$\nu\colon H\times H\to H$
be a morphism of schemes over $k$
such that the diagram
\begin{equation}
\begin{CD}
H\times H@>{\nu}>> H\\
@V{f\times f}VV@VVfV\\
G\times G@>{\mu}>> G
\end{CD}
\label{eqfGH}
\end{equation}
is commutative
and let $e'$ be a $k$-valued point
of $N=H\times_Ge$.
We consider the following
conditions:

{\rm (1)}
The morphism
$\nu\colon H\times H\to H$
is a group structure and
$e'\in H(k)$ is the unit.

{\rm (2)}
The morphism
$\widetilde \nu\colon 
H\times H\to H\times H$ 
defined by $\widetilde \nu(x,y)=(\nu(x,y), y)$
is an isomorphism
and $\nu(e',e')=e'$.

Then, we have
{\rm (1)}$\Rightarrow${\rm (2)}.
If $H$ is geometrically connected
over $k$, 
we have
{\rm (2)}$\Rightarrow${\rm (1)}.
\end{pr}

\proof{
The implication 
{\rm (1)}$\Rightarrow${\rm (2)}
is obvious.
Assume that $H$ is geometrically connected.
The inverse image of
the diagonal $G\subset G\times G$ 
by the isomorphism
$\widetilde \mu\colon 
G\times G\to G\times G$ 
defined by $\widetilde \mu(x,y)=(\mu(x,y), y)$
is $e\times G$.
The condition (2) implies
that the diagram
\begin{equation}
\begin{CD}
H\times H@>{\widetilde\nu}>> H\times H\\
@V{f\times f}VV@VV{f\times f}V\\
G\times G@>{\widetilde \mu}>> G\times G
\end{CD}
\label{eqfGHb}
\end{equation}
is cartesian.
Hence the isomorphism
$\widetilde \nu$ induces
an isomorphism
$N\times H\to H\times_GH$.
This induces an isomorphism
from $e'\times H$
to the diagonal $H$
since $\nu(e',e')=e'$
and $H$ is connected.

Since (\ref{eqfGH})
is commutative,
the pair of morphisms
$H\times H\times H\to H$
defined by
$\nu(\nu(x,y),z)$
and by $\nu(x,\nu(y,z))$
defines a morphism
$H\times H\times H\to H
\times_GH$.
Since $H$ is geometrically connected,
its composition with
the projection
$H\times_G H
\to N\times H\to N$
is constant.
Since $\nu(e',e')=e'$,
its image is $e'$ and 
the morphism $\nu$ satisfies
the associativity law.
Similarly,
$e'$ satisfies the property of the unit element.
Further the composition of
the morphisms
$$\begin{CD}
H@>{(e',1_H)}>> H\times H
@>{\widetilde \nu^{-1}}>>
H\times H
@>{{\rm pr}_1}>>H
\end{CD}$$
gives the inverse $H\to H$.
\qed}

\begin{pr}\label{prisog}
Let $G$ be a smooth geometrically 
connected group scheme over a field $k$
with the group structure
$\mu\colon G\times G\to G$
and the unit $e\in G(k)$.
Let $f\colon H\to G$
be an $A$-torsor
for a finite group $A$
and assume that $H$ is geometrically connected
over $k$.
Assume that $N=f^{-1}(e)$
has a $k$-rational point.
Let $\xi\colon 
{\rm Spec}\, K\to G$ 
and $\eta\colon 
{\rm Spec}\, L\to H$ 
be the generic points.
Let $L'$ be an extension of $L$
and $\eta'\colon {\rm Spec}\, L'\to H$ 
be the composition with $\eta$.

Let $\mu_K\colon G_K=G\times_
{{\rm Spec}\, k}
{\rm Spec}\, K\to G_K$
be the translation $\cdot\xi$
by the generic point
$\xi$ defined as
the morphism induced by
the composition
$\mu\circ(1_G,\xi)\colon
G\times_
{{\rm Spec}\, k}
{\rm Spec}\, K
\to G\times G\to G$.
Similarly for a morphism $\nu\colon 
H\times H\to H$,
let $\nu_{L'}\colon H_{L'}=H\times_
{{\rm Spec}\, k}
{\rm Spec}\, L'\to H_{L'}$
be the morphism induced by
the composition
$\nu\circ(1_H,\eta')\colon
H\times_
{{\rm Spec}\, k}
{\rm Spec}\, L'
\to H\times H\to H$.

Let $M$ be the set of
group structures 
$\nu\colon H\times H\to H$
of schemes over $k$
such that $f\colon H\to G$
is a morphism of group schemes.
Let $M'$ be the set of
isomorphisms 
$\nu'\colon H_{L'}\to H_{L'}$
of schemes over $L'$
such that the diagram
$$
\begin{CD}
H_{L'}@>{\nu'}>>H_{L'}\\
@VVV@VVV\\
G_K@>{\mu_K}>>G_K
\end{CD}$$
is commutative.
Then, the mapping
$M\to M'$
sending $\nu$ to $\nu_{L'}$
is a bijection.
\end{pr}

To prove the proposition,
it suffices to show the following.

\begin{lm}\label{lmisog}
Let the notation be as in Proposition
{\rm \ref{prisog}}.
Let $M_1$ be the set of
morphisms 
$\nu\colon H\times H\to H$
of schemes over $k$
such that the diagram {\rm (\ref{eqfGHb})}
is commutative and that
$\widetilde \nu\colon
H\times H\to H\times H$
defined by $\widetilde \nu(x,y)=(\nu(x,y), y)$
is an isomorphism.

{\rm 1.}
The mapping
$M_1\to M'$
sending $\nu$ to $\nu_{L'}$
is canonically identified
with the mapping
\begin{equation}
{\rm Isom}_{G\times H}
(H\times H,
(G\times H)
\times_GH)
\to
{\rm Isom}_{G_{L'}}
(H_{L'},
G_{L'}\times_GH)
\label{eqGHL}
\end{equation}
where the fiber products
are taken with respect to
$\mu\circ (1_G\times f)
\colon G\times H
\to G$
and its composition with
$G_{L'}\to G\times H$.

{\rm 2.}
The mapping
{\rm (\ref{eqGHL})}
a bijection.

{\rm 3.}
We have $M=M_1$.
\end{lm}

\proof{
1.
For $\nu\colon H\times H\to H$,
in the commutative diagram
$$\begin{CD}
H\times G@<<< H\times H
@>{\widetilde \nu}>> H\times H@>{{\rm pr}_1}>>H\\
@VVV@VVV@VVV@VVV\\
G\times G@<<< G\times H
@>{(\mu\circ(1_G\times f),
{\rm pr}_2)}
>> G\times H@>{{\rm pr}_1}>>G,
\end{CD}$$
the left and the right squares
are cartesian.
Hence $M'$
is identified with
the source of (\ref{eqGHL}).
Since the diagram
\begin{equation}
\begin{CD}
H_{L'}@>{\nu_{L'}}>> H_{L'}\\
@V{1_H\times \eta'}VV
@VV{1_H\times \eta'}V\\
H\times H@>{\widetilde \nu}>> H\times H
\end{CD}
\label{eqnn'}
\end{equation}
is cartesian,
similarly
the set $M_1$
is identified with the target of (\ref{eqGHL}).

2.
Define a finite \'etale scheme $\mu^*H$ over $G\times G$
by the cartesian diagram
$$\begin{CD}
\mu^*H@>>>H\\
@VVV@VVV\\
G\times G@>{\mu}>>G
\end{CD}$$
and consider the scheme $X=
I\!som_{G\times G}(H\times G,\mu^*H)$
of isomorphisms.
Then, the mapping (\ref{eqGHL})
is identified with the canonical mapping
\begin{equation}
{\rm Mor}_{G\times G}(G\times H,X)
\to
{\rm Mor}_{G\times G}(G_{L'},X).
\label{lmGGK}
\end{equation}

We show the bijectivity of
(\ref{lmGGK}) 
by applying Lemma \ref{lmetale}.
Let $L''$ be the
function field of $G_{L'}$.
Then, 
by Lemma \ref{lmetale} applied to
the morphisms
$$\begin{CD}
{\rm Spec}\, L''
@.@. X\\
@VVV@.@VVV\\
G_{L'}@>>>
G\times H@>{1_G\times f}>> G\times G,
\end{CD}$$
the source and the target
of (\ref{lmGGK}) 
are identified with subsets
of ${\rm Mor}_{G\times G}
({\rm Spec}\, L'',X)$.
Hence
the mapping 
(\ref{lmGGK}) is an injection.
Further, its image is
the subset defined by the following condition:

(1)
The image of
the induced morphism
$G_{L'}\times_{G\times H}G_{L'}
=G\times {\rm Spec}\, (L'\otimes_LL')
\to X\times_{G\times G}X$
is contained in the diagonal
$X\subset X\times_{G\times G}X$.

\noindent
Since $G$ is geometrically connected
and since
the diagonal
$X\subset X\times_{G\times G}X$
is an open and closed subscheme,
the condition (1) is equivalent to 
the following condition:

(1$'$) The image of
$e\times {\rm Spec}\, (L'\otimes_LL')
\subset G\times {\rm Spec}\, (L'\otimes_LL')
\to X\times_{G\times G}X$
is contained in the diagonal
$X\subset X\times_{G\times G}X$.

It suffices to show that the
condition (1$'$) is satisfied
for any morphism ${\rm Spec}\, L'
\to X$ over $G\times G$.
By Lemma \ref{lmetale}
applied to the diagram
$$\begin{CD}
{\rm Spec}\, L'
@. X\\
@VVV@VVV\\
H@>{(e,f)}>> G\times G,
\end{CD}$$
the condition (1$'$) is equivalent
to the condition: 

(2) The mapping
\begin{equation}
{\rm Mor}_{G\times G}(H,X)
\to
{\rm Mor}_{G\times G}({\rm Spec}\, L',X)
\label{lmGHK}
\end{equation}
is a bijection.

\noindent
Since the diagrams
$$\begin{CD}
H\times G@<<< N\times H\\
@VVV@VVV\\
G\times G@<{(e,f)}<<H,
\end{CD}
\qquad
\begin{CD}
H\times_GH@.
\longrightarrow@.H\\
@VVV@.@VVV\\
H@>{(e,f)}>>G\times G@>{\mu}>>G
\end{CD}
$$
are cartesian,
the condition (2)
means the condition:

(2$'$)
The mapping
\begin{equation}
{\rm Isom}_H(N\times H,H\times_GH)
\to
{\rm Isom}_{L'}(N_{L'},H\times_G{\rm Spec}\, L')
\label{lmHK}
\end{equation}
is a bijection.

Since $H$ is an $A$-torsor over $G$,
the fiber product
$H\times_GH$ is canonically identified with
$A\times H$.
Since $N(k)\subset H(k)\neq \varnothing$
and $H$ is connected,
the two sets in (\ref{lmHK}) are canonically
identified with the finite
set ${\rm Isom}(N,A)$.
Hence
the mapping
{\rm (\ref{eqGHL})}
a bijection.

3.
By Proposition \ref{prET},
it suffices to show that,
for $\nu\colon H\times H\to H$
in $M_1$,
there exists a unique element
$e\in N$ such that $\nu(e,e)=e$.
Since such an $e$
is the unit element with
respect to the group structure $\nu$,
it is unique.

We show the existence.
Let $a\in N$ and set
$\nu(a,a)=b\in N$.
Then, there exists a unique element
$\sigma\in A$ such that
$b=\sigma(a)$.
Then, by Proposition \ref{prET},
$\nu'=\sigma^{-1}\circ \nu$ is
a group structure and $a$
is the unit.
Since $\sigma$ is a unique automorphism
of $H$ over $G$ sending $a$
to $b$,
we have $\sigma(x)=\nu'(b,x)$.
Define $e\in N$ by $\nu'(b,e)=a$.
Then, since
$\nu'(\nu'(b,e),\nu'(b,e))=
\nu'(a,a)=a=\nu'(b,e)$,
we obtain
$\nu(e,e)=\nu'(b,\nu'(e,e))=e$.
}

\section{Construction of functors}\label{sFun}

We recall the geometric construction
on which the definition of ramification groups is based in Section \ref{ssFun}.
A functorial construction of
finite \'etale coverings
of twisted normal spaces is also recalled.
This plays a fundamental role in
the proof of Theorem \ref{thmgr}.
As a preliminary for the construction
in Section \ref{ssFun}, 
we recall some basic properties
of dilatations in Section \ref{ssDil}.

We prove in Proposition \ref{prisom}
in Section \ref{sskey}
that a finite covering of
a twisted normal space
constructed in Section \ref{ssFun}
satisfies a condition in the criterion
(Proposition \ref{prisog})
to be a morphism of group schemes.
A crucial ingredient in
the proof of Proposition \ref{prisom}
is a descent property
proved in Proposition \ref{prnormal}
of \'etaleness of finite morphism
in Section \ref{ssNorm}.

\subsection{Dilatations}\label{ssDil}

\begin{df}\label{dfdil}
Let $X$ be a scheme,
$D\subset X$ be an effective Cartier
divisor and $Y\subset X$
be a closed subscheme.
Let $j\colon X\sm D=U\to X$
be the open immersion
of the complement
and let ${\mathcal I}_D,{\mathcal I}_Y\subset {\mathcal O}_X$
denote the ideal sheaves
defining $D,Y\subset X$.
We define the dilatation
$X^{[D\cdot Y]}$
to be the scheme affine over $X$
defined by the quasi-coherent
sheaf of ${\mathcal O}_X$-subalgebra
of $j_*{\mathcal O}_U$
generated by ${\mathcal I}_D^{-1}
{\mathcal I}_Y$.
\end{df}

If $X={\rm Spec}\, A$
and if $Y$ and $D$
are defined by an ideal
$I\subset A$
and a non-zero divisor $a\in A$,
we have $X^{[D\cdot Y]}
={\rm Spec}\, A\bigl[\dfrac Ia\bigr]$
for the subring $A\bigl[\dfrac Ia\bigr]
\subset
A\bigl[\dfrac 1a\bigr]$.

The dilatations have the following functoriality.
Let 
\begin{equation}
\begin{CD}
D\times_XX'&\subset D'\subset &X'@<<<Y'\\
@VVV@VfVV@VVV\\
D@>>> X@<<<Y
\end{CD}
\label{eqdilf}
\end{equation}
be a commutative diagram of schemes
such that the horizontal arrows
are closed immersions and that
$D\subset X$ and $D'\subset X'$
are effective Cartier divisors.
Then the pull-back
$f^*D=D\times_XX'\subset X'$
is an effective Cartier divisor 
and the morphism $f\colon X'\to X$
is uniquely lifted to a morphism
$X^{\prime [D'\cdot Y']}
\to X^{[D\cdot Y]}$.
In particular, if
$D\times_XY$ is an effective Cartier divisor
of $Y$, the closed immersion
$Y\to X$ is uniquely lifted to a 
closed immersion $Y\to X^{[D\cdot Y]}$.

\begin{lm}\label{lmdilfl}
Assume that in
the diagram {\rm (\ref{eqdilf})},
the divisor $D'$ equals the pull-back $f^*D$.

{\rm 1.}
Assume that the diagram 
{\rm (\ref{eqdilf})} is
cartesian and
that $f$ is flat.
Then, the canonical morphism
$X^{\prime [D'\cdot Y']}
\to X^{[D\cdot Y]}\times_XX'$
is an isomorphism.

{\rm 2.}
Assume that
the morphism $f$ is smooth
and that $Y'\to Y$ is an isomorphism.
Then the canonical morphism
$X^{\prime [D'\cdot Y']}
\to X^{[D\cdot Y]}$
is also smooth
and the induced morphism
$X^{\prime [D'\cdot Y']}
\times_{X'}D'\to
X^{[D\cdot Y]}\times_XD$
is the projection of a vector bundle.
\end{lm}

\proof{
1.
Let $j\colon X\sm D=U\to X$
and $j'\colon X'\sm D'=U'\to X'$
denote the open immersions
of the complements.
Then, since $f$ is
flat and $D'=f^*D$,
the canonical morphism
$j_*{\mathcal O}_U
\otimes_{{\mathcal O}_X}{\mathcal O}_{X'}
\to
j'_*{\mathcal O}_{U'}$
is an isomorphism
and induces an isomorphism
on the subalgebras.

2.
Since the assertion is local on a neighborhood
of $Y'$,
we may assume that
$X={\rm Spec}\, A$ and
$X'={\rm Spec}\, A'$ are affine
and that there exists $t_1,\ldots,t_n\in A'$
such that $X'\to X'_0={\mathbf A}^n_X
={\rm Spec}\, A[T_1,\ldots,T_n]$
defined by $t_1,\ldots,t_n$
is \'etale and the composition
$Y=Y'\to X'\to X'_0={\mathbf A}^n_X$
is the 0-section.

Assume that $Y\subset X$
and $D\subset X$
are defined by an ideal
$I\subset A$ and a non-zero divisor
$a\in A$.
For $X'_0={\mathbf A}^n_X$
and the $0$-section
$Y\to {\mathbf A}^n_X$,
we have an isomorphism
$X_0^{\prime [D'_0\cdot Y']}
\to {\mathbf A}^n_{X^{[D\cdot Y]}}
={\rm Spec}\, A\Bigl[\dfrac Ia\Bigr]
\Bigl[\dfrac{T_1}a,\ldots,\dfrac{T_n}a\Bigr]$.
Since
$X^{\prime [D'\cdot Y']}
\to
X_0^{\prime [D'_0\cdot Y']}
\times
_{X'_0}X'$ is an isomorphism by 1, 
the morphism
$X^{\prime [D'\cdot Y']}
\to X^{[D\cdot Y]}$ is smooth.

The inverse images
 $X^{\prime [D'\cdot Y']}\times_{X'}D'$
and $X_0^{\prime [D'_0\cdot Y']}\times_{X'_0}D'_0$
are schemes over $Y'$
and 
the morphism
$X^{\prime [D'\cdot Y']}
\to
X_0^{\prime [D'_0\cdot Y']}$
induces an isomorphism
on the inverse images
of $Y'$.
Hence the morphism
$X^{\prime [D'\cdot Y']}\times_{X'}D'
\to 
X_0^{\prime [D'_0\cdot Y']}\times_{X'_0}D'_0$ 
is an isomorphism
and the scheme $X^{\prime [D'\cdot Y']}\times_{X'}D'$ is a vector bundle
over $X^{[D\cdot Y]}\times_XD$.
\qed

}


\begin{pr}\label{prdiltr}
Let $X$ be a scheme
and
let $D\subset X$ be a Cartier divisor.
Let $Y\subset X$ be
a closed subscheme
and let $Z=Y\times_XD$
be the intersection.

{\rm 1.}
The canonical morphism
$X^{[D\cdot Z]}\to 
X^{[D\cdot Y]}$
of dilatations is an isomorphism.

{\rm 2.}
Assume that the immersion
$Y\to X$
is a regular immersion
and that $Z\subset Y$ is a Cartier divisor.
Let $T_YX$ and $T_YX(-D)$
be the normal bundle
and its twist
by the Cartier divisor $D$.
Then, we have a cartesian diagram
\begin{equation}
\begin{CD}
T_YX(-D)\times_YZ@>>>
X^{[D\cdot Y]}\\
@VVV@VVV\\
D@>>>X.
\end{CD}
\label{eqdilC}
\end{equation}
\end{pr}

\proof{
1.
Since the assertion is local on $X$,
we may assume that
$X={\rm Spec}\, A$
and $Y$ and $D$
are defined by $I\subset A$ and $a\in A$.
Then, we have
$X^{[D\cdot Z]}=
{\rm Spec}\, A\Bigl[\dfrac {I+aA}a\Bigr]
=
{\rm Spec}\, A\Bigl[\dfrac {I}a\Bigr]
=
X^{[D\cdot Y]}$.

2.
The twisted normal bundle
$T_YX(-D)$
is the affine scheme over $Y$
defined by the quasi-coherent
${\mathcal O}_Y$-algebra
$\bigoplus_{q=0}^\infty
{\mathcal I}_D^{-q}
\otimes_{{\mathcal O}_X}
{\mathcal I}_Y^q/
{\mathcal I}_Y^{q+1}$.
We show that
the canonical morphism
$\sum_{q=0}^\infty
{\mathcal I}_D^{-q}{\mathcal I}_Y^q
\otimes_{{\mathcal O}_X}
{\mathcal O}_D
\to
\bigoplus_{q=0}^\infty
({\mathcal I}_D^{-q}
\otimes_{{\mathcal O}_X}
{\mathcal I}_Y^q/
{\mathcal I}_Y^{q+1})$ 
induces an isomorphism
$\sum_{q=0}^\infty
{\mathcal I}_D^{-q}{\mathcal I}_Y^q
\otimes_{{\mathcal O}_X}
{\mathcal O}_D
\to
\bigoplus_{q=0}^\infty
({\mathcal I}_D^{-q}
\otimes_{{\mathcal O}_X}
{\mathcal I}_Y^q/
{\mathcal I}_Y^{q+1})
\otimes_{{\mathcal O}_Y}
{\mathcal O}_Z$.
Since the question is local on $X$,
we may assume that
$X={\rm Spec}\, A$,
the closed subscheme
$Y$ is defined by
the ideal $I\subset A$ generated by 
a regular sequence
$t_1,\ldots,t_n\in A$
and the divisor $D$ 
is defined by $a\in A$
such that
$t_1,\ldots,t_n,a$
remain to be a regular sequence.
Then, 
the morphism
$A[T_1,\ldots,T_n]
/(aT_i-t_i$,
$i=1,\ldots,n)
\to A\Bigl[\dfrac 1a\Bigl]$
induces an injection
on the graded quotients
and hence itself is an injection.
Thus we obtain an isomorphism
$A[T_1,\ldots,T_n]
/(aT_i-t_i,\ i=1,\ldots,n)
\to A\Bigl[\dfrac Ia\Bigl]$
and this induces an isomorphism
$A/(aA+I)[T_1,\ldots,T_n]
\to A\Bigl[\dfrac Ia\Bigl]
/aA\Bigl[\dfrac Ia\Bigl]$.
\qed

}

\begin{cor}\label{cordiltr}
Let $X$ be a scheme
and 
let $D\subset X$ be a Cartier divisor.
Let $Y,Y'\subset X$ be
closed subschemes
such that the immersions
are regular immersions.
Assume that
the intersections
$Y\times_X D$ and $Y'\times_X D$
are the same
and are Cartier divisors of
$Y$ and of $Y'$ respectively.
We set
$Z=Y\times_X D=Y'\times_XD$
and identify
$X^{[D\cdot Y]}=X^{[D\cdot Z]}=
X^{[D'\cdot Y]}$
and
$T_YX(-D)\times_YZ
=T_ZD(-D)=
T_{Y'}X(-D)\times_{Y'}Z$.

Let $u\colon T_YX(-D)\times_YZ
\to T_{Y'}X(-D)\times_{Y'}Z$
be the composition 
of the isomorphisms
\begin{equation}
\begin{CD}
T_YX(-D)\times_YZ
@>{\simeq}>{\rm (\ref{eqdilC})}>
X^{[D\cdot Y]}\times_XD
&=X^{[D\cdot Z]}\times_XD@.\\
@.&=
X^{[D\cdot Y']}\times_XD
@<{\simeq}<{\rm (\ref{eqdilC})}<
T_{Y'}X(-D)\times_{Y'}Z
\end{CD}
\label{eqt}
\end{equation}
and
$s\colon
Z\to 
T_YX(-D)\times_YZ$
be the inverse image 
of the $0$-section of 
$T_{Y'}X(-D)\times_{Y'}Z$
by $u$.
Then $u$ 
is the composition
of the translation $+s
\colon
T_YX(-D)\times_YZ
\to 
T_YX(-D)\times_YZ$
with 
the isomorphisms
\begin{equation}
\begin{CD}
T_YX(-D)\times_YZ
@<{\simeq}<<
T_ZD(-D)@>{\simeq}>>
T_{Y'}X(-D)\times_{Y'}Z.
\end{CD}
\label{eqs}
\end{equation}
%
%
%
%
%
%
\end{cor}

\proof{
%
%
%
Since the assertion is local on $X$,
we may assume that
$X={\rm Spec}\, A$ is affine,
$Y={\rm Spec}\, A/I$,
$Y'={\rm Spec}\, A/I'$ and
$D={\rm Spec}\, A/aA$.
Then, 
the assumption $Y\times_X D=Y'\times_X D$
means the equality
$I+aA=I'+aA$.

Since $(I+aA)/(I+aA)^2=
(N_{Y/X}\otimes_{{\mathcal O}_Y}
{\mathcal O}_Z)
\oplus
(N_{D/X}\otimes_{{\mathcal O}_D}
{\mathcal O}_Z)$
and the same equality
with $I$ and $Y$ replaced
by $I'$ and $Y'$,
the equality
$I+aA=I'+aA$
induces an isomorphism
$(N_{Y/X}\otimes_{{\mathcal O}_Y}
{\mathcal O}_Z)
\oplus
(N_{D/X}\otimes_{{\mathcal O}_D}
{\mathcal O}_Z)
\to
(N_{Y'/X}\otimes_{{\mathcal O}_{Y'}}
{\mathcal O}_Z)
\oplus
(N_{D/X}\otimes_{{\mathcal O}_D}
{\mathcal O}_Z)$
of locally free ${\mathcal O}_Z$-modules,
preserving
the subsheaf
$N_{D/X}\otimes_{{\mathcal O}_D}
{\mathcal O}_Z$
but not necessarily the direct sum decomposition.
Hence the morphism
$u$ is the composition
of (\ref{eqs})
with the translation $+s'$ by a section
$s'\colon
Z\to 
T_{Y'}X(-D)\times_{Y'}Z$.
By comparing the composition
with the $0$-section
$Z\to 
T_YX(-D)\times_YZ$,
we obtain $s'=s$ and
the assertion follows.
%
\qed

}

\medskip

The naive functoriality
above is generalized as follows.

\begin{lm}\label{lmdilf}
Let 
\begin{equation}
\begin{CD}
X'@<<<Y'\\
@VfVV@VVV\\
X@<<<Y
\end{CD}
\label{eqIq}
\end{equation}
be a commutative diagram of schemes
such that the horizontal arrows
are closed immersions.
Let
$D\subset X$ and $D',D'_1\subset X'$
be effective Cartier divisors
and $q\geqq 1$ be an integer
such that the defining ideals
satisfy
${\mathcal I}_Y{\mathcal O}_{X'}
\subset
{\mathcal I}_{D'_1}{\mathcal I}_{Y'}
+
{\mathcal I}_{Y'}^q$
and 
${\mathcal I}_{D'}\cdot{\mathcal I}_{D'_1}
+{\mathcal I}_{D'}^q\subset
{\mathcal I}_D{\mathcal O}_{X'}$.

Then,
the morphism
$X'\to X$
is uniquely lifted to
$X^{\prime[D'\cdot Y']}
\to X^{[D\cdot Y]}$.
\end{lm}

The naive functoriality corresponds
to the case $q=1$.

\proof{
Since the assertion is local,
we may assume that
$X={\rm Spec}\, A$
and
$X'={\rm Spec}\, B$
are affine, that
$Y\subset X$ and
$Y'\subset X'$
are defined by ideals
$I\subset A$
and $J\subset B$
and that $D\subset X$
and $D',D'_1
\subset X'$ are defined
by non-zero divisors
$a\in A$ and
$b,c\in B$.
By the assumption,
we have
$IB\subset cJ+J^q$
and $a$ divides $bc$ and $b^q$.
This implies that the image
of $\dfrac Ia\subset A\Bigl[\dfrac 1a\Bigr]$
in $B\Bigl[\dfrac 1b\Bigr]$
is contained in
$\dfrac{cJ}{a}+\dfrac{J^q}{a}
\subset
\dfrac{J}{b}+\dfrac{J^q}{b^q}
\subset
B\Bigl[\dfrac Jb\Bigr]$
and the assertion follows.
\qed

}

\subsection{Normalizations}\label{ssNorm}

Let $X$ be a reduced scheme,
$j\colon U\to X$ be an immersion
of normal dense open subscheme
and $V\to U$ be a finite
\'etale morphism of normal schemes.
Then, the normalization $Y$ of
$X$ in $V$ is the scheme
affine over $X$
corresponding to the
quasi-coherent ${\mathcal O}_X$-algebra
${\mathcal O}_Y$ defined as the integral closure
of ${\mathcal O}_X$
in $j_*{\mathcal O}_V$.
The formation of
normalizations is functorial as follows.

\begin{lm}\label{lmnormal}
Let
$$\begin{CD}
X'@>{\supset}>>U'@<<<V'@>{\subset}>>Y'\\
@VVV@VVV@VVV@.\\
X@>{\supset}>>U@<<<V@<{\subset}<<Y
\end{CD}
$$
be a cartesian diagram of schemes.
Assume that
$X$ and $X'$ are reduced,
$U$ and $U'$ are their normal dense open
subschemes,
the middle horizontal arrows
are finite \'etale,
and $Y$ and $Y'$
are the normalizations of
$X$ and $X'$ in $V$ and $V'$
respectively.
Then, there exists a unique
morphism $Y'\to Y$
that makes a square in
the right of the diagram
commutative.
\end{lm}

\proof{
Let $Z\subset Y\times_XX'$
be the reduced closed subscheme
such that the underlying set
is the closure of
$V\times_UU'\subset Y\times_XX'$.
Then $Y'$ equals the normalization
of $Z$ in the finite \'etale scheme
$V'$ over $V\times_UU'$
and we obtain a morphism
$Y'\to Z\to Y\times_XX'\to Y$.
\qed

}

\begin{pr}\label{prnormal}
{\rm (\cite[Lemma 2.7]{cln},
\cite[Lemma 1.20]{nonlog})}
Let ${\mathcal O}_K$ be a discrete
valuation ring 
and $\eta={\rm Spec}\, K
\in S={\rm Spec}\, {\mathcal O}_K$ 
be the generic point.
Let $X$ be a smooth scheme over $S$,
let $U\subset X_K$
be a dense open subscheme
and let $V\to U$ be a finite \'etale
morphism.
Let $Y$ be the normalization of $X$
in $V$ and
let $t\colon S\to Y$ be a section
such that
$t(\eta)\in V$.

Let ${\mathcal O}_K\to {\mathcal O}_{K'}$
be a quasi-finite morphism
of discrete valuation rings
and let $X'\supset U'\gets V'$
be the base change by
$S'={\rm Spec}\, {\mathcal O}_{K'}\to S$.
Let $Y'$ be the normalization of $X'$
in $V'$ and 
let $t'\colon S'\to Y'$
be the section
induced by $t$.

Assume that
$Y'\to X'$ is \'etale on
a neighborhood of
the image $t'(S')\subset Y'$.
Then, 
$Y\to X$ is \'etale on
a neighborhood of
the image $t(S)\subset Y$.
\end{pr}

Proposition \ref{prnormal}
is a special case of
\cite[Lemma 2.7]{cln},
\cite[Lemma 1.20]{nonlog}.
For the sake of convenience,
we give a proof.

\proof{
By replacing $S$ by its strict localization,
we may assume that
$S$ and $S'$ are strictly local
and that $S'\to S$
is finite.
Let $x,y,x',y'$ be the geometric
points of $X,Y,X',Y'$
defined as the images
of the closed points
by $S\to Y\to X$
and $S'\to Y'\to X'$
and we consider the commutative
diagram
\begin{equation}
\begin{CD}
Y_{(y)}@<<<Y'_{(y')}\\
@VVV@VVV\\
X_{(x)}@<<<X'_{(x')}&\, =
X_{(x)}\times_SS'
\end{CD}
\label{eqSnorm}
\end{equation}
of strict localizations.
It suffices to prove that
the finite morphism
$Y_{(y)}\to X_{(x)}$
is an isomorphism.

Let
$A$ denote the normalization of
the fiber product 
$Y_{(y)}\times_{X_{(x)}}
X'_{(x')}=
Y_{(y)}\times_SS'$.
Since $A$ is finite over the strictly
local scheme $X'_{(x')}$
and since the degree of
$Y_{(y)}$ over $X_{(x)}$
is the same as the degree of
$A$ over $X'_{(x')}$,
it suffices to show that
$A\to X'_{(x')}$ is an isomorphism.

Since $t(\eta)$
is a point of $V$
and $V'=V\times_SS'$,
the inverse image of
$t(\eta)\subset Y_{(y)}$
by the finite morphism
$A\to Y_{(y)}$ of normal schemes
equals
$t'(\eta')\subset A$
for the generic point $\eta'\in S'$.
Hence by the going down theorem,
the inverse image
of the closed point $y\in Y_{(y)}$
by $A\to Y_{(y)}$ consists
of the unique point $y'$.
Since $A$ is a disjoint union
of finitely many strictly local
schemes and $Y'_{(y')}$
is one of its component, we have
$A=Y'_{(y')}$.

By the assumption that
$Y'\to X'$ is \'etale
on a neighborhood of $y'\in t'(S')$,
the morphism
$A=Y'_{(y')}\to X'_{(x')}$
of strictly local schemes
is \'etale and hence
is an isomorphism.
Thus the assertion follows.
\qed
}

\subsection{Construction of functors}\label{ssFun}

We briefly recall the construction
on which the definition of
the filtration by ramification groups
is based.
The definition itself is recalled in
Section \ref{ssRam}.
Let $K$ be a henselian discrete valuation field
with residue field $F$. 
Let $L$ be a finite separable extension
of $K$ with residue field $E$.
Let $Q$ be a smooth scheme over
$S={\rm Spec}\, {\mathcal O}_K$
and let $T={\rm Spec}\, {\mathcal O}_L\to Q$
be a closed immersion.

Let $r>0$ be a rational number.
Let $K'$ be a finite extension of
$K$ of ramification index $e$
with residue field $F'$
and set $S'={\rm Spec}\, {\mathcal O}_{K'}$.
Assume that $er>0$ is an integer.
Let $Q_{S'}^{[er]}=
Q_{S'}^{[erQ_{F'}\cdot T_{S'}]}$ denote the 
dilatation of the base change
$Q_{S'}=Q\times_SS'$
at $T_{S'}=T\times_SS'$
and the pull-back $erQ_{F'}$
of the Cartier divisor
${\rm Spec}\, {\mathcal O}_{K'}/{\mathfrak m}^{er}_{K'}
\subset S'$.
Since $S'\to S$ is flat,
the immersion 
$T_{S'}\to Q_{S'}$ is a regular immersion
and the conormal sheaf
$N_{T_{S'}/Q_{S'}}$ is the pull-back of
$N_{T/S}$.
Hence
the reduced closed fiber
$(Q_{S'}^{[er]}\times_S{\rm Spec}\, F')
_{\rm red}$
denoted abusively by $Q^{[er]}_{F'}$
is canonically identified with
the vector bundle over
$E'=(E\otimes_FF')_{\rm red}$
associated to the dual of the
free $E'$-module
${\rm Hom}_{E'}
({\mathfrak m}_{K'}
^{er}/{\mathfrak m}_{K'}^{er+1}
\otimes_{F'}E'
,N_{T/Q}\otimes_{{\mathcal O}_L}E')$
by Proposition \ref{prdiltr}.2.
If we assume that
$E$ is purely inseparable
over $F$, then $E'$
is also a purely inseparable extension
of $F'$.
In practice, we may and will later
assume that this condition
be satisfied since we are interested
in ramification.

Let
$Q_{S'}^{(er)}\to Q_{S'}^{[er]}$
denote the normalization
and let
$\bar T_{S'}$ be 
the normalization of the fiber product
$T_{S'}=T\times_SS'$.
By the functoriality of dilatations
and normalizations,
the immersion $T\to Q$ defines
a commutative diagram
\begin{equation}
\begin{CD}
\bar T_{S'}@>>>
T_{S'}\\
@VVV@VVV\\
Q_{S'}^{(er)}
@>>>
Q_{S'}^{[er]}.
\end{CD}
\label{eqTTQQ}
\end{equation}

By the reduced fiber
theorem \cite{BLR},
there exists
a finite separable
extension $K'$ such that
the geometric closed fiber
$Q_{S'}^{(er)}\times_{S'}{\rm Spec}\, \bar F$
is reduced.
If $K'$ is an extension
of a Galois closure of $L$ over $K$,
the normalization
$\bar T_{S'}$
is a disjoint union of sections,
parametrized by
morphisms $L\to K'$ over $K$.
For such extensions $K'$,
the construction of the normalization
$Q_{S'}^{(er)}$ commutes with
base change.
Let $\bar K$ be a separable closure
and $\bar F$ be the residue field.
Then, the geometric closed fiber
\begin{equation}
Q_{\bar F}^{(r)}
=
Q_{S'}^{(er)}\times_{S'}{\rm Spec}\, \bar F
\label{eqQFbr}
\end{equation}
is independent of such
$K'\subset \bar K$.

Let 
\begin{equation}
Q_{\bar F}^{(r)}
\to
Q_{\bar F}^{[r]}
=(Q_{S'}^{[er]}\times_{S'}{\rm Spec}\, \bar F)_{\rm red}
\label{eqQFb[r}
\end{equation}
denote the canonical morphism
of the reduced closed fibers
and 
$T_{\bar F}=(T\times_S\, {\rm Spec}\, \bar F)_{\rm red}
={\rm Spec}\,
(E\otimes_F\bar F)_{\rm red}$,
by abuse of notation.
Since the formation of
$Q_{S'}^{[er]}$ commutes with
base change by 
Lemma \ref{lmdilfl}.1,
$Q_{\bar F}^{[r]}$ does not
depend on $S'$ either.

Let ${\mathfrak m}_{\bar K}^r
\supset {\mathfrak m}_{\bar K}^{r+}$
denote
$\{x\in {\bar K}\mid
{\rm ord}_Kx\geqq r\}
\supset \{x\in {\bar K}\mid
{\rm ord}_Kx> r\}$
and define a vector bundle
$N^{(r)}_{T/Q,{\bar F}}$ over $T_{\bar F}$
by 
\begin{equation}
N^{(r)}_{T/Q,{\bar F}}
=
{\rm Hom}_{E_{\bar F,{\rm red}}}
({\mathfrak m}_{\bar K}
^r/{\mathfrak m}_{\bar K}^{r+}
\otimes_{F'}E_{\bar F,{\rm red}}
,N_{T/Q}\otimes_{{\mathcal O}_L}E_{\bar F,{\rm red}})^\vee.
\label{eqNTr}
\end{equation}
The right hand side is defined as
the spectrum of
the symmetric algebra of
the $E_{\bar F,{\rm red}}$-module
with $\vee$ removed.
If we fix a morphism $i_0\colon L\to \bar K$,
we define a connected component
\begin{equation}
N^{(r)\circ}_{T/Q,{\bar F}}
=
{\rm Hom}_{E_{\bar F,{\rm red}}}
({\mathfrak m}_{\bar K}
^r/{\mathfrak m}_{\bar K}^{r+}
,N_{T/Q}\otimes_{{\mathcal O}_L}\bar F)^\vee
\subset
N^{(r)}_{T/Q,{\bar F}}
\label{eqNTrcirc}
\end{equation}
by the induced morphism $\bar i_0
\colon E\to \bar F$.

By Proposition \ref{prdiltr}.2,
we have a canonical isomorphism
\begin{equation}
Q_{\bar F}^{[r]}
\to 
N^{(r)}_{T/Q,{\bar F}}.
\label{eqQNT}
\end{equation}
If we fix a morphism $i_0\colon L\to \bar K$,
in the commutative diagram
(\ref{eqTTQQ}),
the morphism $i_0$
defines a section
$S'\to \bar T_{S'}$.
This defines $\bar F$-points in 
$Q_{\bar F}^{(r)}
\to
Q_{\bar F}^{[r]}$.
Let
\begin{equation}
Q_{\bar F}^{(r)\circ}
\to
Q_{\bar F}^{[r]\circ}
\label{eqcirc}
\end{equation}
denote the connected components
containing the points.
The isomorphism
(\ref{eqQNT})
induces an isomorphism
$Q_{\bar F}^{[r]\circ}
\to 
N^{(r)\circ}_{T/Q,{\bar F}}$.
In the case where
$E$ is purely inseparable over $F$,
we have
$Q_{\bar F}^{[r]\circ}=
Q_{\bar F}^{[r]}$,
or equivalently
$N^{(r)\circ}_{T/Q,{\bar F}}
=N^{(r)}_{T/Q,{\bar F}}$.

We give an explicit description
in monogenic case.

\begin{lm}\label{lmmono}
Let $L$ be a wildly ramified
finite Galois extension of $K$
of degree $n$ of Galois group $G$
and extend the normalized valuation
${\rm ord}_K$ to $L$.
Assume that ${\mathcal O}_L$
is generated by one element
$\alpha$ over ${\mathcal O}_K$.
Let $f\in {\mathcal O}_K[X]$
be the minimal polynomial of $\alpha$
and set $f=\prod_{i=1}^n
(X-\alpha_i)$ so that $\alpha_n=\alpha$
and
${\rm ord}_K(\alpha_i-\alpha_n)$
is increasing.
Define an ideal $I=(\alpha_{n-1}-\alpha_n)
\subsetneqq {\mathcal O}_L$
and let $H=
{\rm Ker}(G\to
{\rm Aut}({\mathcal O}_L/I))
\subset G$
be the smallest non-trivial lower ramification group.
Let $E$ be the residue field of $L$.

{\rm 1}. 
The mapping
$H\to I\otimes_{{\mathcal O}_L}E$
sending $\sigma\in H$
to $\sigma(\alpha)-\alpha$
is an injection of abelian groups.

{\rm 2}.
Define a mapping
$\beta\colon H\to E$
by $\beta(\sigma)=
(\sigma(\alpha)-\alpha)/
(\alpha_{n-1}-\alpha_n)$.
Then,
the polynomial
$b=
\prod_{\sigma\in H}
(Y-\beta(\sigma))\in E[Y]$
is an additive polynomial
and 
we have an exact sequence
$$
\begin{CD}
0
@>>> H
@>\beta>> {\mathbf A}^1_E
@>b>> {\mathbf A}^1_E
@>>>0.\end{CD}$$

{\rm 3}.
Set $\delta
=\prod_{i=1}^{n-1}
(\alpha_n-\alpha_i)=f'(\alpha)$,
$r={\rm ord}_K\, \delta+
{\rm ord}_K(\alpha_{n-1}-\alpha_n)
\in \dfrac1n{\mathbf Z}$,
$m={\rm Card}\, (G\sm H)$
and define
\begin{equation}
c=
\prod_{i=1}^m
(\alpha_n-\alpha_i)
\cdot
(\alpha_{n-1}-\alpha_n)^{{\rm Card}\, H}.
\label{eqpir}
\end{equation}
Then, we have ${\rm ord}_K\, c=r$.

{\rm 4}.
Let $T={\rm Spec}\, {\mathcal O}_L
\to Q={\mathbf A}^1_S
={\rm Spec}\, {\mathcal O}_K[X]$
be the immersion defined by $\alpha$.
Define an isomorphism
${\mathbf A}^1_E\to
Q^{[r]}=N_{T/Q,E}^{(r)}$
by the basis $$\dfrac fc\in
N_{T/Q,E}^{(r)}
={\rm Hom}({\mathfrak m}_L^r/
{\mathfrak m}_L^{r+1},
N_{T/Q}\otimes_{{\mathcal O}_L}E)$$
sending $c$ to $f$.
Then, there exist a
connected component
$Q^{(r)\circ}_E\subset Q^{(r)}_E$,
an isomorphism
${\mathbf A}^1_E\to
Q^{(r)\circ}_E$
and a commutative diagram
\begin{equation}
\begin{CD}
{\mathbf A}^1_E@>>> Q^{(r)\circ}_E\\
@VbVV@VVV\\
{\mathbf A}^1_E@>{f/c}>> Q^{[r]}_E.
\end{CD}
\label{eqmono}
\end{equation}
\end{lm}


\proof{
1.
Since $H$ acts trivially
on $I\otimes_{{\mathcal O}_L}E$,
for $\sigma,\tau\in H$,
we have
$\sigma\tau(\alpha)-\alpha
-((\sigma(\alpha)-\alpha)
+(\tau(\alpha)-\alpha))
=
\sigma(\tau(\alpha)-\alpha)
-(\tau(\alpha)-\alpha)=0$
in $I\otimes_{{\mathcal O}_L}E$.

2.
Since $\alpha_{n-1}-\alpha_n$
is a basis 
of $I\otimes_{{\mathcal O}_L}E$,
the injection $\beta\colon H\to E$
is a morphism of abelian groups by 1.
Hence the assertion follows from
Lemma \ref{lmapol}.

3.
For $i=m+1,\ldots,n-1$,
we have
${\rm ord}_K(\alpha_i-\alpha_n)
={\rm ord}_K(\alpha_{n-1}-\alpha_n)$.

4.
Define a morphism
${\mathbf A}^1_T={\rm Spec}\, 
{\mathcal O}_L[Y]
\to Q={\rm Spec}\, 
{\mathcal O}_L[X]$
by sending $X$ to
$(\alpha_{n-1}-\alpha_n)Y+\alpha_n$.
Then since 
$f((\alpha_{n-1}-\alpha_n)Y+\alpha_n)
=
\prod_{i=1}^n
((\alpha_{n-1}-\alpha_n)Y-(\alpha_i-\alpha_n))$
is divisible by $c$,
by the universalities of dilatation
and normalization,
the morphism
${\mathbf A}^1_T
\to Q$ is uniquely lifted to an morphism
${\mathbf A}^1_T
\to Q^{(r)}$.

Since the morphism
${\mathbf A}^1_T
\to Q^{(r)}$ of normal schemes
is an isomorphism on the generic fiber
and the composition
${\mathbf A}^1_T
\to Q^{[r]}$ is quasi-finite,
the morphism
${\mathbf A}^1_T
\to Q^{(r)}$ is an open immersion.
Since $f/c \equiv b$,
the diagram (\ref{eqmono})
with $\circ$ omitted is commutative.
Since the morphism
$b\colon {\mathbf A}^1_E
\to {\mathbf A}^1_E$ is finite,
the open immersion
${\mathbf A}^1_E
\to Q^{(r)}_E$
is an isomorphism to a connected component.
\qed

}
\medskip


Since the normalization
$\bar T_{S'}$
is a disjoint union of copies of $S'$
indexed by morphisms $L\to K'$ over $K$,
if $K'$ contains the conjugates of $L$,
the closed fiber
$\bar T_{\bar F}$
is canonically identified with
the finite set ${\rm Mor}_K(L,\bar K)$.
Similarly,
$T_{\bar F}$
is canonically identified with
the finite set ${\rm Mor}_F(E,\bar F)$.
The reduced geometric closed fibers
of (\ref{eqTTQQ}) define
the commutative diagram
\begin{equation}
\begin{CD}
{\rm Mor}_K(L,\bar K)=\,&
\bar T_{\bar F}@>>>
T_{\bar F}
&\,={\rm Mor}_F(E,\bar F)
\\
&@VVV@VVV&\\
&Q_{\bar F}^{(r)}
@>>>
Q_{\bar F}^{[r]}.&
\end{CD}
\label{eqQF}
\end{equation}
The construction of (\ref{eqQF}) satisfies the
following functoriality.

\begin{pr}\label{prQTf}
Let ${\mathcal O}_K\to {\mathcal O}_{K'}$
be a morphism of henselian
discrete valuation rings
of ramification index $e$.
Let $L$ and $L'$ be finite separable 
extensions
of $K$ and $K'$ and
let $L\to L'$ be a morphism
compatible with $K\to K'$.
Fix a morphism $\bar K\to \bar K'$
of separable closures
compatible with $K\to K'$.
Let 
\begin{equation}
\begin{CD}
Q'@<<<T'&={\rm Spec}\, {\mathcal O}_{L'}\\
@VVV@VVV\\
Q@<<<T&={\rm Spec}\, {\mathcal O}_L
\end{CD}
\label{eqQQTT}
\end{equation}
be a commutative diagram of schemes
such that $Q$ and $Q'$
are smooth schemes over
$S={\rm Spec}\, {\mathcal O}_K$
and
$S'={\rm Spec}\, {\mathcal O}_{K'}$,
that the horizontal arrows are closed immersions and that
the vertical arrows are compatible with
$S'\to S$.
Let $r>1$ and $r'>1$
be rational numbers.

{\rm 1.}
Let $r'=er$.
Then, the diagram {\rm (\ref{eqQQTT})} 
defines a commutative diagram
\begin{equation}
\begin{CD}
\bar T'_{\bar F'}
@>>>
Q^{\prime(r')}_{\bar F'}
@>>>
Q^{\prime[r']}_{\bar F'}
@>>>
T'_{\bar F'}
\\
@VVV@VVV@VVV@VVV\\
\bar T_{\bar F}
@>>>
Q^{(r)}_{\bar F}
@>>>
Q^{[r]}_{\bar F}
@>>>
T_{\bar F}
\end{CD}
\label{eqQTf}
\end{equation}
of schemes compatible with
the induced morphism
$\bar F\to \bar F'$ of residue fields
of $\bar K$ and $\bar K'$.

{\rm 2.}
Assume that $S'=S$, $T'=T$, $r=r'$
and that $Q'\to Q$ is smooth.
Then 
the middle square in {\rm (\ref{eqQTf})} is cartesian
and the middle right vertical arrow
is a surjection
of vector bundles over $T_{\bar F}$.

{\rm 3.}
Assume that
the discrete valuation rings
${\mathcal O}_K$
and ${\mathcal O}_{K'}$
are essentially of finite type
over a discrete valuation ring
${\mathcal O}_{K_0}$ with perfect residue field $k$
and that 
there exist integers $m\geqq 0$
and $q\geqq 2$ satisfying the condition
{\rm (2)} in Proposition {\rm \ref{prImq}}
such that
$r\geqq \dfrac me\dfrac q{q-1}$
and $r'=er-m
\geqq \dfrac m{q-1}$.
Assume further that {\rm (\ref{eqQQTT})} is 
induced by an essentially cartesian diagram
\begin{equation}
\begin{CD}
Q'_0@<<< T'\\
@VVV@VVV\\
Q_0@<<<T
\end{CD}
\label{eqcar0b}
\end{equation}
of essential immersions.
Then,
the diagram {\rm (\ref{eqQQTT})} 
defines a commutative diagram
{\rm (\ref{eqQTf})} compatible with
the induced morphism
$\bar F\to \bar F'$ of residue fields
of $\bar K$ and $\bar K'$
as in {\rm 1}.
\end{pr}

\proof{
1.
Let 
$$\begin{CD}
S'_1@>>> S'\\
@VVV@VVV\\
S_1@>>>S
\end{CD}$$
be a commutative diagram of
spectra of henselian discrete valuation
rings such that the horizontal
arrows are finite of ramification 
indices $e_1$ and $e'_1$.
Assume that $e_1r$
and $e'_1r'$ are integers.
Then, by the naive functoriality 
of dilatations,
the morphism $Q'\to Q$
is lifted to a morphism
$Q^{\prime [e'_1r']}_{S'_1}
\to
Q^{[e_1r]}_{S_1}$.
This induces a morphism
$Q^{\prime (e'_1r')}_{S'_1}
\to
Q^{(e_1r)}_{S_1}$
of normalizations
and the assertion follows.

2.
By Lemma \ref{lmQT}.3,
if $Q'\to Q$ is smooth,
we may assume that
there is an \'etale morphism
$Q'\to Q'_0={\mathbf A}^n_Q$
such that the composition
$T\to Q'\to {\mathbf A}^n_Q$ is the $0$-section.
Then, in the notation of the proof of 1.,
we have an isomorphism
$Q^{\prime (e_1r)}_{S_1}
\to {\mathbf A}^n_{
Q^{(e_1r)}_{S_1}}\times_{Q'_0}Q'$ 
inducing isomorphisms
$Q^{\prime(r)}_{\bar F}
\to
Q^{\prime(r)}_{0,\bar F}$
and
$Q^{\prime[r]}_{\bar F}
\to
Q^{\prime[r]}_{0,\bar F}$
by Lemma \ref{lmdilfl}.2.
Hence the assertion follows.

3.
Let ${\mathcal I}_T\subset {\mathcal O}_Q$
and ${\mathcal I}_{T'}\subset {\mathcal O}_{Q'}$
be the ideals defining
the closed subschemes
$T\subset Q$ and $T'\subset Q'$.
Then, by Proposition \ref{prImq},
we have ${\mathcal I}_T{\mathcal O}_{Q'}
\subset {\mathfrak m}_{K'}^m{\mathcal I}_{T'}
+{\mathcal I}_{T'}^q$.
Hence by the generalized functoriality 
Lemma \ref{lmdilf} of dilatations,
the assertion follows as in 1.
\qed

}

\medskip

Since the left vertical arrow in (\ref{eqQF})
is induced by the immersion $T\to Q$,
its composition with
(\ref{eqQNT})
is the $0$-section
$T_{\bar F}\to N_{T/Q,\bar F}^{(r)}$
of the vector bundle.
The surjectivity of the mappings
$\varphi_{\bar x}^{D+}\colon
\bar Y_{\bar x}\to
F^{D+}_{\bar x}(Y/X)
\to
F^{D}_{\bar x}(Y/X)$
in \cite[Proposition 3.1.2 (1)]{TJM}
means exactly that the mappings 
\begin{equation}
\bar T_{\bar F}
\to 
Q_{\bar F}^{(r)}
\times_{
Q_{\bar F}^{[r]}}
T_{\bar F}
\to
\pi_0(Q_{\bar F}^{(r)})
\label{eqgd}
\end{equation}
defined by the diagram (\ref{eqQF}) 
are surjections.
This fact is a consequence
of the going down theorem.

\begin{cor}\label{corFr}
There exist surjections
\begin{equation}
\bar T_{\bar F}
\to F^{r+}(L)
\to F^{r}(L)
\label{eqcorFr}
\end{equation}
of finite sets satisfying the following property:
Let $T\to Q$
be an immersion to a smooth scheme over $S$.
Then, there exists a commutative diagram
\begin{equation}
\xymatrix{
\bar T_{\bar F}\ar[r]\ar[rd]
&
Q^{ (r)}_{\bar F}
\times_{Q^{ [r]}_{\bar F}}
T_{\bar F}
\ar[r]\ar[d]
&
\pi_0(Q^{ (r)}_{\bar F})
\ar[d]
\\
& F^{r+}(L)
\ar[r]
&F^{r}(L)
}
\label{eqcorFr2}
\end{equation}
such that the vertical arrows are bijections.
\end{cor}

\proof{
Let $T\to Q$ and $T\to Q'$
be immersions to smooth schemes
over $S$.
Then by Proposition \ref{prQTf}.2,
the smooth projections
$Q\times_SQ'\to Q$
and 
$Q\times_SQ'\to Q'$
induces bijections of the finite sets
compatible with the surjections
from $\bar T_{\bar F}$.
Hence the assertion follows.
\qed

}

\medskip

By Proposition \ref{prQTf}.1,
the constructions of
$F^r(L)$
and $F^{r+}(L)$
are functorial in $L$.

\begin{df}\label{dfFr}
Let
\begin{equation}
F\colon 
\text{\rm
(Finite separable extensions of $K$})
\to
\text {\rm
(Finite $G_K$-sets)}
\label{eqF}
\end{equation}
denote the fiber functor
defined by 
$F(L)=
\bar T_{\bar F}$.
Let $r>0$ be a rational number.

{\rm 1.}
We define functors
\begin{equation}
F^r, F^{r+}\colon 
\text{\rm
(Finite separable extensions of $K$})
\to
\text {\rm
(Finite $G_K$-sets)}
\label{eqFr}
\end{equation}
sending $L$ to
$F^r(L)$
and 
$F^{r+}(L)$ 
and surjective morphisms
\begin{equation}
F\to F^{r+}\to F^r
\label{eqFr+}
\end{equation}
of functors
as in
Corollary {\rm \ref{corFr}}.

{\rm 2.}
Let $L$ be a finite separable extension 
of $K$.
We say that the ramification is bounded by
$r$ (resp.\ by $r+$) if the mapping
$F(L)\to F^r(L)$
(resp.\ $F(L)\to F^{r+}(L)$)
is a bijection.
\end{df}
\medskip

We identify $F^1(L)=T_{\bar F}$
by \cite[Proposition 3.7]{AS},\,\cite[Proposition 3.3.5]{TJM}.
If $L$ is a finite Galois extension,
the Galois group $G={\rm Gal}(L/K)$
acts on the finite $G_K$-set
$F(L)={\rm Mor}_K(L,\bar K)$
and its quotients
$F^r(L)$ and $F^{r+}(L)$
by the functoriality.
This action commutes with
the $G_K$-action.
If $M$ is a subextension of $L$
and if we fix $L\to \bar K$,
the $G_K$-action
factors through the quotient
$G$.

The functor $F$ has the following
functoriality.
Let ${\mathcal O}_K\to {\mathcal O}_{K'}$
be a morphism of henselian discrete
valuation field
of ramification index $e$
and fix a morphism $\bar K\to 
\bar K'$ of separable closures.
Let $F'$ denote the functor $F$
for finite separable extensions of $K'$
and let $L\to L'$ be a morphism
of finite separable extensions
of $K$ and of $K'$.
Then, for rational numbers
$r>0$ and $r'=er$ or $r'=er-m>0$
as in Proposition \ref{prQTf}, 
the commutative diagram
(\ref{eqQTf})
defines a commutative
diagram
\begin{equation}
\begin{CD}
F'(L')@>>>F^{\prime r'+}(L')@>>>F^{\prime r'}(L')\\
@VVV@VVV@VVV\\
F(L)@>>>F^{r+}(L)
@>>>F^r(L).
\end{CD}
\label{eqFF'}
\end{equation}

\begin{lm}\label{lmRed0}
Let ${\mathcal O}_K\to {\mathcal O}_{K'}$
be a morphism of henselian discrete valuation
rings of ramification index $e\geqq 1$.
Let $r>0$ and $r'>0$
be rational numbers as in Proposition {\rm \ref{prQTf}}
satisfying one of the following conditions
{\rm (1)} and {\rm (2):}

{\rm (1)}
$r'=er$.

{\rm (2)}
The discrete valuation rings
${\mathcal O}_K$
and ${\mathcal O}_{K'}$
are essentially of finite type
over a discrete valuation ring
${\mathcal O}_{K_0}$ with perfect residue field $k$.
There exist integers $m\geqq 0$
and $q\geqq 2$ satisfies the condition
{\rm (1)} in Proposition {\rm \ref{prImq}}
such that
$r\geqq \dfrac me\dfrac q{q-1}$
and $r'=er-m
\geqq \dfrac m{q-1}$.

Let $L\to L'$ be a morphism
of finite separable extensions
of $K$ and of $K'$.
If the ramification of
$L$ is bounded by $r$ (resp.\ $r+$),
then ramification of
$L'$ is bounded by $r'$ (resp.\ $r'+$).
\end{lm}

\proof{
In the commutative diagram
(\ref{eqFF'}),
if the ramification of
$L$ over $K$ is bounded by
$r$ (resp.\ by $r+$),
then the composition of
lower horizontal arrows
(resp. the lower left horizontal arrow)
is a bijection.
Since the left verticl
arrow is an injection, 
the composition of
upper horizontal arrows
(resp. the upper left horizontal arrow)
is an injection and 
the ramification of
$L'$ over $K'$ is bounded by
$r'$ (resp.\ by $r'+$).
\qed

}
\medskip

\begin{pr}\label{prbdd}
{\rm (\cite{AS}, \cite{TJM})}
Let $K$ be a henselian discrete valuation
field and $L$ be a finite
separable extension.
Let
\begin{equation}
\begin{CD}
Q@<<< T&\, ={\rm Spec}\, {\mathcal O}_L\\
@VVV@VVV\\
P@<<< S&\, ={\rm Spec}\, {\mathcal O}_K\\
\end{CD}
\label{eqbdd}
\end{equation}
be a cartesian diagram of schemes
over $S$ where
the horizontal arrows are closed immersions
to smooth schemes over $S$
and the vertical arrows are
quasi-finite and flat.
Let $r>0$ be a rational number
and let $K'$ be
a finite separable extension
of ramification index $e$
of $K$ such that
$er$ is an integer
and that 
the morphism 
$Q^{(er)}_{S'}\to 
S'={\rm Spec}\, {\mathcal O}_{K'}$
has reduced geometric fibers.

{\rm 1.}
The following conditions are equivalent:

{\rm (1)}
The ramification of $L$ over $K$
is bounded by $r+$.

{\rm (2)}
The ramification of $L$ over $K$
is bounded by $s$
for every rational number $s>r$.

\noindent
If $r>1$,
these conditions are further equivalent to
the following condition.

{\rm (3)}
The finite morphism
$Q^{(r)}_{\bar F}\to
Q^{[r]}_{\bar F}$
is \'etale.

{\rm (4)}
The quasi-finite morphism
$Q^{(er)}_{S'}\to
P^{(er)}_{S'}$
is \'etale
on neighborhood of the
closed fibers.

{\rm 2.}
The following conditions are equivalent:

{\rm (1)}
The ramification of $L$ over $K$
is bounded by $r$.

{\rm (2)}
The finite morphism
$Q^{(r)}_{\bar F}\to
Q^{[r]}_{\bar F}$
is a split \'etale covering.

{\rm (3)}
The finite morphism
$Q^{(r)}_{\bar F}\to
P^{(r)}_{\bar F}$
is a split \'etale covering.
\end{pr}

\proof{
1.
The implication
(1)$\Rightarrow$(2)
is \cite[Lemma 3.1.6]{TJM}.
The implication
(2)$\Rightarrow$(1)
is \cite[Theorem 3.2.6]{TJM}.

If (3) is satisfied,
the diagram (\ref{eqQF})
is cartesian and (1) is satisfied.
The implication (1)$\Rightarrow$(4)
is \cite[Proposition 3.2.10]{TJM}.
We show
(4)$\Rightarrow$(3).
For an immersion 
$T\to Q$ to a smooth scheme over $S$,
there exists
a cartesian diagram (\ref{eqbdd})
by Lemma \ref{lmQT}.1.
By the reduced fiber theorem
\cite{BLR},
there exists
a finite separable extension
$K'$ of ramification 
index $e$ such that
$er$ is an integer
and that 
the morphism 
$Q^{(er)}_{S'}\to 
S'={\rm Spec}\, {\mathcal O}_{K'}$
has reduced geometric fibers.
Since the diagram (\ref{eqbdd})
is cartesian,
the finite morphism
$Q^{[r]}_{\bar F}
\to
P^{[r]}_{\bar F}=
P^{(r)}_{\bar F}$ is \'etale.
Hence
by taking the geometric closed fiber
of $Q^{(er)}_{S'}\to
P^{(er)}_{S'}$,
we obtain (3) from (4).
The implication
(2)$\Rightarrow$(3)
is also proved in \cite[Theorem 7.2]{AS}.

2.
We show that (2) and (3) are equivalent
to each other.
As in the proof of 1
(4)$\Rightarrow$(3),
there exists a cartesian diagram
(\ref{eqbdd}).
Since $Q\to P$ is flat,
the finite morphism
$Q^{[r]}_{\bar F}\to
P^{[r]}_{\bar F}$
is a split \'etale covering.
Hence
the conditions (2) and (3)
are equivalent to each other.

We show that (1) and (3) are equivalent
to each other.
As in the proof of 1
(4)$\Rightarrow$(3),
we have a cartesian diagram
(\ref{eqbdd})
and a finite separable extension
$K'$ loc.~cit.
Since the assertion is \'etale local,
we may assume that the
morphism $Q\to P$ is finite.
Hence 
$Q^{(er)}_{S'}\to 
P^{(er)}_{S'}$ is a finite
morphism of degree $[L:K]$
and the conditions (1)
and (3) 
are equivalent to each other.
\qed

}

\medskip
For a rational number $r>0$,
let
\begin{equation}
\Theta^{(r)}_{L/K,\bar F}=
{\rm Hom}_{E_{\bar F,{\rm red}}}
({\mathfrak m}_{\bar K}
^r/{\mathfrak m}_{\bar K}^{r+}
\otimes_{\bar F}E_{\bar F,{\rm red}}
,{\rm Tor}_1^{{\mathcal O}_L}
(\Omega^1_{{\mathcal O}_L/{\mathcal O}_K},
E_{\bar F,{\rm red}}))^\vee
\label{eqThLK}
\end{equation}
denote the vector bundle 
over $T_{\bar F}$.
If we fix a morphism
$i_0\colon L\to \bar K$,
we define a connected component
\begin{equation}
\Theta^{(r)\circ}_{L/K,\bar F}=
{\rm Hom}_{\bar F}
({\mathfrak m}_{\bar K}
^r/{\mathfrak m}_{\bar K}^{r+}
,{\rm Tor}_1^{{\mathcal O}_L}
(\Omega^1_{{\mathcal O}_L/{\mathcal O}_K},
\bar F))^\vee
\subset
\Theta^{(r)}_{L/K,\bar F}
\label{eqThLK0}
\end{equation}
where ${\rm Tor}$ is defined
with respect to the morphism
${\mathcal O}_L\to \bar F$
induced by $i_0\colon L\to \bar K$.
If $E$ is a purely inseparable extension of
$F$, we have
$\Theta^{(r)\circ}_{L/K,\bar F}
=\Theta^{(r)}_{L/K,\bar F}.$

The canonical injection
${\rm Tor}_1^{{\mathcal O}_L}
(\Omega^1_{{\mathcal O}_L/{\mathcal O}_K},
E)\to N_{T/Q}\otimes_{{\mathcal O}_L}E$
(\ref{eqTor})
induces a surjection
\begin{equation}
N_{T/Q,\bar F}^{(r)}\to 
\Theta^{(r)}_{L/K,\bar F}
\label{eqQTh}
\end{equation}
of vector bundles on $T_{\bar F}.$
If $T\to Q$ is minimal,
by Lemma \ref{lmTorH}.2,
the morphism (\ref{eqQTh})
is an isomorphism.

\begin{lm}\label{lmXL}

Let $r>1$ and assume that
the ramification of $L$ over $K$
is bounded by $r+$.
There exists a finite \'etale morphism
\begin{equation}
\begin{CD}
X_{L/K}^{(r)}
@>>>
\Theta^{(r)}_{L/K,\bar F}
\end{CD}
\label{eqXrdf}
\end{equation}
endowed with a bijection of
finite sets $\bar T_{\bar F}
\to X^{(r)}_{L/K,\bar F}
\times_ 
{\Theta^{(r)}_{L/K,\bar F}}
T_{\bar F}$
satisfying the following property:
Let $T={\rm Spec}\, {\mathcal O}_L\to Q$
be an immersion
to smooth schemes $Q$ over 
$S={\rm Spec}\, {\mathcal O}_K$.
Then, there exists a commutative diagram
\begin{equation}
\xymatrix{
\bar T_{\bar F}\ar[r]\ar[rd]
&
Q^{(r)}_{\bar F}
\ar[r]\ar[d]
&
Q^{[r]}_{\bar F}
\ar[d]
\\
& X^{(r)}_{L/K,\bar F}
\ar[r]
&
\Theta^{(r)}_{L/K,\bar F}
}
\label{eqXL2}
\end{equation}
such that the square is cartesian
and that the right vertical arrow
is the composition of
$Q^{[r]}_{\bar F}
\to
N_{T/Q,\bar F}^{(r)}$
{\rm (\ref{eqQNT})}
and
$N_{T/Q,\bar F}^{(r)}
\to
\Theta^{(r)}_{L/K,\bar F}$
{\rm (\ref{eqQTh})}.
\end{lm}

\proof{
If $T\to Q$ is minimal,
the morphism
$N_{T/Q,\bar F}^{(r)}
\to
\Theta^{(r)}_{L/K,\bar F}$
(\ref{eqQTh})
is an isomorphism
and the condition requires that
the vertical arrow
$Q^{(r)}_{\bar F}
\to
X^{(r)}_{L/K,\bar F}$
in (\ref{eqXL2}) is an isomorphism.

Let $T\to Q_0$ be a minimal
immersion to a smooth scheme over $S$
and set $Q^{(r)}_{0,\bar F}
=
X^{(r)}_{L/K,\bar F}$.
Then by Lemma \ref{lmQT}.2 and
Lemma \ref{lmmin}.2,
after replacing $Q$ by an \'etale neighborhood
of $T$,
there is a smooth morphism
$Q\to Q_0$ compatible with the immersions
of $T$.
Hence we obtain a required
commutative diagram
(\ref{eqXL2})  with cartesian square.
\qed

}
\medskip

Since a morphism of connected
finite \'etale schemes
is uniquely determined on the 
induced mapping on the fibers
of a point,
the finite \'etale scheme
$X_{L/K,\bar F}^{(r)}
\to
\Theta^{(r)}_{L/K,\bar F}$
(\ref{eqXrdf}) is characterized
uniquely up to a unique isomorphism
by the property in Lemma \ref{lmXL}.
If we fix a morphism $i_0\colon L\to \bar K$,
similarly as
$Q_{\bar F}^{(r)\circ}
\to
Q_{\bar F}^{[r]\circ}$,
we define connected components
\begin{equation}
X_{L/K,\bar F}^{(r)\circ}
\to
\Theta^{(r)\circ}_{L/K,\bar F}.
\label{eqXcirc}
\end{equation}

Since the construction of
$X^{(r)}_{L/K,\bar F}$ is functorial in $L$,
by sending $L$ to
(\ref{eqXrdf}),
we obtain a functor
\begin{align}
&\text{(Finite separable extensions of $K$
of ramification bounded by $r+$)}
\label{eqXT}
\\
&\to 
\text{(Finite \'etale morphism
of smooth schemes over ${\bar F}$).}
\nonumber
\end{align}
We identify 
$F^r(L)=\pi_0(
X^{(r)}_{L/K,\bar F})$ and
$F(L)=F^{r+}(L)=X^{(r)}_{L/K,\bar F}
\times_{\Theta^{(r)}_{L/K,\bar F}}T_{\bar F}$.

Assume that $L$ is a Galois extension.
Then the Galois group $G$
has a natural action on
$
X^{(r)}_{L/K}
\to
\Theta_{L/K,\bar F}^{(r)}$
(\ref{eqXrdf}).
By the injection
${\rm Tor}_1^{{\mathcal O}_L}
(\Omega^1_{T/S},E)
\to H_1(L_{E/S})$ (\ref{eqTorL}),
the inertia subgroup
$I\subset G$
acts trivially on
${\rm Tor}_1^{{\mathcal O}_L}
(\Omega^1_{T/S},E)$
and hence 
on $\Theta_{L/K,\bar F}^{(r)}$.

Let ${\mathcal O}_K\to {\mathcal O}_{K'}$
be a morphism of henselian discrete valuation
rings of ramification index $e_{K'/K}=e\geqq 1$.
Let $r>1$ be a rational number
and set $r'=er$ or
$r'=er-m>1$ as in Lemma
\ref{lmRed0}.
Let $L\to L'$ be a morphism of
finite separable extension of 
$K$ and of $K'$.
Assume that the ramification
of $L$ over $K$ is
bounded by $r+$ and 
and that the ramification
of $L'$ over $K'$ is
by $r'+$.
Then, the commutative diagram
(\ref{eqQTf})
gives a commutative diagram
\begin{equation}
\begin{CD}
X^{(r')}_{L'/K',\bar F'}
@>>> \Theta^{(r')}_{L'/K',\bar F'}
\\
@VVV@VVV\\
X^{(r)}_{L/K,\bar F}
@>>> \Theta^{(r)}_{L/K,\bar F}.
\end{CD}
\label{eqXrf}
\end{equation}

\subsection{Construction of an automorphism}\label{sskey}

Let $Q$ be a scheme
and $T,T'\subset Q$
be closed subschemes.
For a point
$t\in Q$
and an integer $n\geqq 0$,
we say that
$T\equiv T'
\bmod {\mathfrak m}_t^n$,
if for the morphism
$Q_n={\rm Spec}\,
{\mathcal O}_{Q,t}/
{\mathfrak m}_t^n\to Q$
of schemes,
we have an equality
$T\times_QQ_n=
T'\times_QQ_n$
of closed subschemes of $Q_n$.

\begin{lm}\label{lmQ12n}
Let $K$ be a henselian discrete valuation
field and $L$ and $L'$ be
finite separable extensions.
Let $Q$ be a smooth scheme over
$S={\rm Spec}\, {\mathcal O}_K$
and $T={\rm Spec}\, {\mathcal O}_{L}
\to Q$ and
$T'={\rm Spec}\, {\mathcal O}_{L'}
\to Q$ be closed immersions.
Let $n\geqq 2$ be an integer
and assume 
$T\equiv T'
\bmod {\mathfrak m}_t^n$
for the closed point
$t\in T\subset Q$.
Let $Q_F=Q\times_S{\rm Spec}\, F
\subset Q$
be the closed fiber and let
$Q^{[n]}=Q^{[nQ_F\cdot T]}$ and
$Q^{[n]\prime}=Q^{[nQ_F\cdot T']}$ 
denote the dilations.

{\rm 1.}
The canonical isomorphism
{\rm (\ref{eqNTt1})}
induces isomorphisms
\begin{align}
N_{T/Q}\otimes_{{\mathcal O}_{T}}E
\to&\,
{\rm Ker}
({\mathfrak m}_{Q,t}/
{\mathfrak m}^2_{Q,t}
\to
{\mathfrak m}_{T,t}/
{\mathfrak m}^2_{T,t})
\nonumber
\\&\,
=
{\rm Ker}
({\mathfrak m}_{Q,t}/
{\mathfrak m}^2_{Q,t}
\to
{\mathfrak m}_{T',t}/
{\mathfrak m}^2_{T',t})
\gets
N_{T'/Q}\otimes_{{\mathcal O}_{T'}}E.
\label{eqNT12}
\end{align}

{\rm 2.}
We have
an equality of dilatations
$Q^{[n]}
=Q^{[n]\prime}$.

{\rm 3.}
Let $s\colon T_{\bar F}=
{\rm Spec}\,
(E\otimes_F\bar F)_{\rm red}\to
N_{T/Q,\bar F}^{(n)}$
denote the section of
the vector bundle
defined by the morphism
$T'\to Q^{[n]\prime}=Q^{[n]}$
lifting $T'\to Q$
and the isomorphism
$Q_{F}^{[n]}
\to
N_{T/Q,F}^{(n)}$
{\rm (\ref{eqQNT})}.
Let $+s\colon
N_{T/Q,F}^{(n)}\to
N_{T/Q,F}^{(n)}$
denote the translation by $s$
and let
\begin{equation}
u\colon
N_{T/Q,F}^{(n)}
\to
N_{T'/Q,F}^{(n)}.
\label{eqN12}
\end{equation}
be the isomorphism
induced by {\rm(\ref{eqNT12})}.
Then, the diagram
\begin{equation}
\begin{CD}
Q_{F}^{[n]}
@>{\rm (\ref{eqQNT})}>>
N_{T/Q,F}^{(n)}
\\
@|@VV{u\circ +s}V\\
Q_{F}^{[n]\prime}
@>{\rm (\ref{eqQNT})}>>
N_{T'/Q,F}^{(n)}
\end{CD}
\label{eqcd1}
\end{equation}
is commutative.
\end{lm}

\proof{1.
The arrows in (\ref{eqNT12})
are the isomorphisms in Lemma \ref{lmNQt}.1.
By the assumption that
${\rm Spec}\, {\mathcal O}_{L}/
{\mathfrak m}^n_{L}=
{\rm Spec}\, {\mathcal O}_{L'}/
{\mathfrak m}^n_{L'}$ 
and $n\geqq2$,
we have the equality
in (\ref{eqNT12}).

2.
Also by
${\rm Spec}\, {\mathcal O}_{L}/
{\mathfrak m}^n_{L}=
{\rm Spec}\, {\mathcal O}_{L'}/
{\mathfrak m}^n_{L'}$ 
and $n\geqq2$,
we have
$Q^{[n]}=Q^{[n]\prime}$ by
Proposition \ref{prdiltr}.1.

3.
This follows from
Corollary \ref{cordiltr}.
\qed

}

\begin{lm}\label{lmQ12}
Let $K$ be a henselian discrete valuation
field and $L$ and $L'$ be
finite separable extensions.
Let $Q$ be a smooth scheme over
$S={\rm Spec}\, {\mathcal O}_K$
and $T={\rm Spec}\, {\mathcal O}_{L}
\to Q$ and
$T'={\rm Spec}\, {\mathcal O}_{L'}
\to Q$ be closed immersions.
Let $n\geqq 2$ be an integer
such that 
$T\equiv T'
\bmod {\mathfrak m}_t^n$
for the closed point
$t\in T\subset Q$.
Let
$Q^{[n]}=Q^{[nQ_F\cdot T]}$ and
$Q^{[n]\prime}=Q^{[nQ_F\cdot T']}$ 
be the dilations 
and $u\colon
N_{T/Q,F}^{(n)}
\to
N_{T'/Q,F}^{(n)}$
be the isomorphism
{\rm (\ref{eqN12})}
as in Lemma {\rm \ref{lmQ12n}}.
Let $f\colon T'\to T$
be a morphism over $S$
inducing the identity on the residue field
$E$.

Then, there exist an \'etale neighborhood
$Q'\to Q$ of $T'$
and a morphism $\widetilde f\colon Q'\to Q$
lifting $f\colon T'\to T$
such that under the identification
$Q^{\prime[n]}_{F}=
Q^{\prime[nQ'_F\cdot T']}=
Q^{[n]\prime}_{F}$
induced by the \'etale
morphism $Q'\to Q$,
the diagram
\begin{equation}
\begin{CD}
Q_{F}^{[n]\prime}
@>{\rm (\ref{eqQNT})}>>
N_{T'/Q,F}^{(n)}
\\
@V{\widetilde f}VV@VV{u^{-1}}V\\
Q_{F}^{[n]}
@>{\rm (\ref{eqQNT})}>>
N_{T/Q,F}^{(n)}
\end{CD}
\label{eqcd2}
\end{equation}
is commutative.

The commutative diagrams
{\rm (\ref{eqcd1})}
and 
{\rm (\ref{eqcd2})}
define a commutative diagram
\begin{equation}
\begin{CD}
Q_{\bar F}^{(n)}
@>>>
N_{T/Q,\bar F}^{(n)}
\\
@VVV@VV{+s}V\\
Q_{\bar F}^{(n)}
@>>>
N_{T/Q,\bar F}^{(n)}.
\end{CD}
\label{eqcd3}
\end{equation}
Further if $f\colon T'\to T$
is an isomorphism,
the left vertical arrow is an isomorphism
and the diagram
{\rm (\ref{eqcd3})} is cartesian.
\end{lm}

\proof{
Let $D\subset Q$ be a regular effective
Cartier divisor meeting
$T$ transversally at the closed
point $t={\rm Spec}\, E$.
Since $T\equiv T'
\bmod 
{\mathfrak m}^2_t$,
the divisor $D$ meets $T'$
also transversally at $t$.
Since 
$f\colon T'\to T$
induces the identity on
$t={\rm Spec}\, E=
T\cap D\to T'\cap D$,
the morphism $f$
and the identity $1_D$
are extended to a morphism
$\widetilde f_1\colon T'\cup D
\to T\cup D$ of
closed subschemes of $Q$.
Since $Q$ is smooth over $S$,
there exists an \'etale neighborhood
$Q'\to Q$ of $t$ and
a morphism $\widetilde f\colon Q'\to Q$
lifting $\widetilde f_1$.

We show the commutativity of
the square (\ref{eqcd2}).
We identify
${\mathfrak m}_{Q,t}/
{\mathfrak m}_{Q,t}^2
=
{\mathfrak m}_{Q',t}/
{\mathfrak m}_{Q',t}^2$
and
${\mathfrak m}_{T,t}/
{\mathfrak m}_{T,t}^2=
{\mathfrak m}_{T',t}/
{\mathfrak m}_{T',t}^2$.
Then by (\ref{eqNT12}),
it suffices to show that
$\widetilde f\colon Q'\to Q$
induces the identity
on 
${\rm Ker}({\mathfrak m}_{Q,t}/
{\mathfrak m}_{Q,t}^2
\to
{\mathfrak m}_{T,t}/
{\mathfrak m}_{T,t}^2)$.
This follows from the isomorphism
(\ref{eqNTt3})
and the fact that
$\widetilde f_1$ induces
the identity on $D$
and hence on $N_{t/D}$.

Combining
{\rm (\ref{eqcd1})}
and 
{\rm (\ref{eqcd2})}
together with the functoriality of
{(\ref{eqQFb[r})},
we obtain a commutative diagram
$$\begin{CD}
Q_{\bar F}^{(n)}
@>{(\ref{eqQFb[r})}>>
Q_{\bar F}^{[n]}
@>{\rm (\ref{eqQNT})}>>
N_{T/Q,\bar F}^{(n)}
\\
@|@|@VV{u\circ +s}V\\
Q_{\bar F}^{\prime(n)}
@>{(\ref{eqQFb[r})}>>
Q_{\bar F}^{\prime[n]}
@>{\rm (\ref{eqQNT})}>>
N_{T'/Q,\bar F}^{(n)}
\\
@V{\widetilde f}VV@V{\widetilde f}VV
@VV{u^{-1}}V
\\
Q_{\bar F}^{(n)}
@>{(\ref{eqQFb[r})}>>
Q_{\bar F}^{[n]}
@>{\rm (\ref{eqQNT})}>>
N_{T/Q,\bar F}^{(n)}.
\end{CD}
$$
Thus, we obtain (\ref{eqcd3}).

By the assumption that
$T$ meets $D$ transversally
at $t$,
the canonical morphism
${\mathfrak m}_{Q,t}/
{\mathfrak m}_{Q,t}^2
\to
{\mathfrak m}_{T,t}/
{\mathfrak m}_{T,t}^2
\times
{\mathfrak m}_{D,t}/
{\mathfrak m}_{D,t}^2$
is an isomorphism.
Similarly,
${\mathfrak m}_{Q',t}/
{\mathfrak m}_{Q',t}^2
\to
{\mathfrak m}_{T',t}/
{\mathfrak m}_{T',t}^2
\times
{\mathfrak m}_{D,t}/
{\mathfrak m}_{D,t}^2$
is an isomorphism.
Hence if $f\colon T'\to T$ is an isomorphism,
then $\widetilde f\colon Q'\to Q$
is \'etale.
This implies that
the left vertical arrow of
{\rm (\ref{eqcd3})} is an isomorphism
and the diagram itself is cartesian.
\qed

}
\medskip

Let $f,g\colon S={\rm Spec}\, A\to X$
be morphisms of schemes
and $I\subset A$ be an ideal.
If the restrictions of $f$ and $g$
to the closed subscheme
${\rm Spec}\, A/I\subset S$
are the same,
we write $f\equiv g \bmod I$.

\begin{lm}\label{lmcong}
Let $K$ be a henselian
discrete valuation field
and let $L$ be a finite
separable extension
of $K$.
Let
\begin{equation}
\begin{CD}
Q@<<<T&={\rm Spec}\, {\mathcal O}_L\\
@VVV@VVV\\
P@<i<< S&={\rm Spec}\, {\mathcal O}_K
\end{CD}
\label{eqPQTSS}
\end{equation}
be a cartesian diagram
of separated schemes over $S$
such that the horizontal
arrows are closed immersions
to smooth schemes over $S$
and that the vertical
arrows are finite and flat.
Let $t\in T\subset Q$
be the closed point
and let $n\geqq 2$ be an integer.

{\rm 1.}
Let $i'\colon S\to P$ be
a closed immersion over $S$
such that $i'\equiv i
\bmod 
{\mathfrak m}_K^n$.
Let $T'=Q\times_PS$
be the fiber product 
with respect to $i'\colon S\to P$
and set $T'={\rm Spec}\, A$.

Then, we have
$T\equiv T'\bmod 
{\mathfrak m}_{t}^n$
and $A$ is a discrete valuation ring.

{\rm 2.}
Let $P^{[n]}=P^{[nP_F\cdot S]}$ denote
the dilatation and
let $S_1={\rm Spec}\, {\mathcal O}_{K_1}
\to P^{[n]}$
be the henselization
at the generic point $\xi$
of the closed fiber $P^{[n]}_F
=P^{[n]}\times_S
{\rm Spec}\, F$.
Let
$i_1,
i'_1\colon S_1\to P_1=P\times_SS_1$
be the sections
induced by the compositions
$S_1\to S\overset i\to P$
and by
$S_1\to P^{[n]}\to P$
respectively
and consider the cartesian diagram
\begin{equation}
\begin{CD}
T_1@>>> Q_1=Q\times_SS_1
@<<<T'_1\\
@VVV@VVV@VVV\\
S_1@>{i_1}>> P_1=P\times_SS_1
@<{i'_1}<<S_1.
\end{CD}
\label{eqTT1}
\end{equation}
Let $t_1\in T_1\subset Q_1$
denote the closed point.

Then, ${\mathcal O}_{K_1}$
is a discrete valuation ring
and the residue field
$F_1=k(\xi)$ is a purely transcendental
extension of $F$.
We have $i_1\equiv i'_1
\bmod {\mathfrak m}_{K_1}^n$
and $T_1\equiv T'_1
\bmod {\mathfrak m}_{t_1}^n$
and the schemes
$T_1$ and $T'_1$
are spectrums of discrete valuation rings.
\end{lm}

\proof{1.
We consider $S$ as
a closed subscheme
of $P$ by the section
$i$ and
let $s\in S\subset P$ be the closed point.
Let $S'\subset P$ denote
the closed subscheme defined
by the section $i'\colon S\to P$.
Then, since
$i\equiv i'\bmod 
{\mathfrak m}_K^n$,
we have
$S\equiv S'\bmod 
{\mathfrak m}_s^n$.
Hence
we have
$T\equiv T'\bmod 
{\mathfrak m}_{t}^n$.

By $T\equiv T'\bmod 
{\mathfrak m}_{t}^n$,
the maximal ideal of $A/{\mathfrak m}^n_A$
is generated by one element.
Since $n\geqq 2$,
by Nakayama's lemma,
${\mathfrak m}_A$ is generated by
one element.
Since $A$ is finite flat over ${\mathcal O}_K$,
we have
$\dim A=1$
and hence $A$ is a discrete valuation ring.

2.
By Lemma \ref{lmdilfl}.2,
$P^{[n]}$ is smooth over $S$
and the closed fiber
$P^{[n]}_F$ is a vector space over $F$.
Hence ${\mathcal O}_{K_1}$
is a discrete valuation ring
and the residue field
$F_1=k(\xi)$ is a purely transcendental
extension of $F$.
By the definition of dilatations,
the morphism
$P^{[n]}\to P$ is congruent to
the composition $P^{[n]}\to P\to S
\overset i\to P$ modulo
the pull-back of ${\mathfrak m}_K^n$.
Hence, we have
$i_1\equiv i'_1
\bmod {\mathfrak m}_{K_1}^n$.

The fiber product
$T_1=Q_1\times_{P_1}S_1$ 
is equal to
$T\times_SS_1$
with respect to the composition
$S_1\to P^{[n]}\to P\to S$.
Since $P^{[n]}\to S$ is smooth
and the residue field
$F_1=k(\xi)$ is a purely transcendental
extension of $F$,
we have
$T_1=T\times_SS_1={\rm Spec}\, {\mathcal O}_{L_1}$
for a discrete valuation ring ${\mathcal O}_{L_1}$.
By 1.\ applied to
$i_1\equiv i'_1
\bmod {\mathfrak m}_{K_1}^n$,
we have
$T_1\equiv T'_1
\bmod {\mathfrak m}_{t_1}^n$
and 
$T'_1={\rm Spec}\, {\mathcal O}_{L'_1}$
for a discrete valuation ring ${\mathcal O}_{L'_1}$.
\qed

}

\begin{lm}\label{lmtrs}
Let $K$ be a henselian
discrete valuation field
and let $L$ be a finite
separable extension
of $K$.
Let $n\geqq 2$ be an integer
and cartesian diagrams
{\rm (\ref{eqPQTSS})} and 
{\rm (\ref{eqTT1})} be as in Lemma
{\rm \ref{lmcong}}.
Assume that the ramification
of $L$ over $K$ is bounded by $n+$.

{\rm 1.}
There exists
an open neighborhood
$W\subset P^{[n]}$
of the closed fiber
$P^{[n]}_F\subset P^{[n]}$
such that on the generic fiber
$W_K =W\times_S{\rm Spec}\, K$,
the base change
$Q_{W_K}=
Q\times_PW_K
\to W_K$ is \'etale.

{\rm 2.}
There exists a finite separable extension
$L'_1$ of $K_1$
such that $T'_1={\rm Spec}\,
{\mathcal O}_{L'_1}$.

{\rm 3.}
Let $T_{W_K}=T\times_SW_K$
and define a cartesian diagram
\begin{equation}
\begin{CD}
A@<<<
A_K&\,=I\!som_{W_K}
(Q_{W_K},T_{W_K})\\
@VVV@VVV\\
W@<<<W_K
\end{CD}
\label{eqAW}
\end{equation}
by letting $A$ denote
the normalization of $W$
in $A_K$.
Let $A^{\text {\rm\'et}}\subset A$
denote the largest open subscheme
\'etale over $W$.

Then, the intersection
${\rm Im}(A^{\text {\rm\'et}}\to W)
\cap P^{[n]}_F$
is dense in the closed fiber $P^{[n]}_F$.
Further, the universal isomorphism
$Q\times_PA_K
\to
T\times_SA_K$
is extended to an isomorphism
$Q\times_PA^{\text {\rm\'et}}
\to
T\times_SA^{\text {\rm\'et}}$.
\end{lm}

\proof{
1.
By the reduced fiber theorem
\cite{BLR},
there exists
a finite separable extension
$K'$ of ramification 
index $e$ such that
the morphism 
$Q^{(en)}_{S'}\to 
S'={\rm Spec}\, {\mathcal O}_{K'}$
has reduced geometric fibers.
By the assumption that
the ramification
of $L$ over $K$ is bounded by $n+$
and by Proposition \ref{prbdd}
(1)$\Rightarrow$(4),
the finite morphism
$Q^{(en)}_{S'}
\to 
P^{[en]}_{S'}$
is \'etale on a neighborhood
$W'$ of the closed fiber
$P^{[en]}_{F'}$.
Since $S'\to S$
is faithfully flat,
the image $W
\subset P^{[n]}$
of $W'$ satisfies
the required condition.

2.
Let $i_1\colon S_1\to P_1=P\times_SS_1$
denote the section
induced by $S_1\to S\overset i\to P$
and
let $i'_1\colon S_1\to P_1$
denote the section
induced by $S_1\to P^{[n]}\to P$
as in Lemma
{\rm \ref{lmcong}}.
Then, by Lemma \ref{lmcong}.2,
we have $T'={\rm Spec}\, {\mathcal O}_{L'_1}$
for a finite morphism
${\mathcal O}_{K_1}\to {\mathcal O}_{L'_1}$
of discrete valuation rings.

Since $\xi\in P^{[n]}_F\subset W$,
the canonical morphism $S_1\to P^{[n]}$
factors through $S_1\to W$.
Since
${\rm Spec}\, L'_1=
Q_{W_K}\times_{W_K}{\rm Spec}\, K_1$,
the finite extension $L'_1$ of $K_1$
is separable.

3.
We apply Proposition
\ref{prnormal}
to $W\supset W_K\gets 
A_K$.
Similarly as in the proof of 2.,
the lifting $S\to P^{[n]}$
of $i\colon S\to P$ factors through
$S\to W$.
The cartesian diagram
{\rm (\ref{eqPQTSS})}
induces an isomorphism
$Q_{W_K}\times_{W_K}{\rm Spec}\, K
\to
T_{W_K}\times_{W_K}{\rm Spec}\, K$.
This define a section of
the finite \'etale scheme
$A_K\to W_K$
at the point ${\rm Spec}\, K\subset
S\to P^{[n]}$.
By the universality of the normalization,
the section 
${\rm Spec}\, K\to A_K$
is extended to a section
$S\to A$.

Let $S'$ be as in the proof of 1.
Changing $W'$ if necessary,
we may assume that
the finite morphism $Q^{(en)}_{S'}\to 
P^{[en]}_{S'}$ is \'etale on
$W'=W\times_SS'$.
We may further assume that
$K'$ contains a Galois closure of
$L$.
Then the normalization $\bar T_{S'}
\to S'$ of $T\times_SS'$
is a split finite \'etale covering.
Hence $A'=
I\!som_{W'}(Q^{(en)}_{S'}
\times_{P^{(en)}_{S'}}W',
\bar T_{S'}
\times_{S'}W')$
is a finite \'etale scheme
over $W'$.
Since $A'\to S'$
is the normalization of
$W'\subset P^{[n]}\times_SS'$
in  
$A_K \times_{W_K}W'_{K'}$,
it follows that
$A^{\text {\rm\'et}}$ is an open neighborhood
of the section $S\to W$
by Proposition \ref{prnormal}.
Since $P^{[n]}_F$ is irreducible,
the non-empty intersection
${\rm Im}(A^{\text {\rm\'et}}\to W)
\cap P^{[n]}_F$ is dense.

Since $A^{\text {\rm\'et}}\to W$
is \'etale,
$Q\times_PA^{\text {\rm\'et}}$
and 
$T\times_SA^{\text {\rm\'et}}$
are the normalizations of
$A^{\text {\rm\'et}}$
in 
$Q\times_PA_K$
and 
$T\times_SA_K$.
Hence the universal isomorphism
$Q\times_PA_K
\to
T\times_SA_K$
is extended to an isomorphism
$Q\times_PA^{\text {\rm\'et}}
\to
T\times_SA^{\text {\rm\'et}}$.
\qed

}
\medskip

We prepare a crucial construction
in the proof of additivity
Theorem \ref{thmlin}.

\begin{pr}\label{prisom}
Let $K$ be a henselian
discrete valuation field
and let $L$ be a finite
separable extension of $K$.
Let
\begin{equation}
\begin{CD}
Q@<<<T&={\rm Spec}\, {\mathcal O}_L\\
@VVV@VVV\\
P@<i<< S&={\rm Spec}\, {\mathcal O}_K
\end{CD}
\label{eqPQTSSp}
\end{equation}
be a cartesian diagram
of schemes over $S$
such that the horizontal
arrows are closed immersions
to smooth schemes over $S$
and that the vertical
arrows are finite and flat.
Let $n\geqq 2$ be an integer 
such that the ramification of
$L$ over $K$ is bounded by $n+$.

Assume that the residue
field $E$ of $L$ is a purely inseparable
extension of the residue field $F$ of $K$.
Let $F_1=k(\xi)$ be the residue field
at the generic point
$\xi\in P^{(n)}_F=
P^{(n)}\times_S{\rm Spec}\, F$ of
the closed fiber
and $\bar F\to \bar F_1$
be a morphism of algebraic closures.
Then, there exists
a cartesian diagram
\begin{equation}
\begin{CD}
Q_{\bar F}^{(n)}
\times_{\bar F}\bar F_1
@>>>
N_{T/Q,\bar F}^{(n)}
\times_{\bar F}\bar F_1
\\
@VVV@VV{+\xi}V\\
Q_{\bar F}^{(n)}
\times_{\bar F}\bar F_1
@>>>
N_{T/Q,\bar F}^{(n)}
\times_{\bar F}\bar F_1
\end{CD}
\label{eqcd4}
\end{equation}
where the horizontal arrows are
the base change of 
the composition of
{\rm (\ref{eqQFb[r})}
and {\rm (\ref{eqQNT})}.
\end{pr}

\proof{
We use the notation in
{\rm (\ref{eqTT1})} in Lemma
{\rm \ref{lmcong}}.
We identify $F_1$
with the residue field at the closed
point of $S_1$.
The horizontal arrow
$Q_{\bar F}^{(n)}
\times_{\bar F}\bar F_1
\to
N_{T/Q,\bar F}^{(n)}
\times_{\bar F}\bar F_1$
(\ref{eqcd4})
is canonically identified with
$Q_{1,\bar F_1}^{(n)}
\to
N_{T_1/Q_1,\bar F_1}^{(n)}$
defined for
$T_1\to Q_1$ over $S_1$ in (\ref{eqTT1}).

The finite extension
$L'_1$ is a separable extension of $K_1$
by Lemma \ref{lmtrs}.2.
Hence by Lemma \ref{lmQ12},
in order to define a cartesian
diagram (\ref{eqcd4}),
it suffices to define a finite
unramified extension
$K_2$ of $K_1$
and an isomorphism
${\mathcal O}_{L_2}\to {\mathcal O}_{L'_2}$
for $L_2=L_1K_2$ and $L'_2=L'_1K_2$
inducing the identity on
the residue fields
$E_2$ of ${\mathcal O}_{L_2}$.
To construct them,
we apply Lemma \ref{lmtrs}.3.

%

As in Lemma {\rm \ref{lmtrs}.3},
let $A$ be the normalization
of $W\subset P^{[n]}$
in $A_K=I\!som_{W_K}
(Q_{W_K},T_{W_K})$
over the generic fiber $W_K\subset W$
and let $A^{\text {\rm\'et}}\subset A$
denote the largest open subscheme
\'etale over $W$.
By  Lemma {\rm \ref{lmtrs}.3},
the intersection
${\rm Im}(A^{\text {\rm\'et}}\to W)
\cap P^{[n]}_F$
is dense in $P^{[n]}_F$.
Take a point
$\xi_2\in A^{\text {\rm\'et}}$ above
the generic point $\xi \in P^{[n]}_F$
and 
let $S_2={\rm Spec}\, {\mathcal O}_{K_2}$ be the
henselization of
$A$ at $\xi_2$.
Then $K_2$ is a finite unramified
extension of $K_1$.

Since $T_1=T\times_SS_1$
and $T'_1=Q\times_PS_1$,
by pulling back the universal isomorphism
$Q\times_PA^{\text {\rm\'et}}
\to T\times_SA^{\text {\rm\'et}}$
on $A^{\text {\rm\'et}}$
by $S_2\to A^{\text {\rm\'et}}$,
we obtain an isomorphism
$f\colon T'_2=T'_1\times_{S_1}S_2
\to
T_2=T\times_SS_2$.
Since $E$ is assumed to be
purely inseparable over $F$,
the residue field $E_2
=EF_2$ of $T_2$
is also purely inseparable over
the residue field $F_2$ of $K_2$
and hence the morphism
$f$ induces the identity on $E_2$.
Thus by applying Lemma \ref{lmQ12}
to the isomorphism $f$,
we obtain a cartesian diagram
(\ref{eqcd4}).
%
\qed

}

\section{Ramification groups}\label{sRam}


\subsection{Ramification groups}\label{ssRam}

\begin{thm}\label{thmFr}
{\rm (\cite{AS}, \cite{TJM})}
Let $K$ be a henselian discrete valuation
field. Let $L$ be a finite Galois extension 
of $K$ of Galois group $G$.

{\rm 1. (\cite[Theorems 3.3]{AS},\,
\cite[Theorem 3.3.1 (1)]{TJM})}
There exists a unique decreasing filtration
$(G^r)$ indexed by rational numbers $r>0$ by normal subgroups
of $G$ indexed by positive rational
numbers such that,
for every rational number
$r>0$ and for every intermediate
extension $K\subset M\subset L$,
the morphism $F\to F^r$ of
functors induces a bijection
\begin{equation}
F(M)/G^r\to F^r(M)
\label{eqFMr}
\end{equation}
of finite $G_K$-sets.

{\rm 2. (\cite[Theorems 3.8]{AS},\,
\cite[Theorem 3.3.1 (2)]{TJM})}
There exists a finite increasing sequence
$0=r_0<r_1<\cdots<r_n$ of
rational numbers 
such that $G^r$ is constant in the
semi-open intervals $(r_{i-1},r_i]$
for $i=1,\ldots, n$ and
in $(r_n,\infty)$.

{\rm 3. (\cite[Theorem 3.3.1 (1)]{TJM})}
For a rational number $r>0$, 
set $G^{r+}=\bigcup_{s>r}G^s$.
Then, 
for every rational number
$r>0$ and for every intermediate
extension $K\subset M\subset L$,
the morphism $F\to F^{r+}$ of
functors induces a bijection
\begin{equation}
F(M)/G^{r+}\to F^{r+}(M)
\label{eqFMr+}
\end{equation}
of finite $G_K$-sets.

{\rm 4. (\cite[Proposition 3.7]{AS},\,\cite[Proposition 3.3.5]{TJM})}
The subgroups
$G^1\supset G^{1+}$
equal the inertia subgroup
$I$ and its $p$-Sylow subgroup $P$.
\end{thm}

\begin{cor}\label{cortors}
Let $L$ be a finite Galois extension
of $K$.
Let $r>1$ be a rational number
and assume that the ramification of
$L$ over $K$ is bounded by $r+$.
Let $T={\rm Spec}\, {\mathcal O}_L
\to Q$ be an immersion
to a smooth scheme over 
$S={\rm Spec}\, {\mathcal O}_K$.
Fix a morphism $i_0\colon L\to \bar K$
and let
$Q^{(r)\circ }_{\bar F}\subset
Q^{(r)}_{\bar F}$ 
and 
$Q^{[r]\circ }_{\bar F}\subset
Q^{[r]}_{\bar F}$ {\rm (\ref{eqcirc})}
denote the connected
components.

{\rm 1.}
$Q^{(r)\circ }_{\bar F}$ 
is a $G^r$-torsor over
$Q^{[r]\circ}_{\bar F}$ 
and 
$X_{L/K,\bar F}^{(r)\circ}$
{\rm (\ref{eqXcirc})}
is a $G^r$-torsor over
$\Theta^{(r)\circ}_{L/K,\bar F}$.

{\rm 2.}
The $G^r$-torsor
$Q^{(r)\circ }_{\bar F}$ 
over $Q^{[r]\circ }_{\bar F}$ 
is additive if and only if
the $G^r$-torsor
$X_{L/K,\bar F}^{(r)\circ}$
over
$\Theta^{(r)\circ}_{L/K,\bar F}$
is additive.
\end{cor}

\proof{
1.
Since $r>1$, the subgroup
$G^r\subset G^{1+}=P$
acts trivially on $\Theta^{(r)}_{L/K,\bar F}$
and acts on
$X_{L/K,\bar F}^{(r)\circ}$
as an automorphism
over $\Theta^{(r)\circ}_{L/K,\bar F}$. 
Since the fiber 
$F^{r+}(L)\cap X_{L/K,\bar F}^{(r)\circ}
\subset F^{r+}(L)=F(L)$
of 
$X_{L/K,\bar F}^{(r)\circ}
\to \Theta^{(r)\circ}_{L/K,\bar F}$
at $0$ 
is a $G^r$-torsor,
the finite \'etale
scheme $X_{L/K,\bar F}^{(r)\circ}$
over $\Theta^{(r)\circ}_{L/K,\bar F}$
is a $G^r$-torsor.
By the cartesian square
in (\ref{eqXL2}),
$Q^{(r)\circ }_{\bar F}$ 
is also a $G^r$-torsor over
$Q^{[r]\circ}_{\bar F}$.

2.
This follows from
the cartesian square
in (\ref{eqXL2})
and Corollary \ref{corlinb}.3.
\qed

}
\medskip

The filtrations 
$(G^r)_{r>0}$ and
$(G^{r+})_{r>0}$
satisfy the following
functoriality.
Let ${\mathcal O}_K\to {\mathcal O}_{K'}$
be a morphism of henselian discrete
valuation field
of ramification index $e$.
Let $L\to L'$ be a
morphism of finite Galois extensions
of $K$ and of $K'$.
Then, for rational numbers
$r>0$ and $r'=er$ or $r'=er-m>0$
as in Proposition \ref{prQTf}, 
the commutative diagram
(\ref{eqFF'})
implies that the morphism
$G'={\rm Gal}(L'/K')
\to
G={\rm Gal}(L/K)$
induces mappings
$G'/G^{\prime r'+}\to G/G^{r+}$
and
$G'/G^{\prime r'}\to G/G^r$
and hence morphisms
\begin{equation}
\begin{CD}
G^{\prime r'+}
\to
G^{r+},\qquad
G^{\prime r'}
\to
G^r.
\end{CD}
\label{eqGG'}
\end{equation}

We define the total dimension
of a Galois representation.
Let $C$ be a field of characteristic
different from $p$
and let $V$ be a representation 
of $G$ on a $C$-vector space
of finite dimension.
Then, since $G^{1+}$
is a $p$-group,
there exists a unique decomposition
\begin{equation}
V=\bigoplus_{r\in {\mathbf Q},r\geqq1}
V^{(r)}
\label{eqVr}
\end{equation}
such that
for every rational number
$r>1$ and for
the $G^r$-fixed part,
we have
\begin{equation}
V^{G^r}=
\bigoplus_{1\leqq s<r}V^{(s)}.
\label{eqVGr}
\end{equation}
The decomposition
(\ref{eqVr})
is called the slope decomposition.
The total dimension is defined by
\begin{equation}
{\rm dimtot}(V)=\sum_{r\in {\mathbf Q},r\geqq1}
r\cdot \dim V^{(r)}
\label{eqdtV}
\end{equation}
as a rational number $\geqq 
\dim V$.
The equality is equivalent to
the condition that
the action of $P=G^{1+}$
on $V$ is trivial.
In the classical case where
the residue field is perfect,
the total dimension
is the sum of the dimension
and the Swan conductor:
${\rm dimtot}(V)=\dim V+
{\rm Sw}\, V$
and is known to be an integer
by
\cite[Th\'eor\`eme 1$'$, Section 2, Chapitre VI]{CL}.

\subsection{Reduction steps}\label{ssdvr}

We recall some facts on
discrete valuation rings.

\begin{lm}[{\cite[Lemme 2.2.1]{Artin}}]\label{lmGal}
Let $K$ be a henselian discrete valuation
field and $\hat K$ be the completion.
Then, the completion defines an
equivalence of categories
\begin{equation}
(\text {\rm Finite separable extensions of } K)
\to 
(\text{\rm Finite separable extensions of } \hat K).
\end{equation}
For a finite separable extension
$L$ over $K$,
the canonical morphism
\begin{equation}
{\mathcal O}_L
\otimes_{{\mathcal O}_K}
{\mathcal O}_{\hat K}
\to
{\mathcal O}_{\hat L}
\label{eqOLcomp}
\end{equation}
is an isomorphism.
\end{lm}

\begin{pr}\label{prdvr}
Let $K$ be a complete discrete valuation
field with residue field $F$
of characteristic $p>0$.
Let $k\subset F$ be a perfect subfield.
If $K$ is of characteristic $0$,
let ${\mathcal O}_{K_0}
=W(k)\subset {\mathcal O}_K$ be the Witt ring.
If $K$ is of characteristic $p>0$,
we canonically identify $k$ with a subfield of $K$ and
let ${\mathcal O}_{K_0}
=k[[t]]\subset {\mathcal O}_K$ for a uniformizer $t$
of $K$.

Then, there exist an inductive system
$({\mathcal O}_{K_\lambda})
_{\lambda\in \Lambda}$
of henselian discrete valuation rings
essentially of finite type and
generically separable over ${\mathcal O}_{K_0}$
and an isomorphism
$\varinjlim_{\lambda\in \Lambda}
{\mathcal O}_{K_\lambda}
\to {\mathcal O}_K$
over ${\mathcal O}_{K_0}$.
\end{pr}

\proof{
We follow the argument in
\cite[(0.5)]{SGA7}.
First, we prove the case where
$K$ is of characteristic $0$.
Let $C$ be the Cohen ring with
residue field $F$ \cite[Chapitre 0, D\'efinition (19.8.4)]{EGA4}.
Since the residue field $F=C/{\mathfrak m}_{K_0}C$
is a separable extension of 
a perfect field $k$
and ${\rm Frac}\, C$ is a separable extension of $K_0$,
there exist an inductive system
$(A_\lambda)
_{\lambda\in \Lambda}$
of henselian discrete valuation rings
essentially smooth over ${\mathcal O}_{K_0}$
and an isomorphism
$\varinjlim_{\lambda\in \Lambda}
A_\lambda\to C$
over ${\mathcal O}_{K_0}$
by N\'eron's desingularization
\cite[Theorem (4.5)]{appr}.

There exists an injection
$C\to {\mathcal O}_K$
lifting the identity of $F$ by
\cite[Chapitre 0, Th\'eor\`eme (19.8.6)]{EGA4}.
Since ${\mathcal O}_K$ is of finite
presentation over $C$,
the assertion follows in this case.

Next, we prove the case where
$K$ is of characteristic $p>0$.
Then, since $\Omega^1_{K_0/k}$
is generated by $dt$
and since $dt\neq 0$ in $\Omega^1_{K/k}$,
the canonical morphism
$\Omega^1_{K_0/k}\otimes_{K_0}K\to 
\Omega^1_{K/k}$ is an injection and hence
$K$ is a separable extension of $K_0$
by \cite[Chapter V \S 16, No.\ 4, Corollary of Theorem 3]{NB}.
Hence, 
the assertion follows
by N\'eron's desingularization.
\qed}

\begin{lm}\label{lmOL}
Let $K\subset K'$ be an extension of
henselian discrete valuation fields.
Let $L$ be a finite separable extension of
$K$ and set $L'=L\otimes_KK'$.
Assume that the canonical morphism
${\mathcal O}_L
\otimes_{{\mathcal O}_K}
{\mathcal O}_{K'}
\to
{\mathcal O}_{L'}$
is an isomorphism.
Then, for any intermediate extension
$K\subset M\subset L$
and $M'=M\otimes_KK'$,
the canonical morphism
${\mathcal O}_M
\otimes_{{\mathcal O}_K}
{\mathcal O}_{K'}
\to
{\mathcal O}_{M'}$
is an isomorphism.
\end{lm}

\proof{
Since the inclusion
${\mathcal O}_M\to {\mathcal O}_L$
is faithfully flat,
the induced morphism
${\mathcal O}_M
\otimes_{{\mathcal O}_K}
{\mathcal O}_{K'}\to 
{\mathcal O}_L
\otimes_{{\mathcal O}_K}
{\mathcal O}_{K'}$
is also faithfully flat.
Since ${\mathcal O}_L
\otimes_{{\mathcal O}_K}
{\mathcal O}_{K'}={\mathcal O}_{L'}$
is regular,
${\mathcal O}_M
\otimes_{{\mathcal O}_K}
{\mathcal O}_{K'}$
is also regular by \cite[Chapitre 0, Proposition (17.3.3) (i)]{EGA4}.
\qed}

\medskip

We prepare reduction steps
of the proof of the main results.

\begin{pr}\label{prRed1}
Let ${\mathcal O}_K\to {\mathcal O}_{K'}$
be a morphism of henselian discrete valuation
rings of ramification index $e\geqq 1$.
Let $L$ be a finite Galois extension of $K$
of Galois group $G$.
Let $L'=LK'$ be a composition field
and let $G'={\rm Gal}(L'/K')$ be the Galois group.
Let $r>1$ and $r'>1$
be rational numbers
and assume one of the following conditions
{\rm (1)}--{\rm (3):}

{\rm (1)}
${\mathcal O}_K\to {\mathcal O}_{K'}$
is dominant on the tangent spaces
and $r'=r$.

{\rm (2)}
The discrete valuation rings
${\mathcal O}_K$
and ${\mathcal O}_{K'}$
are essentially of finite type and
generically separable 
over a discrete valuation ring
${\mathcal O}_{K_0}$ with perfect residue field $k$.
Integers $m\geqq 0$,
$q\geqq 2$ and 
$d=\dim_F\Omega^1_F$
and $d'=\dim_{F'}\Omega^1_{F'}$
satisfy the condition
{\rm (2)} in Proposition {\rm \ref{prImq}}
and {\rm (2)} in Proposition {\rm \ref{prTorTS}.4}
and $r\geqq \dfrac me\dfrac q{q-1}$
and $r'=er-m
\geqq \dfrac m{q-1}$.

{\rm (3)}
The morphism
${\mathcal O}_L\otimes_{{\mathcal O}_K}
{\mathcal O}_{K'}\to {\mathcal O}_{L'}$
is an isomorphism
and $r'=er$.

\noindent
Fix compatible
morphisms $i_0\colon L\to \bar K$
and $i'_0\colon L'\to \bar K'$
to separable closures.
Let $X^{(r)\circ }_{L/K,\bar F}\subset
X^{(r)}_{L/K,\bar F}$ 
and 
$\Theta^{(r)\circ }_{L/K,\bar F}\subset
\Theta^{(r)}_{L/K,\bar F}$ 
denote the connected
components containing the images of
$i_0\in F(L)\to F^{r+}(L)
\subset X^{(r)}_{L/K,\bar F}
\to
\Theta^{(r)}_{L/K,\bar F}$
and similarly for $L'/K'$.
Assume that 
the ramification of
$L$ is bounded by $r+$.

{\rm 1.}
The injection $G'={\rm Gal}(L'/K')\to
G$ induces an isomorphism
$G^{\prime r'}\to G^r$.

{\rm 2.}
In the cases {\rm (2)} and {\rm (3)},
the $G^r$-torsor
$X_{L/K,\bar F}^{(r)\circ}$
over
$\Theta_{L/K,\bar F}^{(r)\circ}$
is additive
if and only if
the $G^{\prime r'}$-torsor
$X_{L'/K',\bar F'}^{(r')\circ}$ over
$\Theta_{L'/K',\bar F'}^{(r')\circ}$ is additive.
\end{pr}

\proof{
Let 
\begin{equation}
\begin{CD}
{\rm Spec}\, {\mathcal O}_{L'}=\,&T'@>>> Q'\\
&@VVV@VVV\\
{\rm Spec}\, {\mathcal O}_L=\,&T@>>>Q
\end{CD}
\label{eqTTQQp}
\end{equation}
be a commutative diagram
of schemes where
the horizontal arrows are
closed immersions to smooth schemes
over $S$ and over $S'$
and the vertical arrows are
compatible with $S'\to S$.
In the case (1),
we assume that $T\to Q$
and $T'\to Q'$ are minimal.
In the case (2),
we assume that
(\ref{eqTTQQp}) is induced by
an essentially cartesian diagram
\begin{equation}
\begin{CD}
T'@>>> Q'_0\\
@VVV@VVV\\
T@>>>Q_0
\end{CD}
\label{eqTTQQ0p}
\end{equation}
of divisorial essential immersions
to smooth separated schemes over $S_0=
{\rm Spec}\, {\mathcal O}_{K_0}$.
In the case (3),
we assume that
the diagram (\ref{eqTTQQp}) is cartesian.

By Lemma \ref{lmRed0},
the ramification of
$L'$ is bounded by $r^{\prime+}$
and the injection $G'\to G$
induces an injection $G^{\prime r'}\to G^r$.
Let
\begin{equation}
\begin{CD}
\bar T'_{\bar F'}@>>>
Q^{\prime (r')}_{\bar F'}
@>>> Q^{\prime [r']}_{\bar F'}
\\
@VVV@VVV@VVV\\
\bar T_{\bar F}@>>>
Q^{(r)}_{\bar F}
@>>> Q^{[r]}_{\bar F}
\end{CD}
\label{eqTTQQr}
\end{equation}
be the commutative diagram
induced by (\ref{eqTTQQp}).
Since $i'_0\colon
L'\to \bar K'$
and $i_0\colon L\to K$ are compatible,
we obtain a commutative diagram
\begin{equation}
\begin{CD}
Q^{\prime (r')\circ}_{\bar F'}
@>>> Q^{\prime [r']\circ}_{\bar F'}
\\
@VVV@VVV\\
Q^{(r)\circ}_{\bar F}
@>>> Q^{[r]\circ}_{\bar F}
\end{CD}
\label{eqTTQQr0}
\end{equation}
of the connected components.
The horizontal arrows are
a $G^r$-torsor and
a $G^{\prime r'}$-torsor respectively
and 
the left vertical arrow
compatible with
$G^{\prime r'}\to G^r$.

In the case (1),
since $T\to Q$
and $T'\to Q'$ are assumed minimal,
the right vertical arrow of
(\ref{eqTTQQr0}) is
identified with
$$\Theta^{(r')\circ}_{L'/K',\bar F'}
={\rm Spec}\,
S({\rm Tor}_1^{{\mathcal O}_{L'}}
(\Omega^1_{{\mathcal O}_{L'}/{\mathcal O}_{K'}},\bar F'))
\to \Theta^{(r)\circ}_{L/K,\bar F}
={\rm Spec}\,
S({\rm Tor}_1^{{\mathcal O}_L}
(\Omega^1_{{\mathcal O}_L/{\mathcal O}_K},\bar F)).$$
Since the morphism
$S'\to S$ is assumed to
be dominant on the tangent spaces
at the closed point $s'\in S'$,
the morphism
$S(H_1(L_{\bar F/S}))
\to S(H_1(L_{\bar F'/S'}))$
of symmetric algebras
is an injection.
By Lemma \ref{lmTorH}.1,
the morphism
$S({\rm Tor}_1^{{\mathcal O}_L}
(\Omega^1_{{\mathcal O}_L/{\mathcal O}_K},\bar F))
\to 
S({\rm Tor}_1^{{\mathcal O}_{L'}}
(\Omega^1_{{\mathcal O}_{L'}/{\mathcal O}_{K'}},\bar F'))$
is also an injection.
Hence, the right vertical arrow
of {\rm (\ref{eqTTQQr0})} is dominant.
By Corollary \ref{corlinb}.1
applied to the commutative diagram
(\ref{eqTTQQr0}),
the injection
$G^{\prime r'}\to G^r$
of finite $p$-groups
is an isomorphism.

In the case (2),
by Proposition \ref{prTorTS}.4
(2)$\Rightarrow$(1),
the $E'$-linear mapping
$N_{T/Q}\otimes_{{\mathcal O}_T} E'\to
{\mathfrak m}_{K'}^m
N_{T'/Q'}\otimes_{{\mathcal O}_{T'}} 
E'$ is an injection.
Hence the right vertical arrow
of {\rm (\ref{eqTTQQr0})}
induces a surjection
$Q^{\prime [r']\circ}_{\bar F'}
\to
Q^{[r]\circ}_{\bar F}
\times_{\bar F}\bar F'$.
Hence by Corollary \ref{corlinb}.3,
the injection
$G^{\prime r'}
\to G^r$
is an isomorphism.
The equivalence of
the additivities
also follows from
Corollary \ref{corlinb}.3 and
Corollary \ref{cortors}.

In the case (3),
since the diagram {\rm (\ref{eqTTQQr})} 
is cartesian,
the right vertical arrow
of {\rm (\ref{eqTTQQr0})}
induces an isomorphism
$Q^{\prime [r]\circ}_{\bar F'}
\to
Q^{[r]\circ}_{\bar F}
\times_{\bar F}\bar F'$.
Hence the assertion follows
from Corollary \ref{corlinb}.3 and
Corollary \ref{cortors}
as in the case (2). 
\qed

}

\begin{cor}\label{corRed1}
Let the assumption be the
same as in Proposition {\rm \ref{prRed1}}
except that
we do {\em not}
assume that 
the ramification of
$L$ is bounded by $r+$.
Then, the canonical injection
$G'\to G$ induces
isomorphisms
${\rm Gr}^{\prime r'}G\to 
{\rm Gr}^rG$
and
$G^{\prime r'}\to G^r$.
\end{cor}

\proof{
Let $M\subset L$
be the subextension
corresponding to the subgroup
$G^{r+}\subset G$
and $M'=MK'\subset L'$.
Since the construction 
of the ramification
groups is compatible
with quotient \cite[Corollary 1.4.3]{TJM},
the subquotient
${\rm Gr}^rG=G^r/G^{r+}$
is identified with
$H^r
\subset {\rm Gal}(M/K)=H=G/G^{r+}$.
In the case (3),
the morphism
${\mathcal O}_M\otimes_{{\mathcal O}_K}
{\mathcal O}_{K'}\to {\mathcal O}_{M'}$
is an isomorphism
by Lemma \ref{lmOL}.
Since the ramification
of $M$ over $K$ is bounded by $r+$,
Applying Proposition \ref{prRed1}.1
to $M$ and
to $H'={\rm Gal}(MK'/K')
=G'/G^{\prime r'+}
\subset 
H={\rm Gal}(M/K)=G/G^{r+}$,
we see that
${\rm Gr}^{r'}G'
\to {\rm Gr}^rG$ is an isomorphism.

Hence for every rational number 
$s\geqq r$
and $s'=es\geqq r'$ in the cases
(1) and (3) or
$s'=es-m\geqq r'$
in the case (2),
the morphism
${\rm Gr}^{s'}G'
\to {\rm Gr}^sG$ is an isomorphism.
Thus
$G^{\prime r'+}\to G^r$
is also an isomorphism.
\qed

}

\subsection{Graded quotients}\label{ssGr}

\begin{thm}\label{thmgr}
Let $K$ be
a henselian discrete valuation field
with residue field of characteristic
$p>0$ and let
$L$ be a finite Galois
extension of $K$
of Galois group $G$.

{\rm 1.}
Let $r>1$ be a rational number.
Then, the graded quotient
${\rm Gr}^rG=G^r/G^{r+}$
is an abelian group
and is an ${\mathbf F}_p$-vector space.

{\rm 2.}
Let $C$ be a field of
characteristic different from $p$
and let $V$ be a representation
of $G$ on a $C$-vector space
of finite dimension.
Then,
the total dimension ${\rm dimtot}(V)$
is an integer.
\end{thm}

\proof{
In the classical case
where the residue field is perfect,
assertion 1 is proved in
\cite[Corollaire 1 of
Proposition 7, Section 2, Chapitre IV]{CL}
and assertion 2
is proved in 
\cite[Th\'eor\`eme 1$'$, Section 2, Chapitre VI]{CL}.
Hence Theorem 
follows from the existence
of tangentially dominant
extension with perfect residue field
Proposition \ref{prperf}
and the reduction step
Corollary \ref{corRed1} (1).
\qed}
\medskip

For an abelian extension,
we obtain the following.

\begin{cor}\label{corab}
Let the notation be as in 
Theorem {\rm \ref{thmgr}}
and assume that
$G$ is abelian.
Let $r>0$ be a rational number.
If $G^r\supsetneqq G^{r+}$,
then $r$ is an integer.
\end{cor}

\proof{
It suffices to apply
Theorem \ref{thmgr}.2
to a character $\chi\colon G\to C^\times$
such that ${\rm Ker}\chi\supset G^{r+}$
and ${\rm Ker}\chi\not \supset G^r$.
\qed
}
\medskip



\begin{thm}\label{thmlin}
Let $K$ be a henselian discrete
valuation field
and let $L$ be a finite
Galois extension of $K$
of Galois group $G$.
Let $r>1$ be a rational number
and assume that
the ramification of $L$ over $K$
is bounded by $r+$.
Fix a morphism $i_0\colon L\to \bar K$
to a separable closure
and let $X^{(r)\circ }_{L/K,\bar F}\subset
X^{(r)}_{L/K,\bar F}$ 
and 
$\Theta^{(r)\circ }_{L/K,\bar F}\subset
\Theta^{(r)}_{L/K,\bar F}$ 
{\rm (\ref{eqXcirc})}
denote the connected
components.

{\rm 1.}
The $G^r$-torsor
$X^{(r)\circ }_{L/K,\bar F}$ 
over $\Theta^{(r)\circ}_{L/K,\bar F}$
is additive.

{\rm 2.}
The subgroup $G^r\subset G$
is an ${\mathbf F}_p$-vector space.
The class $[X^{(r)\circ }_{L/K,\bar F}]$ 
defines an injection
\begin{equation}
{\rm Hom}(G^r,{\mathbf F}_p)
\to
{\rm Hom}_{\bar F}(
{\mathfrak m}_{\bar K}^r/
{\mathfrak m}_{\bar K}^{r+},
{\rm Tor}_1^{{\mathcal O}_L}
(\Omega^1_{{\mathcal O}_L/{\mathcal O}_K},
\bar F))
\label{eqgrT}
\end{equation}
of abelian groups,
where ${\mathcal O}_L\to \bar F$
is induced by $i_0$.

{\rm 3.}
Let $G_K={\rm Gal}(\bar K/K)
\to G={\rm Gal}(L/K)$
be the surjection defined by
$i_0\colon L\to \bar K$.
Then, the injection
{\rm (\ref{eqgrT})} is compatible
with $G_K\to G$
with respect to the conjugate
action of $G$ on the source
and the natural action of $G_K$
on the target.
\end{thm}

In the proof of Theorem \ref{thmlin},
we do {\em not} use
Theorem \ref{thmgr} or
Propositon \ref{prperf}.

\proof{
1.
First, we reduce the assertion 
to the case where $r>1$ is an integer.
By Proposition \ref{prRed1}.2 (3) 
and
Lemma \ref{lmGal},
we may assume that $K$ is complete
by replacing $K$ by the completion.
Further by Proposition \ref{prRed1}.2 (3)
and
Proposition \ref{prdvr},
we may assume that 
there exists
a discrete valuation subring
${\mathcal O}_{K_0}$ with perfect residue
field $k={\mathcal O}_{K_0}/{\mathfrak m}_{K_0}$
such that
${\mathcal O}_K$ is essentially
of finite type
and generically separable
over ${\mathcal O}_{K_0}$.

By Proposition \ref{prbsc} and
Lemma \ref{lmcar0}.3,
there exists a cartesian diagram
\begin{equation}
\begin{CD}
Q_0@<<< T\\
@VVV@VVV\\
P_0@<<< S
\end{CD}
\label{eqP0Q0}
\end{equation}
of schemes over $S_0
={\rm Spec}\, {\mathcal O}_{K_0}$
such that the vertical arrows 
are finite flat
and the horizontal arrows
are divisorial essentially immersions
to smooth separated schemes.
Let $m>1$ be an integer
such that $n=mr\geqq 2$ is an integer
and that the $p$-power part $q>1$ of $m$
satisfies $r\geqq \dfrac q{q-1}$.
Note that $r>1=\displaystyle{\lim_{q\to \infty}}\dfrac q{q-1}$.
We define 
a morphism ${\mathcal O}_K\to {\mathcal O}_{K'}$
of discrete valuation rings
as in Example \ref{eg129}.
Hence by 
Proposition \ref{prRed1}.2 (2),
Proposition \ref{prTorTS}.4 and
Example \ref{eg135},
by replacing $K$ by $K'$,
we may assume that $r=n>1$ is an integer.
By replacing $K$ by an unramified
extension,
we may assume that the
residue field $E$ of $L$
is a purely unramified
extension of the residue field $F$ of $K$.

We assume that $r=n\geqq 2$ is an integer
and that the
residue field $E$ of $L$
is a purely inseparable
extension of the residue field $F$ of $K$.
By Lemma \ref{lmmin}.1
and Lemma \ref{lmQT}.1,
there exists a cartesian diagram
\begin{equation}
\begin{CD}
Q@<<< T\\
@VVV@VVV\\
P@<<< S
\end{CD}
\label{eqPQt}
\end{equation}
of schemes over $S$
such that the vertical arrows 
are finite flat,
the horizontal arrows
are closed immersions
to smooth schemes over $S$
and that $T\to Q$ is minimal.
Let $F'$ be the residue field
of the generic point $\xi$
of the vector space $\Theta_{L/K,\bar F}^{(r)}$.

By Proposition \ref{prisom},
we obtain a commutative diagram
\begin{equation}
\begin{CD}
X_{L/K,\bar F}^{(r)}\times_{\bar F}{\bar F'}
@>>>
X_{L/K,\bar F}^{(r)}\times_{\bar F}{\bar F'}
\\
@VVV@VVV\\
\Theta_{L/K,\bar F}^{(r)}
\times_{\bar F}{\bar F'}
@>{+\xi}>>
\Theta_{L/K,\bar F}^{(r)}
\times_{\bar F}{\bar F'}.
\end{CD}
\label{eqtrsm}
\end{equation}
By the assumption that $E$ is purely inseparable over $F$,
the scheme
$\Theta_{L/K,\bar F}^{(r)}
=
\Theta_{L/K,\bar F}^{(r)\circ}$
is connected.
By replacing $f\colon T'_2\to T_2$
in the construction of (\ref{eqtrsm})
in the proof of 
Proposition \ref{prisom}
by the composition of
the automorphism of $T_2$
induced by an element of $G$,
we may assume that
the automorphism
$X_{L/K,\bar F}^{(r)}\times_{\bar F}{\bar F'}$
(\ref{eqtrsm})
preserves the connected component
$X_{L/K,\bar F}^{(r)\circ}
\times_{\bar F}{\bar F'}$
defined by $i_0\colon L\to \bar K$.

Hence by Proposition \ref{prisog},
the scheme $X_{L/K,\bar F}^{(r)\circ}$ has a 
structure of group scheme over $\bar F$
such that the finite \'etale morphism
$X_{L/K,\bar F}^{(r)\circ}
\to \Theta^{(r)}_{L/K,\bar F}$
of connected schemes
is a morphism of group schemes.
Consequently, by Proposition \ref{prlin}
(2)$\Rightarrow$(1),
the $G^r$-torsor
$X_{L/K,\bar F}^{(r)\circ}
\to \Theta^{(r)}_{L/K,\bar F}$
is additive.

2.
The assertion follows from
1. and Proposition \ref{prlin}
(1)$\Rightarrow$(3).

3.
The actions of $G_K$
and $G$ on
$\Theta^{(r)}_{L/K,{\bar F}}$
and on
$X^{(r)}_{L/K,{\bar F}}$
defined by the action of
$G$ on $L$ and 
the action of $G_K$
on $\bar K$
commute to each other.
Let $\widetilde \tau\in G_K$
and let $\tau\in G$ be its image.
Then, the automorphisms
$\widetilde \tau^{-1} \circ \tau$
of $\Theta^{(r)}_{L/K,{\bar F}}$
and
of $X^{(r)}_{L/K,{\bar F}}$
are compatible and preserve
the connected components
of $\Theta^{(r)\circ }_{L/K,{\bar F}}$
and
of $X^{(r)\circ}_{L/K,{\bar F}}$.
For $\sigma \in G$,
we have a commutative diagram
$$
\begin{CD}
X^{(r)\circ}_{L/K,{\bar F}}
@>{\sigma}>>
X^{(r)\circ}_{L/K,{\bar F}}
\\
@V{\widetilde \tau^{-1} \circ \tau}VV
@VV{\widetilde \tau^{-1} \circ \tau}V\\
X^{(r)\circ}_{L/K,{\bar F}}
@>{\tau\sigma\tau^{-1}}>>
X^{(r)\circ}_{L/K,{\bar F}}.
\end{CD}$$
Since the action of
$\widetilde \tau^{-1} \circ \tau$
on $\Theta^{(r)\circ }_{L/K,{\bar F}}$
is compatible with the action of
$\bar \tau$ on
${\rm Hom}_{\bar F}(
{\mathfrak m}_{\bar K}^r/
{\mathfrak m}_{\bar K}^{r+},$
${\rm Tor}_1^{{\mathcal O}_L}
(\Omega^1_{{\mathcal O}_L/{\mathcal O}_K},
\bar F))$,
the latter is compatible
with the conjugate action of
$\tau$ of $G^r$.
\qed

}

\medskip

The injection 
(\ref{eqgrT}) has the following functoriality.
Let $K\to K'$ be an extension
of henselian discrete valuation fields
of ramification index $e$.
Let $L\to L'$ be a morphism
of finite Galois extensions of
$K$ and of $K'$ 
and $G'={\rm Gal}(L'/K')\to G$ be the canonical morphism.
Let $r'=er$ or
$r'=er-m$ as in Lemma
\ref{lmRed0}.
Then, by the commutative diagram
(\ref{eqXrf}), the diagram
\begin{equation}
\begin{CD}
{\rm Hom}(G^r,{\mathbf F}_p)
@>>>
{\rm Hom}_{\bar F}(
{\mathfrak m}_{\bar K}^r/
{\mathfrak m}_{\bar K}^{r+},
{\rm Tor}_1^{{\mathcal O}_L}
(\Omega^1_{{\mathcal O}_L/{\mathcal O}_K},
E))
\\
@VVV@VVV\\
{\rm Hom}(G^{\prime r'},{\mathbf F}_p)
@>>>
{\rm Hom}_{\bar F}(
{\mathfrak m}_{\bar K'}^{r'}/
{\mathfrak m}_{\bar K'}^{r'+},
{\rm Tor}_1^{{\mathcal O}_{L'}}
(\Omega^1_{{\mathcal O}_{L'}/{\mathcal O}_{K'}},
E'))
\end{CD}
\label{eqgrf}
\end{equation}
is commutative.

\begin{cor}\label{corgr}
Let $K$ be a henselian discrete
valuation field
and let $L$ be a finite
Galois extension of $K$
of Galois group $G$.
Let $r>1$ be a rational number
and
fix a morphism $i_0\colon L\to \bar K$
to a separable closure.

{\rm 1.}
The graded quotient
${\rm Gr}^rG$ is an
${\mathbf F}_p$-vector space
and the injection
{\rm (\ref{eqgrT})} defines
an injection
\begin{equation}
{\rm Hom}({\rm Gr}^rG,{\mathbf F}_p)
\to
{\rm Hom}_{\bar F}(
{\mathfrak m}_{\bar K}^r/
{\mathfrak m}_{\bar K}^{r+},
{\rm Tor}^{{\mathcal O}_L}_1
(\Omega^1_{{\mathcal O}_L/
{\mathcal O}_K},\bar F)).
\label{eqgrL}
\end{equation}

{\rm 2.}
Let $G_K={\rm Gal}(\bar K/K)
\to G={\rm Gal}(L/K)$
be the surjection defined by
$i_0\colon L\to \bar K$.
Then, the injection
{\rm (\ref{eqgrL})} is compatible
with $G_K\to G$
with respect to the conjugate
action of $G$ on the source
and the natural action of $G_K$
on the target.

{\rm 3.}
The action of the
wild inertia group $P=G^{1+}$
on ${\rm Gr}^rG$ is trivial.
If ${\rm Gr}^rG\neq 1$,
the prime-to-$p$ part of
the denominator of $r$
divides the ramification index
$e=e_{L/K}$.
\end{cor}

\proof{
1.
Let $M\subset L$
be the intermediate extension
corresponding to
the subgroup $G^{r+}\subset G$
and identify
$G'={\rm Gal}(M/K)=G/G^{r+}$.
Then since the filtrations
are compatible with quotients by
\cite[Corollary 1.4.3]{TJM},
the composition of the injection
(\ref{eqgrT}) for $M$
and for $G^{\prime r}=
{\rm Gr}^rG$
with the injection
${\rm Tor}_1^{{\mathcal O}_M}
(\Omega^1_{{\mathcal O}_M/{\mathcal O}_K},
\bar F)
\to 
{\rm Tor}_1^{{\mathcal O}_L}
(\Omega^1_{{\mathcal O}_L/{\mathcal O}_K},
\bar F)$
(\ref{eqML})
defines an injection
\begin{align*}
{\rm Hom}({\rm Gr}^rG,{\mathbf F}_p)
=
{\rm Hom}(G^{\prime r},{\mathbf F}_p)
&\, \to
{\rm Hom}_{\bar F}(
{\mathfrak m}_{\bar K}^r/
{\mathfrak m}_{\bar K}^{r+},
{\rm Tor}^{{\mathcal O}_M}_1
(\Omega^1_{{\mathcal O}_M/
{\mathcal O}_K},\bar F))\\
&\, \to
{\rm Hom}_{\bar F}(
{\mathfrak m}_{\bar K}^r/
{\mathfrak m}_{\bar K}^{r+},
{\rm Tor}^{{\mathcal O}_L}_1
(\Omega^1_{{\mathcal O}_L/
{\mathcal O}_K},\bar F)).
\end{align*}

2.
This follows from Theorem \ref{thmlin}.3.

3.
Since the wild inertia subgroup
$P_K\subset G_K$ acts trivially on
the source of
(\ref{eqgrL}),
the action of the
wild inertia group $P$
on ${\rm Gr}^rG$ is trivial
by the injection (\ref{eqgrL}) and 2.

Assume ${\rm Gr}^rG\neq 1$
and let $m$ be 
the prime-to-$p$ part of
the denominator of $r$.
Then, the action of the tame inertia
group $I_K/P_K$ on
the image of ${\rm Gr}^rG$
induces an faithful action of
the quotient $\mu_m$.
Hence by 2,
$\mu_m$ is a quotient of
the tame inertia $I/P$ 
and $m$ divides $e=e_{L/K}$.
\qed}
\medskip

By the functoriality (\ref{eqgrf}),
the injection (\ref{eqgrL})
satisfies a similar functoriality.
Namely,
for an extension
of henselian discrete valuation fields,
we have a commutative diagram
\begin{equation}
\begin{CD}
{\rm Hom}({\rm Gr}^rG,{\mathbf F}_p)
@>>>
{\rm Hom}_{\bar F}(
{\mathfrak m}_{\bar K}^r/
{\mathfrak m}_{\bar K}^{r+},
{\rm Tor}^{{\mathcal O}_L}_1
(\Omega^1_{{\mathcal O}_L/
{\mathcal O}_K},\bar F))
\\
@VVV@VVV\\
{\rm Hom}({\rm Gr}^{r'}G',{\mathbf F}_p)
@>>>
{\rm Hom}_{\bar F}(
{\mathfrak m}_{\bar K'}^{r'}/
{\mathfrak m}_{\bar K'}^{r'+},
{\rm Tor}^{{\mathcal O}_{L'}}_1
(\Omega^1_{{\mathcal O}_{L'}/
{\mathcal O}_{K'}},\bar F))
\end{CD}
\label{eqgrLf}
\end{equation}
for $r'=er$ or
$r'=er-m>1$ as in Lemma
\ref{lmRed0}.

The composition of (\ref{eqgrL})
with the injection
${\rm Tor}^{{\mathcal O}_L}_1
(\Omega^1_{{\mathcal O}_L/
{\mathcal O}_K},\bar F)
\to
H_1(L_{\bar F/{\mathcal O}_K})$
(\ref{eqTorL})
defines an injection
\begin{equation}
{\rm Hom}({\rm Gr}^rG,{\mathbf F}_p)
\to
{\rm Hom}_{\bar F}(
{\mathfrak m}_{\bar K}^r/
{\mathfrak m}_{\bar K}^{r+},
H_1(L_{\bar F/{\mathcal O}_K})).
\label{eqgrK}
\end{equation}
The morphism (\ref{eqgrK})
satisfies a similar functoriality
as (\ref{eqgrLf}).
Taking the limit on $L$,
we obtain an injection
\begin{equation}
{\rm Hom}({\rm Gr}^rG_K,{\mathbf F}_p)
\to
{\rm Hom}_{\bar F}(
{\mathfrak m}_{\bar K}^r/
{\mathfrak m}_{\bar K}^{r+},
H_1(L_{\bar F/{\mathcal O}_K}))
\label{eqgrLK}
\end{equation}
for the absolute Galois group
$G_K={\rm Gal}(\bar K/K)$.

If $S$ is essentially of finite type
and generically separable
over $S_0={\rm Spec}\, {\mathcal O}_{K_0}$,
the composition of (\ref{eqgrK})
with the morphism
$H_1(L_{E/S})
\longrightarrow 
\Omega^1_{S/S_0}
\otimes_{{\mathcal O}_K}E$
(\ref{eqH1LOm})
defines a morphism
\begin{equation}
{\rm Hom}({\rm Gr}^rG,{\mathbf F}_p)
\to
{\rm Hom}_{\bar F}(
{\mathfrak m}_{\bar K}^r/
{\mathfrak m}_{\bar K}^{r+},
\Omega^1_{S/S_0}
\otimes_{{\mathcal O}_S}\bar F).
\label{eqgrOm}
\end{equation}
This is an injection if
$e_{K/K_0}\neq 1$ by Proposition
\ref{prLOm}.3.
For $r'=er$ or
$r'=er-m>1$ as in Lemma
\ref{lmRed0},
we have a commutative diagram
\begin{equation}
\begin{CD}
{\rm Hom}({\rm Gr}^rG,{\mathbf F}_p)
@>>>
{\rm Hom}_{\bar F}(
{\mathfrak m}_{\bar K}^r/
{\mathfrak m}_{\bar K}^{r+},
\Omega^1_{S/S_0}
\otimes_{{\mathcal O}_S}\bar F)
\\
@VVV@VVV\\
{\rm Hom}({\rm Gr}^{r'}G',{\mathbf F}_p)
@>>>
{\rm Hom}_{\bar F}(
{\mathfrak m}_{\bar K'}^{r'}/
{\mathfrak m}_{\bar K'}^{r'+},
\Omega^1_{S'/S'_0}
\otimes_{{\mathcal O}_{S'}}\bar F').
\end{CD}
\label{eqgrOmf}
\end{equation}



\begin{thebibliography}{99}

\bibitem{AS}
A.\ Abbes, T.\ Saito,
{\em Ramification of local fields with imperfect residue fields},
Amer.\ J.\ of Math., 124.5 (2002), 879-920.


\bibitem{cln}
A.\ Abbes, T.\ Saito,
{\em Ramification and cleanliness}, Tohoku Math.\ J.\ (2) 63 (2011), no. 4, 775-853.

\bibitem{appr}
M.\ Artin,
{\em Algebraic approximation of structures over complete local rings,}
Publ.\ Math.\ IHES {\bf 36}, 23-58 (1969).

\bibitem{Artin}
M.\ Artin,
{\em Dimension cohomologique:
Premiers r\'esultats},
SGA 4 III, Expos\'e X,
43-63, Th\'eorie des Topos et Cohomologie Etale des Sch\'emas
Tome 3, Springer Lecture Notes in Math., 305.

\bibitem{Bo}
J.\ Borger,
{\em Conductors and the moduli of residual perfection},
{\tt https://arxiv.org/}
{\tt abs/math/0112305}

\bibitem{BLR}
S.\ Bosch, W.\ L\"utkebohmert, M.\ Raynaud,
{\em Formal and rigid geometry. IV. The Reduced Fibre Theorem},
Inventiones Mathematicae,
119, (1995) no.\ 2, 361-398


\bibitem{NB}
N.\ Bourbaki,
{\sc Algebra I},
Chapters 1-3,
1989, Springer Verlag,
Berlin Heidelberg.

\bibitem{SGA7}
P.\ Deligne,
{R\'esum\'e des
premiers expos\'es de
A. Grothendieck},
Expos\'e I, p.\ 1-24,
Seminaire de Geometrie Algebrique du Bois-Marie 1967-1969, (SGA 7 I) Springer Lecture notes in math.\ 1972.

\bibitem{GR}
O.\ Gabber, L.\ Ramero,
{\em Foundations for almost ring theory -- Release 7.5},
{\tt https://}
{\tt arxiv.org/abs/math/0409584}.

\bibitem{EGA4}
A.\ Grothendieck,
{\sc \'El\'ements de g\'eom\'etrie alg\'ebrique IV}, 
\'Etude locale des sch\'emas et
des morphismes de sch\'emas, 
Publ.\ Math.\ IHES 20, 24, 28, 32 (1964-67).

\bibitem{17}
A.\ Grothendieck, 
{\em Cat\'egories cofibr\'ees additives et Complexe cotangent relatif}, Lecture Notes in Mathematics 79, Springer-Verlag, Berlin, Heidelberg, New York 1968.


\bibitem{Ill}
L.\ Illusie,
{\em Complexe cotangent et d\'eformations I},
Springer Lecture Notes in Math., 239,
Springer-Verlag, Berlin, Heidelberg, New York 1971.

\bibitem{Kato}
K.\ Kato,
{\em  Swan conductors for characters of degree 
one in the imperfect residue field case},  
Contemporary Math. {\bf 83}
(1989), 101--131. 





\bibitem{III}
T.\ Saito,
{\em Ramification of local fields with imperfect residue fields III}, 
Mathematische Annalen, 352, Issue 3 (2012), 567-580.

\bibitem{nonlog}
-----, 
{\em Wild Ramification and the Cotangent Bundle}, 
Journal of Algebraic Geometry,
26 (2017), 399-473.



\bibitem{TJM}
-----, 
{\em Ramification groups of coverings and valuations},
Tunisian J.\ of Math.,
Vol.\ 1, No.\ 3, 373-426, 2019.

\bibitem{TATA}
-----, 
{\em
A characterization of ramification groups
of local fields with imperfect residue fields},
{\tt arXiv:2004.03770}



\bibitem{LEX}
-----, 
{\em Cotangent bundles
and micro-supports
in mixed characteristic case},
{\tt arxiv:2006. 00448}



\bibitem{CL}
J-P.~Serre,
{\sc Corps Locaux,}
Hermann, Paris, 1968.



\bibitem{X}
L.\ Xiao, 
{\em On ramification filtrations and $p$-adic differential equations, II: mixed characteristic case},
Compositio Math.,
{\bf 148}, 2 (2012), pp.\ 415-463.
%
\end{thebibliography}
\end{document}